\documentclass[oneside, a4paper,11pt,reqno]{amsart}
\newcommand{\cal}{\mathcal}

\usepackage[T1]{fontenc}
\usepackage{verbatim}
\usepackage{amsmath}
\usepackage{amsfonts}
\usepackage{color}
\usepackage{multicol}
\usepackage{graphicx}
\usepackage{amssymb}
\usepackage{grffile}
\usepackage{amsmath}
\usepackage{bbm}
\usepackage{mdframed}
\usepackage{xcolor}
\usepackage{float}

\usepackage{amsthm}\usepackage{epsfig}

\newtheorem{theo}{Theorem}
\newtheorem{prop}[theo]{Proposition}
\newtheorem{cor}[theo]{Corollary}

\newcommand {\pare}[1] {\left( {#1} \right)}
\newcommand {\cro}[1] {\left[ {#1} \right]}
\newcommand {\acc}[1] {\left\{ {#1} \right\}}
\newcommand {\nor}[1] { \left\| {#1} \right\|}

\newcommand {\bra}[1] { \left\langle {#1} \right\rangle_{\mu}}

\newcommand {\abs}[1] {\left\lvert {#1} \right\rvert}
\newcommand {\va}[1] {\left| {#1} \right|}

\def \E {\mathbb{E}}
 
\def \P {\mathbb{P}}
\def \R  {\mathbb{R}} 
\def \Z  {\mathbb{Z}} 
\def \C  {\mathbb{C}} 

\def \AA {{\mathcal A}}
\def \EE {{\mathcal E}}
\def \FF {{\mathcal F}}
\def \PP {{\mathcal P}}
\def  \RR {{\mathcal R}}
\def \VV {{\mathcal X}}

\def  \VV {{\mathcal V}}
\def \SS {{\mathcal S}}

\def \GG {{\mathcal G}}
\def \BB {{\mathcal B}}
\def \RR {{\mathcal R}}

\def \ind {\hbox{ 1\hskip -3pt I}}
\def \bx {\bar{x}}

\def \bV {\bar{\VV}}
\def \bP {\bar{P}}
\def \bw {\bar{w}}
\def \by {\bar{y}}

\def \baf {\bar{f}}
\def \bR {\bar{R}}

\def \bLL {\bar{L}}

\def \brV {\breve{\VV}}
\def \brx {\breve{x}}

\def \brf {\breve{f}}
\def \brR {\breve{R}}

\def \Id {\mbox{Id}}
\def \Tr {\mbox{Trace}}

\begin{document}

\title[Random Forests and Networks Analysis]{Random Forests and Networks Analysis}
\author{Luca Avena} 
\address{Leiden University.}
\email{l.avena@math.leidenuniv.nl}
\today

\author{Fabienne Castell} 
\address{
 Aix-Marseille Universit\'e, CNRS, Centrale Marseille. I2M UMR CNRS 7373. 39, rue Joliot Curie. 13 453 Marseille Cedex
13. France.}
\email{fabienne.castell@univ-amu.fr}

\author{Alexandre Gaudilli\`ere} 
\address{
 Aix-Marseille Universit\'e, CNRS, Centrale Marseille. I2M UMR CNRS 7373. 39, rue Joliot Curie. 13 453 Marseille Cedex
13. France.}
\email{alexandre.gaudilliere@math.cnrs.fr} 

\author{Clothilde M\'elot} 
\address{
 Aix-Marseille Universit\'e, CNRS, Centrale Marseille. I2M UMR CNRS 7373. 39, rue Joliot Curie. 13 453 Marseille Cedex
13. France.}
\email{clothilde.melot@univ-amu.fr}


\begin{abstract} D. Wilson~\cite{[Wi]} in the 1990's described a simple and efficient algorithm based on loop-erased random walks to sample uniform spanning trees and more generally weighted trees or forests spanning a given graph. This algorithm provides a powerful tool in analyzing structures on networks and along this line of thinking, in recent works~\cite{AG1,AG2,ACGM1,ACGM2} we focused on applications of spanning rooted forests on finite graphs. The resulting main conclusions are reviewed in this paper by collecting related theorems, algorithms, heuristics and numerical experiments. A first foundational part on determinantal structures and efficient sampling procedures is followed by four main
applications: 1) a random-walk-based notion of well-distributed points in a graph 2) how to describe metastable dynamics in finite settings by means of  Markov intertwining dualities 3) coarse graining schemes for networks and associated processes 4) wavelets-like pyramidal algorithms for graph signals.
\end{abstract}

\maketitle

\section{\large{Introduction: networks, trees and forests}}\label{intro}

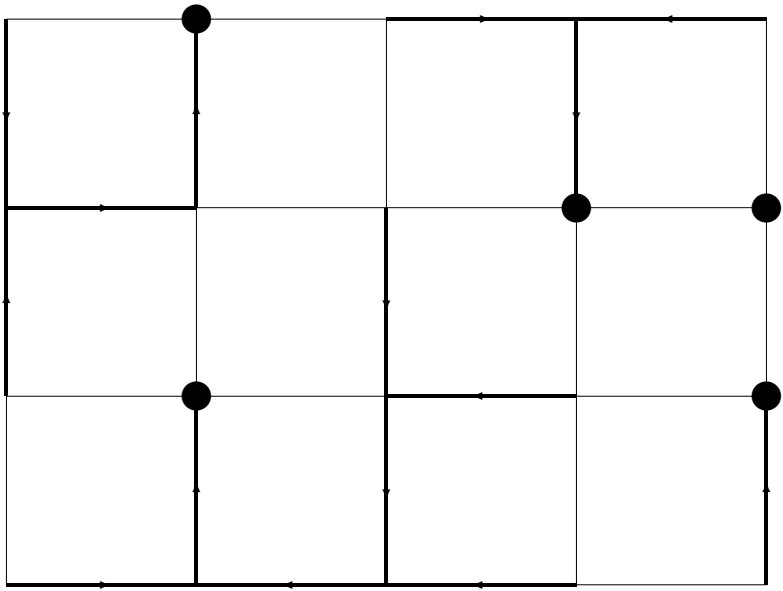
\begin{figure}
\begin{center}

\setlength{\unitlength}{2.5cm}
\begin{picture}(6,3)
\linethickness{0.01mm}
\multiput(1,0)(1,0){5}{\line(0,1){3}}

\linethickness{0.01mm}
\multiput(1,0)(0,1){4}
{\line(1,0){4}}

\put(2,1){\circle*{0.15}}
\put(2,3){\circle*{0.15}}
\put(4,2){\circle*{0.15}}
\put(5,2){\circle*{0.15}}
\put(5,1){\circle*{0.15}}

\linethickness{0.5mm}
\put(1,0){\vector(1,0){0.55}}
\put(1.55,0){\line(1,0){0.45}}

\linethickness{0.5mm}
\put(1,2){\vector(1,0){0.55}}
\put(1.55,2){\line(1,0){0.45}}

\linethickness{0.5mm}
\put(3,3){\vector(1,0){0.55}}
\put(3.55,3){\line(1,0){0.45}}


\linethickness{0.5mm}
\put(3,0){\vector(-1,0){0.55}}
\put(2.45,0){\line(-1,0){0.45}}

\linethickness{0.5mm}
\put(4,0){\vector(-1,0){0.55}}
\put(3.45,0){\line(-1,0){0.45}}

\linethickness{0.5mm}
\put(4,1){\vector(-1,0){0.55}}
\put(3.45,1){\line(-1,0){0.45}}

\linethickness{0.5mm}
\put(5,3){\vector(-1,0){0.55}}
\put(4.45,3){\line(-1,0){0.45}}

\linethickness{0.5mm}
\put(1,1){\vector(0,1){0.55}}
\put(1,1.55){\line(0,1){0.45}}

\linethickness{0.5mm}
\put(2,0){\vector(0,1){0.55}}
\put(2,0.55){\line(0,1){0.45}}

\linethickness{0.5mm}
\put(5,0){\vector(0,1){0.55}}
\put(5,0.55){\line(0,1){0.45}}

\linethickness{0.5mm}
\put(2,2){\vector(0,1){0.55}}
\put(2,2.55){\line(0,1){0.45}}

\linethickness{0.5mm}
\put(1,3){\vector(0,-1){0.55}}
\put(1,2.45){\line(0,-1){0.45}}

\linethickness{0.5mm}
\put(4,3){\vector(0,-1){0.55}}
\put(4,2.45){\line(0,-1){0.45}}

\linethickness{0.5mm}
\put(3,2){\vector(0,-1){0.55}}
\put(3,1.45){\line(0,-1){0.45}}

\linethickness{0.5mm}
\put(3,1){\vector(0,-1){0.55}}
\put(3,0.45){\line(0,-1){0.45}}
\end{picture}

\end{center}
\caption{An example of an element of $\FF$ with 5 roots
on a $2$-dimensional $5 \times 4$ box of $\Z^2$.}
\label{torus}
\end{figure}

The aim of this paper is to survey some recent results~\cite{AG1,AG2,ACGM1,ACGM2} on a certain measure on spanning forests of a given graph and its applications within the context of networks analysis.
We call a {\bf network} on $n\in\mathbb{N}$ vertices a {\bf directed and weighted} graph $$\GG=(\VV,\EE, w),$$  
where $\VV$ denotes a {\bf finite vertex set} of size $|\VV|=n$, $\EE$ stands for a {\bf directed edge set} seen as a prescribed collection of ordered pairs of vertices $\{ (x,y) \in \VV\times\VV \}$,  
and $w: \VV \times \VV \mapsto \R^+$ is a {\bf weight function}, which associates to each ordered pair $(x,y)\in\EE$  a strictly positive weight $w(x,y)$. 
We will consider {\bf irreducible} networks where for every two distinct vertices $x,y\in\VV$,  there is always a directed path connecting them, that is, a sequence $\{e_i=(x_i,y_i)\}_{i=1}^{l}\subset \EE$ for some $l\in\mathbb{N}$ such that $x_1=x, y_l=y$ and $y_i=x_{i+1}$ for every $i\leq l-1$.
Let us introduce the measure at the core of this work. 
A {\bf rooted spanning forest} $\phi$ is a subgraph of $\GG$ without cycle,
with $\VV$ as set of vertices
and such that, for each $x \in \VV$, there is at most one $y \in \VV$
such that $(x, y)$ is an edge of $\phi$.
The {\bf root set} $\RR(\phi)$ of the forest $\phi$ is the set of points $x \in \VV$
for which there is no edge $(x, y)$ in $\phi$; the connected components of $\phi$
are {\bf trees}, each of them having edges that are oriented towards its own root.
We call $\FF$ the {\bf set of all rooted spanning forests} and we see each element $\phi$ in $\FF$ as a subset of $\EE$. See Figure~\ref{torus}. 
For fixed positive $q\in \R^+$, we are interested in the random forest $\Phi_q$ defined as the random variable with values in $\FF$ with law: 
\begin{equation}\label{nu}
\P(\Phi_q=\phi)=\frac{w(\phi) q^{|\RR(\phi)|}}{Z(q)},
\quad\quad \phi \in \FF,
 \end{equation}  
where $ w(\phi) = \prod_{e \in \phi} w(e)$ is the weight associated to the forest $\phi\in \FF$,  $|\RR(\phi)|$ is the number of trees, which is also the number of roots, and 
$Z(q) = \sum_{\phi \in \FF} w(\phi) q^{|\RR(\phi)|}$ is the {\bf normalizing partition function}.
In particular $\emptyset \in \FF$ is the spanning forest made
of $n$ degenerate trees reduced to simple roots and $w(\emptyset) = 1$.
We can include the case $q = +\infty$ in our definition by setting
$\Phi_\infty = \emptyset \in \FF$ in a deterministic way.
In the sequel we denote by $\E$ {\bf expectation w.r.t. the random forest law} $\P$.

\vspace*{1cm}

\begin{figure}
\centering
\includegraphics[height=8cm]{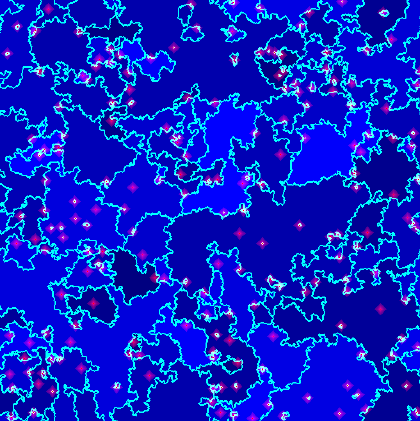}

\caption{A sample of $\PP(\Phi_q)$ and $R(\Phi_q)$ 
on a 2-dimensional box of $\Z^2$ with constant unitary nearest-neigbour weights and periodic boundary conditions. 
We used different shades of blue for blocks in the partition identified by different trees. Leaves are cyan, and roots are at the centers of red diamonds.}
\label{2dtorus}
\end{figure}

\bigskip Let us notice that the random forest $\Phi_q$ induces a partition of the graph into trees, and hence the measure in~\eqref{nu} can be seen on the one hand as a clustering measure similar in spirit to the well-known FK-percolation~\cite{G10}. On the other hand, the forest $\Phi_q$ is rooted and the set of roots $\RR(\Phi_q)$ forms an interesting random subset of vertices whose distribution can be explicitly characterized. Figure~\ref{torus} and~\ref{2dtorus}, respectively, show a realization of $\Phi_q$ in a simple geometrical setting and how the associated {\bf partition} denoted by $\PP(\Phi_q)$ and set of roots  $\RR(\Phi_q)$ look like. 
As we will show,  the presence of the tuning parameter $q$, controlling size and number of trees, and related efficient sampling algorithms make this measure particularly flexible and suitable for applications. 

\subsection{\large{Content of the paper}} 
This survey is organized in two parts. The first three sections constitute a first foundational part, followed by a second part on applications presented in the remaining sections. 
We will start by presenting basic properties of the measure in~\eqref{nu},
section~\ref{theory}, and some sampling algorithmic counterparts, section~\ref{algo}.
We then move to three main applications:
\begin{enumerate}
\item[{\bf (I)}]
In section~\ref{app1}, we will show how the set of roots $\RR(\Phi_q)$ can be used to define a 
probabilistic notion of {\bf \emph{subset of well distributed points in a graph}} and to practically sample it.
\item[{\bf (II)}] 
As a second application, a {\bf \emph{network coarsening scheme}} is presented in section~\ref{app2}, based on the forest $\Phi_q$ and the notion of Markov chain intertwining.
Motivations stem from questions in metastability and in signal processing and we provide two different algorithms and related theorems to control the quality of the resulting coarse grained network.
\item[{\bf (III)}] 
As a third application, supported by theorems and related experiments, in section~\ref{app3} we propose a {\bf \emph{wavelet basis construction and a multiscale algorithm for signal   processing  on arbitrary networks}}. 
\end{enumerate}
To conclude this introductory part, we describe briefly the origin of this measure and related literature, section~\ref{UST}, and introduce some basic crucial objects and notation needed along the paper, section~\ref{Markov.sec}. 
Apart from Wilson's algorithm, Theorem~\ref{genevieve}, Propositions~\ref{march} and ~\ref{SchurProp}, and parts of Theorems~\ref{borgorosa} and~\ref{poireau}, all other statements, algorithms and experiments are original and they were recently derived by the authors. 

\subsection{\large{Uniform spanning tree and a zoo of random combinatorial models}}\label{UST}
The Uniform Spanning Tree (UST) constitutes a by now classical topic in probability theory and it consists of a random spanning tree sampled uniformly among all possible spanning trees of a given graph.
The analysis of this object can be traced back at least to the work of Kirchhoff~\cite{K847} where the number of spanning trees of a graph was characterized in terms of determinants of minors of the corresponding discrete Laplacian matrix (matrix-tree theorem).
In the last decades, UST has been playing a central role in probability and statistical mechanics due to its deep relation with Markov chain theory and its surprising connection to  a number of challenging random combinatorial objects of current interest: e.g. loop erased random walks, percolation, dimers, sandpile models, Gaussian fields. We refer to~\cite{BLPS01,G10, LJ11,LP17,LPW08, Jarai, BP93,F93,K11,KW} for an overview on the vast literature on the subject.
What makes this object particularly interesting is its determinantal nature, namely, related local statistics have a closed-form expression in terms of determinants of certain kernels, together with the fact that there's an efficient random-walk-based-algorithm due to Wilson~\cite{[Wi]} for sampling, which will be presented in section~\ref{Wilson}.
The forest measure in~\eqref{nu} is actually a simple variant of the UST measure and it appeared already as a remark in~\cite{[Wi]} where the author mentions how to sample it.  
As for the UST measure, there are many interesting questions related to scaling and infinite volume limits of observables of the forest measure in~\eqref{nu}. Nonetheless, our focus here is on applications within the context of networks analysis and for this reason we will not insist on this very interesting fundamental line of investigation and we will restrict only to networks with finite number of vertices.

\subsection{\large{Basic objects and notation: random walks, graph Laplacian, Green's function}}\label{Markov.sec} 
Given the finite, irreducible, directed and weighted graph $\GG=(\VV,\EE, w)$ on $n$ vertices, let $X = \{X(t) : t \geq 0\}$ be {\bf the irreducible continuous-time Markov process} with state space ${\VV}$ and {\bf generator} $L$ given by
\begin{equation}\label{gen}
    (L f)(x) = \sum_{y \in \VV} w(x, y)\bigl[f(y) - f(x)\bigr],
    \qquad f:\VV \rightarrow \C,
    \qquad x \in \VV,
\end{equation}
In view of the finiteness and irreducibility assumptions, there exists a unique {\bf invariant measure for the Markov process} $X$ which will be denoted by $\mu$.
Averages of functions w.r.t. $\mu$ will be denoted by $\mu(f)$. We recall that the invariance of $\mu$ is equivalent to $\mu(L f)=0$ for arbitrary functions $f$.
We will denote by $P_x$ and $E_x$, respectively, {\bf law and expectation w.r.t. the random walk} $X$ starting from $x\in\VV$. 
Note that $L$ acts on functions as the matrix, still denoted by $L$:
\begin{equation}\label{Lmat} L(x,y)=w(x,y)  \mbox { for } x \neq y \, ; \, \, L(x,x)= -\sum_{y\neq x} w(x,y)\,. 
\end{equation}
We refer to the operator $-L$ as {\bf the weighted graph Laplacian} and we set 
\begin{equation}\label{wmax}w_{max}=\max_{x\in\VV} -L(x,x).\end{equation} 
 We will denote by $\hat X$ the {\bf discrete-time skeleton chain} associated to $L$ defined as the Markov chain with state space $\VV$ and transition matrix 
\begin{equation}\label{skeleton}P= \frac{L}{w_{max}}+\Id_{\VV},\end{equation}
with $\Id_{\VV}$ being the {\bf identity matrix} of size $|\VV|$.

\noindent 
In the sequel we deal with restrictions of various matrices: for any $\AA, \BB \subset \VV$ and any matrix $M = \bigl(M(x, y) : x, y \in \VV\bigr)$,
we write $[M]_{\AA, \BB}$ for the {\bf restricted matrix}
$$ [M]_{\AA, \BB} = \bigr(M(x, y) : x \in \AA, y \in \BB \bigl) \, .$$
In case $\AA=\BB$, we will simply write $[M]_{\AA}$.


\noindent 
For a subset $\AA \subset \VV$
\begin{equation}\label{hit}T_{\AA}=\inf\{t \geq 0: X(t)\in \AA \} \text{       is {\bf the hitting time} of the set } \AA,\end{equation}
with the convention $T_{\emptyset}=\infty$.

\noindent 
Finally, for a given (possibly random) time $T$ and arbitrary $x,y\in\VV$, we write :
\begin{equation}\label{Green} G_T(x,y)= E_x\Biggl[\int_0^{T}{\mathbbm 1}_{\{X(t) = y\}}\,dt\Biggr] \text{    for {\bf the Green's function up to time }}  T,\end{equation} 
i.e. the mean time $X$ spends in $y$ up to time $T$ when starting from $x$. In case $T$ is a random variable independent of the random walk, we will slightly abuse notation and still use $E_x$ for expectation w.r.t. the random walk and the extra randomness.

\noindent
Let us conclude with basic notation for normed spaces. We will denote by 
\begin{equation}\label{pspace}\ell_{p}(\VV, \mu)= \{  f: \VV \rightarrow \R \mid \nor{f}_{p,\VV} <\infty \}, \quad p\geq 1,\end{equation}
the {\bf $\ell_p$--space} of functions on $\VV$ w.r.t. the norm 
\begin{equation}\label{pnorm}\nor{f}_{p,\VV}  = \pare{\sum_{x \in \VV} \va{f(x)}^p  \mu(x)}^{1/p}\, .\end{equation}

Further, for arbitrary probability measures $\nu_1$ and $\nu_2$ on $\VV$, their {\bf total variation distance} is given by  
\begin{equation}\label{TV} d_{TV}(\nu_1, \nu_2) = \frac 1 2 \sum_{x \in \VV} \abs{\nu_1(x)-\nu_2(x)} \, .
\end{equation}

\section{\large{Random spanning forests: Laplacian spectrum and determinantality}} \label{theory}

We start here an account of the basic fundamental results characterizing the distribution of the main objects related to the random forest $\Phi_q$.
It will be convenient for the sequel to consider the following {\bf generalized version of the random forest}.
For any $\BB\subset \VV$ 
we denote by $\Phi_{q, \BB}$ a random variable in $\FF$
with the law of $\Phi_q$ conditioned on the event $\bigl\{\BB\subset \RR(\Phi_q)\bigl\}$.
We then have, for any $\phi$ in $\FF$,
\begin{equation}\label{extendForest}
    {\mathbb P}\bigl(\Phi_{q, \BB} = \phi\bigr)
    = {w(\phi) q^{|\RR(\phi)| - |\BB|} \over Z_\BB(q)} {\mathbbm 1}_{\{\BB\subset \RR(\phi)\}}
\end{equation}
with 
\begin{equation}
    Z_\BB(q) = \sum_{\phi : \RR(\phi) \supset \BB} w(\phi) q^{|\RR(\phi)| - |\BB|}.
    \label{farfalla}
\end{equation}
The original definition in Equation~\eqref{nu} is recovered by simply setting $\BB=\emptyset$, so that $\Phi_q=\Phi_{q, \emptyset}$ and $Z(q)=Z_\emptyset(q)$.  
This extended law is non-degenerate even for $q = 0$, provided that $\BB$ is non-empty.
And if $\BB$ is a singleton $\{r\}$, then $\Phi_{0, \{r\}}$ is {\bf the classical random spanning tree
with a given root} $r$, namely, a spanning tree $\tau$ rooted at $r$ sampled with probability proportional to $\prod_{e \in \tau} w(e)$.  
Let us emphasize that actually, for $q > 0$, $\Phi_q = \Phi_{q, \emptyset}$ itself is also a special case of the usual random spanning tree
on the extended weighted graph $\GG' = (\VV',\EE', w')$
obtained by addition of an extra point $r$ to $\VV$ ---to form $\VV' = \VV \cup \{r\}$---
and by adding extra edges to $r$ with weights $ w'(x, r) = q$ and $w'(r, x) = 0$ for all $x$ in $\VV$.
Indeed, to get $\Phi_q$ from the random spanning tree on $\VV'$ rooted in $r$,
one only needs to remove all the edges going from ${\VV}$ to $r$.

As a first result we characterize the partition function in terms of the Laplacian spectrum. 
To this end we denote by  
$\lambda_{0, \BB}$, $\lambda_{1, \BB}$,~\dots, $\lambda_{l - 1, \BB}$,
with $l = |\VV \setminus \BB|$ and some $\BB \subset \VV$, {\bf the eigenvalues }
of $[-L]_{\VV \setminus \BB}$.

\begin{theo}[{\bf Partition function and Laplacian spectrum}]\label{genevieve}
 For any $\BB\subset \VV$, $Z_\BB$ is the characteristic polynomial of $[L]_{\VV \setminus \BB}$, i.e., 
    \begin{equation*}
         Z_\BB(q)
        = \det\left[q \Id_{\VV \setminus \BB} - [L]_{\VV \setminus \BB}\right]
        = \prod_{j < l}(q + \lambda_{j, \BB}),
        \quad q \in {\mathbb R}.
        \label{aria}
    \end{equation*}
\end{theo}
The above result can be seen as a version of the well-known matrix-tree-theorem (see for instance~\cite{AT}). The proof can be derived in several ways. We refer to~\cite{AG1} for an elementary proof by a classical argument based on loop-erased random walks and using the notation herein.    

As mentioned above, one of the nice features of the forest measure (as well as for the random spanning tree measure) is its determinantal structure which allows for explicit computations. 
We start by recalling the determinantality of the edge set for which we need some more notation. 

For an oriented edge $e=(x,y)$ we denote  the {\bf starting and ending vertex}, respectively, as $e_-=x$ and $e_+=y$.
For any given $\BB\subset \VV$ and $q>0$, we write $G_{q,\BB}$ for the Green's function in Equation~\eqref{Green} with $T=T_q \wedge T_{\BB}$ , the minimum between 
the hitting time of $\BB$ (see Equation~\eqref{hit}) and time $T_q$ denoting an {\bf independent exponential time of parameter $q$}.
It is not difficult to see that $G_{q,\BB}$ can be identified with the operator $[q \Id- L]_{\VV\setminus {\BB}}^{-1}$ so that,   
$$G_{q,\BB}(x,y)=E_x\left[
		\int_0^{T_q \wedge T_{\BB}}{\mathbbm 1}_{\{X(t)=y\}} dt\right]
		= \left\{ \begin{array}{ll}
		[q \Id- L]_{\VV\setminus {\BB}}^{-1}(x,y) & \mbox{ for } x,y\in\VV\setminus {\BB} \, 
			\\
			0 & \mbox{otherwise. } 
			\end{array} \right.
$$
For $x\in \VV$ and $e$ in $\EE$,
we call 
\begin{equation*}
	J_{q,\BB}^+(x, e) = G_{q,\BB}(x, e_-) w(e)
\end{equation*}
the {\bf expected number of crossings
of the (oriented) edge $e$ up to time $T_q \wedge T_{\BB}$}, and
\begin{equation*}
	J_{q,\BB}(x, e) = J_{q,\BB}^+(x, e) - J_{q,\BB}^+(x, -e)
	\,,
\end{equation*}
the {\bf net flow through $e$
starting from $x$}.

\begin{theo}[{\bf Determinantal edges: transfer-current}]\label{borgorosa}
Fix $\BB\subset \VV$ and $q>0$. Then, for any $\AA_k=\{e_1, \dots, e_k\} \subset \EE:$
	\begin{equation*}\label{tenda}
		{\mathbb P}\left(
			\AA_k \subset \Phi_{q,\BB}
		\right)
		= {\mathbb P}\left(
			e_1, e_2, \dots, e_k \in \Phi_{q,\BB}
		\right)
		= {\rm det} \left[I_{q,\BB}^+\right]_{\AA_k}
	\end{equation*}
	with 
	\begin{equation}\label{current}
		I_{q,\BB}^+(e, e') = J_{q,\BB}^+(e_-, e') - J_{q,\BB}^+(e_+, e')
		\,, \quad e, e' \in \EE
		\,.
	\end{equation}
	In addition,
	denoting by $\{\pm e_1, \dots, \pm e_k \in \Phi_{q,\BB}\}$
	the event that for all $i \leq k$ either $e_i$ or $-e_i$
	belong to $\Phi_{q,\BB}$, it holds
	\begin{equation}\label{anatra}
		{\mathbb P}\left(
			\pm e_1, \dots, \pm e_k \in \Phi_{q,\BB}
	\right)
	={\rm det}\left[I_{q,\BB}\right]_{\AA_k}
	\end{equation}	
	with 
	\begin{equation}
		I_{q,\BB}(e, e') = J_{q,\BB}(e_-, e') - J_{q,\BB}(e_+, e')
		\,, \quad e, e' \in \EE
		\,.
	\end{equation}
\end{theo}

The above theorem is a version of the celebrated transfer-current theorem due to Burton and 
Pemantle~\cite{BP93}. In its original form, this theorem was proven in an undirected graph, 
extensions to the directed setup have appeared for e.g. in the recent Chang~\cite{C13} and the statement in Theorem~\ref{borgorosa} is nothing but a probabilistic reformulation of Thm.~5.2.3 and Coro.~5.2.4 in~\cite{C13}. 
For a simple proof using our notation, we refer the reader to Prop. 3.1 in~\cite{AG2}.

Theorem~\ref{borgorosa} says that $\Phi_{q,\BB}$ is a determinantal process with kernel $I_{q,\BB}^+$ interpretable in terms of random-walk-flow. If from a computational point of view, this allows to get explicit formulas, from a phenomenological perspective, being determinantal means that the corresponding objects tend to repel each other, more precisely, they are {\bf negatively correlated}:
$$\P(e_1, e_2 \in  \Phi_{q,\BB})\leq \P(e_1\in  \Phi_{q,\BB})\P(e_2 \in  \Phi_{q,\BB})\quad \text{for any } e_1,e_2\in\EE.$$

Inherited by the determinantal nature of $\Phi_{q,\BB}$, also the set of roots $\RR(\Phi_{q,\BB})$ is determinantal with a remarkable stochastic kernel given by the random walk $X$ killed at time $T_q \wedge T_{\BB}$: 
\begin{theo}[{\bf Determinantal roots with killed random walk kernel}]\label{detRoots}
Fix $\BB\subset \VV$ and $q>0$. Then, for any $\AA\subset \VV:$
	\begin{equation*}
		{\mathbb P}\left( \AA\subset \RR(\Phi_{q,\BB})\right)
		= {\rm det} \left[K_{q,\BB}\right]_{\AA},
	\end{equation*}
	with 
	\begin{equation*}
		K_{q,\BB}(x,y):=qG_{q,\BB}(x,y)=P_x\bigl(X(T_q \wedge T_{\BB}) = y\bigr), \quad x, y \in \VV
		\,.
	\end{equation*}
In case $\BB=\emptyset$,  $\Phi_{q,\emptyset}=\Phi_q$ and we simply write 
\begin{equation}\label{kernel2}
		K_{q}(x,y):=P_x\bigl(X(T_q) = y\bigr), \quad x, y \in \VV.
\end{equation}
\end{theo}

This theorem has been derived in~\cite{AG1}, see Prop.~2.2 therein.
We next move to the characterization of $|\RR(\Phi_{q,\BB})|$, that is, the number of roots/connected components/trees.
The next statement corresponds to Prop.~2.1 in~\cite{AG1}.
  
\begin{theo}{\bf ({Number of roots})}\label{RootsNumber}
Fix $\BB\subset \VV$ and $q\geq0$ and let $l = |\VV \setminus \BB|$.
Set  \begin{equation}\label{pj}
 p_j(q) = {q \over q + \lambda_{j, \BB}}\,,\quad 0\leq j \leq l-1,\end{equation}

Decompose 
$$
    J_0 = \bigl\{j \leq l-1 : \lambda_{j, \BB} \in {\mathbb R}\bigr\},\qquad
    J_+ = \bigl\{j \leq l-1 : {\rm Im}(\lambda_{j, \BB}) > 0\big\},\qquad
    J_- = \bigl\{j  : {\rm Im}(\lambda_{j, \BB}) < 0\big\},
$$
and define independent random variables $B_j$'s and  $C_j$'s, respectively, with laws: 
$$
    \P(B_j = 1) =  p_j(q), \quad \P(B_j = 0) = 1 - p_j(q), \quad j \in J_0,$$
and for $ j \in J_+$, 
$${\mathbb P}(C_j = 2) = |p_j(q)|^2,
    \quad {\mathbb P}(C_j = 1) = 2 {\rm Re}\left(p_j(q)\right) - 2|p_j(q)|^2,
    \quad {\mathbb P}(C_j = 0) = 1 - 2{\rm Re}\left(p_j(q)\right) + |p_j(q)|^2 \, .$$

Then, whenever $q>0$ or $\BB\neq\emptyset$, the random variable $|\RR(\Phi_{q,\BB})|$ is distributed as 
$$S_{q, \BB} = |\BB| + \sum_{j \in J_0} B_j + \sum_{j \in J_+} C_j.$$
\end{theo}

Notice that in case the spectrum of the graph Laplacian $-L$ is real and $\BB$ is empty,  $|\RR(\Phi_{q})|$ is simply given by the sum of independent Bernoulli's $B_j$'s ($J_+$ being empty). 
In particular, since $\lambda_0=0$, $B_0 = 1$ and we recover the fact that $ \abs{\RR(\Phi_q)} \geq 1$ a.s. 

Further, we emphasize that in general momenta of the $|\RR(\Phi_{q})|$ have simple expressions and can be easily obtained by differentiating w.r.t. $q$ the normalizing partition function $Z(q)$. 
For example, mean and variance are given by \begin{equation}\label{2moments}\E\bigl[ |\RR(\Phi_{q}) | \bigr]
    = \sum_{j < n} {q \over q + \lambda_j}\quad\text{  and  }\quad
Var\bigl[ |\RR(\Phi_{q}) | \bigr]
    = \sum_{j < n} {q \over q + \lambda_j}- \Bigl({q \over q + \lambda_j}\Bigr)^2.\end{equation}

\subsection{\large{Dynamics: forests, roots and partitions.}}
Before moving to sampling algorithms, it is worth mentioning that it is possible to construct a stochastic process with values in $\FF$ which allows to couple at once all $\Phi_q$'s as $q$ varies in $\R^+$. 
A few comments on this coupling are postponed to section~\ref{traj} and Figure~\ref{chocolat}. We state here the main theoretical results and collect some related remarks.
The following statement corresponds to Thm. 2 in~\cite{AG1}.

\begin{theo}[{\bf Forest dynamics: coupling all q's}]\label{alain}
    There exists a (non-homogeneous) continuous-time Markov process $F=\{F(s): s\geq 0\}$ with state space $\FF$
    that couples together all forests $\Phi_q$ for $q > 0$ as follows: 
    for all $s \geq 0$ and $\phi \in \FF$ it holds
    $$
        \P(F(s) = \phi) = \P(\Phi_{1 / t} = \phi) = \P(\Phi_q = \phi)
    $$
    with $t = 1/q$, $s = \ln(1 + w_{max} t)$ and $w_{max}$ as in~\eqref{wmax}.
\end{theo}

The coupling $t \mapsto \Phi_{1/t} = F(\ln(1 +  w_{max} t))$
is associated with a fragmentation and coalescence process,
for which coalescence is strongly predominant,
and at each jump time
one component of the partition is fragmented into pieces 
that possibly coalesce with the other components.
In particular, the process $F$ starts from $F(0)=\Phi_\infty=\emptyset \in \FF$, that is,  the degenerate spanning forest made
of $n$ trees reduced to simple roots, and eventually reaches in finite time a unique spanning rooted tree. 

As a corollary of the above coupling theorem, we get a determinantal characterization of  the finite-dimensional distributions of the process $t \mapsto \RR(\Phi_{1 / t})$, which can be seen as a dynamical extension of
Theorem~\ref{detRoots}.

\begin{cor}[{\bf Dynamic roots distribution}]
 For any choice $0 < t_1 < \cdots < t_k < t_{k + 1} = 1 / q_{k + 1}$ and any sequence 
  $\AA_1$, \dots, $\AA_k$, $\AA_{k + 1}$ of subsets of $\VV$,
    it holds
    \begin{equation}
        \renewcommand{\=}[1]{\AA_{#1} \subset \RR(\Phi_{1 / t_{#1}})}
        \begin{split}
            {\mathbbm P}&\bigl(\={k + 1} \bigm| \=k, \dots, \=1\bigr)\\
            &= \sum_{\BB_k \subset \AA'_k} 
            \sum_{\BB_{k - 1} \subset \AA'_{k - 1}} 
            \cdots
            \sum_{\BB_{1} \subset \AA'_{1}} 
            \prod_{i = 1}^k \biggl({t_i \over t_{k + 1}}\biggr)^{|\BB_i|} \biggl(1 - {t_i \over t_{k + 1}}\biggr)^{|\AA'_i \setminus \BB_i|}
            \det \bigl[K_{q_{k + 1}, \BB}\bigr]_{\AA_{k + 1}} \\
            \hbox{with}& 
            \quad \AA'_k = \AA_k,
            \quad \AA'_{k - 1} = \AA_{k - 1} \setminus \AA_k,
            \quad \dots
            \quad \AA'_1 = \AA_{1} \setminus (\AA_k \cup \AA_{k - 1} \cup \cdots \cup \AA_2)\\
            \quad \hbox{and}&
            \quad \BB = \bigcup_{i = 1}^k \BB_i.
        \end{split}
        \label{blandine}
    \end{equation}
\end{cor}

This statement corresponds to Prop. 2.4 in~\cite{AG1}. We do not have a similar characterization 
for the partition process  $t \mapsto {\mathcal P}(\Phi_{1 / t})$. 
More generally, as we saw, while a precise understanding and characterization of $\Phi_q$ and $\RR(\Phi_q)$ is possible, 
we know very little about the induced partition ${\mathcal P}(\Phi_{q}).$ 

We conclude this part on the relevant theoretical results by mentioning a last property of the root set, namely, by conditioning on the induced partition ${\mathcal P}(\Phi_{q})$ 
the roots are distributed according to the equilibrium measure $\mu$ of the random walk $X$ restricted to each component of the partition:

\begin{theo}[{\bf Roots at restricted equilibria}]\label{maude}
Let $[\AA_1, \dots, \AA_m]$ be denoting an arbitrary partition of $\VV$ in $m \leq n$ subsets and fix $r_i\in\AA_i, i=1,\ldots,m$. Then
    $$ 
        \P\Bigl(\RR(\Phi_q) = \{r_1, \dots, r_m\} \Bigm| {\mathcal P}(\Phi_q) = [\AA_1, \dots, \AA_m]\Bigr)
        = \prod_{i = 1}^m \mu_{\AA_i}(r_i)
    $$ 
provided that the conditioning has non-zero probability, with $\mu_{\AA_i}$ denoting the invariant measure of the restricted dynamics with generator $L_i$ defined by
$$
    (L_i f)(x) = \sum_{y \in \AA_i} w(x, y)\bigl[f(y) - f(x)\bigr],
    \quad x \in \AA_i,
    \quad f: \AA_i \rightarrow {\mathbbm C}.
$$     
\end{theo}
This statement is a consequence of the well-known Markov chain tree theorem (cf. e.g.~\cite{AT}), see Prop. 2.3 in~\cite{AG1}.

\section{\large{Sampling algorithms: Wilson's \& Co. }}\label{algo}
The flourishing literature around the random spanning tree theme is mainly due to Wilson's algorithm (cf.~\cite{[Wi]}), which  is
not only a practical procedure to sample $\Phi_{q,\BB}$, but actually also a powerful tool to analyze its law. 
The reader not acquainted with this topic is invited to look into the proofs of the results presented in section~\ref{theory} which heavily exploit 
the power of this algorithm in action.
We will start by recalling it, section~\ref{Wilson}. We then explain how to get an approximate sample of a forest with a prescribed number of roots, section~\ref{approxM}, and we conclude this sampling 
algorithmic part with some comments about sampling the forest dynamics in Theorem~\ref{alain}, section~\ref{traj}. 

\subsection{\large{Wilson's algorithm: sampling a forest for fixed $q$}}\label{Wilson} 
The following algorithm due to Wilson~\cite{[Wi]} samples 
$\Phi_{q, \BB}$ for $q > 0$ or $\BB\neq \emptyset$:
\begin{itemize}
\item[a.] start from $\BB_0 = \BB$ and $\phi_0 = \emptyset$, choose $x$ in $\VV \setminus \BB_0$ and set $i = 0$;
\item[b.] run the Markov process starting at $x$ up to time $T_q \wedge T_{\BB_i}$
    with $T_q$ an independent exponential random variable with parameter $q$
    (so that $T_q = +\infty$ if $q = 0$) and $T_{\BB_i}$ the hitting time of $\BB_i$;
\item[c.] with 
    $$
        \Gamma^x_{q, \BB_i} = (x_0, x_1, \dots, x_k)
        \in \{x\} \times \bigl(\VV \setminus (\BB_i \cup \{x\})\bigr)^{k - 1} \times \bigl(\VV \setminus \{x\}\bigr)
    $$
    the loop-erased trajectory obtained from $X : [0, T_q \wedge T_{\BB_i}] \rightarrow \VV$,
    set $\BB_{i + 1} = \BB_i \cup \{x_0, x_1, \dots, x_k\}$
    and $\phi_{i + 1} = \phi_i \cup \{(x_0, x_1), (x_1, x_2), \dots, (x_{k - 1}, x_k)\}$
    (so that $\phi_{i + 1} = \phi_i$ if $k = 0$);
\item[d.]\label{d} if $\BB_{i + 1} \neq \VV$, choose $x$ in $\VV \setminus \BB_{i + 1}$
    and repeat b--c with $i + 1$ in place of $i$,
    and, if $\BB_{i + 1} = \VV$, set $\Phi_{q, \BB} = \phi_{i + 1}$.
\end{itemize}

It is worth stressing that in steps a.~and d.~{\bf {\em the choice of the starting points $x$  is arbitrary}}, 
a remarkable fact which represents the main strength of this algorithm.

There are at least two ways to prove that this algorithm indeed samples $\Phi_{q, \BB}$ with the desired law.
One option is to  follow Wilson's original proof in~\cite{[Wi]},
which makes use of the so-called Diaconis-Fulton stack representation of Markov chains, cf.~\cite{[DF]}.
An alternative option is to follow Marchal
who first computes in~\cite{[Ma]} the law of the loop erased trajectory $\Gamma^x_{q, \BB}$
obtained from the random trajectory $X : [0, T_q \wedge T_\BB] \rightarrow \VV$
started at $x \in \VV \setminus \BB$ and stopped in $\BB$ or at an exponential time $T_q$
if $T_q$ is smaller than the hitting time $T_\BB$.
One has indeed:

\begin{prop}[{\bf Distribution of loop-erased walks }]\label{march}
    For any self-avoiding path $(x_0, x_1, \dots, x_k) \in \VV^{k + 1}$ 
    such that $x_0 = x, x_1, \cdots, x_{k-1}$  in $\VV \setminus \BB$, 
    it holds
    $$
        {\mathbb P}\bigl(\Gamma^x_{q, \BB} = (x_0, x_1, \dots, x_k)\bigr)
        = \begin{cases}
            \displaystyle \prod_{j < k} w(x_j, x_{j + 1})
            {\det[q \Id - L]_{\VV \setminus (\BB\cup \{x_0, \dots, x_{k - 1}\})}
            \over \det[q \Id- L]_{\VV \setminus \BB}} & \hbox{if $x_k \in \BB$,}\\
            \displaystyle q \prod_{j < k} w(x_j, x_{j + 1})
            {\det[q \Id- L]_{\VV \setminus (\BB\cup \{x_0, \dots, x_k\})}
            \over \det[q \Id - L]_{\VV \setminus \BB}} & \hbox{if $x_k \not\in \BB$.}
       \end{cases} 
    $$
\end{prop}

From this result one can easily compute the law of $\Phi_{q, \BB}$ following the steps of the algorithm above to get the law in Equation~\eqref{extendForest}.
Further, from this statement we see how determinants of the Laplacian emerge.
Concerning {\bf {\em the average running time of Wilson's algorithm}}, it is in general polynomial in the number of nodes $n$ and typically much smaller than the random walk cover time. In particular, it can be explicitly characterized in spectral terms as the sum of the inverse of the $n-|\BB|$ eigenvalues of $[-L]_{\VV\setminus\BB}$, see e.g. Prop. 1 in~\cite{[Ma]}.  

\subsection{\large{Forests with a prescribed number of roots: approximate sampling}}\label{approxM}
We have seen that Wilson's algorithm provides a practical way to sample $\Phi_q$. In applications, one might be interested in sampling 
 $\Phi_q$ conditioned on having a prescribed number of roots, that is, conditioned on $\bigl\{|\RR(\Phi_q)| = m\bigr\}$ for fixed $m<n$. 
Unfortunately, we do not know any efficient algorithm providing such an outcome. 
Nevertheless we can exploit Theorem~\ref{RootsNumber} to get a procedure to sample $\Phi_q$ with approximately $m$ roots, with an error of order $\sqrt m$ at most.
In fact, it is not difficult to check that ${\rm Var}\bigl(|\RR(\Phi_q)|\bigr)  \leq  2 \E\bigl[|\RR(\Phi_q)|\bigr]$ and, in view of~\eqref{2moments}, it suffices to find 
the solution $q^*$ of the equation \begin{equation}
    \sum_{j < n} {q \over q + \lambda_j} = m.
    \label{pascale}
\end{equation}

However, solving Equation~(\ref{pascale}) requires to compute the eigenvalues $\lambda_j$'s of $-L$ which is in general computationally costly especially if we are dealing with a big size network.
One way to find an approximate value of the solution $q^*$ is to use, on the one hand, the fact that $q^*$ is the only one stable attractor
of the recursive sequence defined by $q_{k + 1} = f(q_k)$ with 
$$
    f : q  \mapsto q \times {m \over \sum_{j < n} {q \over q + \lambda_j}} = {m \over \sum_{j < n} {1 \over q + \lambda_j}}\,,
$$
and on the other hand, the fact that $|\RR(\Phi_q)|$ and $\E\bigl[|\RR(\Phi_q)|\bigr] $
are typically of the same order, at least when $\E\bigl[|\RR(\Phi_q)|\bigr]$, i.e. $q$, is large enough,
since ${\rm Var}\bigl(|\RR(\Phi_q)|\bigr) / \E^2\bigl[|\RR(\Phi_q)|\bigr] \leq  2 / \E\bigl[|\RR(\Phi_q)|\bigr]$.
We then propose the following algorithm to sample $\Phi_q$ with $m \pm 2\sqrt m$ roots.
\begin{itemize}
\item[a.] Start from any $q_0 > 0$, for example $q_0 =w_{max}=\max_{x\in\VV} -L(x,x)$,
    and set $i = 0$.
\item[b.] Sample $\Phi_{q_i}$ with Wilson's algorithm.
\item[c.] If $|\RR(\Phi_{q_i})| \not\in \bigl[m - 2 \sqrt m, m + 2 \sqrt m\bigr]$,
    set $q_{i + 1} = m q_i / |\RR(\Phi_{q_i})|$ and repeat b with $i + 1$ instead of $i$.
    If $|\RR(\Phi_{q_i})| \in \bigl[m - 2 \sqrt m, m + 2 \sqrt m\bigr]$, then return $\Phi_{q_i}$.
\end{itemize}
We refer the reader to section 2.2 in~\cite{AG1} to argue that indeed this algorithm rapidly stops.

\subsection{\large{Coalescence-fragmentation process: sampling for different $q$'s at once}}\label{traj}
The Markov process $F$ in Theorem~\ref{alain} is based on the construction of a coalescence-fragmentation process with values in $\FF$ 
making use of Diaconis-Fulton's stack representation of random walks. 
For a detailed account on this algorithm and a number of related open questions, we refer the reader to section 2.3 in~\cite{AG1}. 

We mention that this algorithm allows to couple forests for different values of $q$'s.  The corresponding coupling is not monotone, in the sense that if $q'<q$, 
it is not true that $|\RR(\Phi_{q'})|\leq |\RR(\Phi_{q})|$ a.s. under the coupling measure, despite the fact that   
$\E\bigl[|\RR(\Phi_{q'})|\bigr]  < \E\bigl[|\RR(\Phi_{q})|\bigr]$, see e.g.~Equation~\eqref{2moments}. 
Yet this coupling is a very valuable tool in applications. In fact, it allows to practically sample $\Phi_{q'}$ starting from a sampled $\Phi_{q}$ for any $q'<q$, and 
more generally, by running this algorithm once in a chosen interval $[0,t^*]$, we get {\bf {\em samples of the whole forest trajectory}}  $(\Phi_{1 / t} : t \leq t^*)$.

\section{\large{Applications: Well distributed points in a network}}\label{app1} 
Given a map of a city modeled as a network $\GG$ where road crossings are identified with vertices, assume that we are interested in locating a number of energy plants at some crossings in such a way that the energy flowing out of these plants takes on average the same amount of time to reach every vertex of the city. 
If the energy flows according to a random walk $X$, we can rephrase the above problem by finding a subset $\BB\subset\VV$ for the locations of the energy plants such that for any $x\in\VV$: 
\begin{equation}\label{detectors}
E_x[T_{\BB}] \text{             is independent of  } x.
\end{equation}  
In other words, $\BB\subset\VV$ would constitute a set of {\bf \emph{well distributed points in the network}}.
We immediately realize that unless we have as many energy plants as the number of crossings, $\BB=\VV$, there is no deterministic proper subset $R$ satisfying the property in~\eqref{detectors}. 
Hence this notion of well distributed points is in principle meaningless but, by thinking in terms of {\bf \emph{disorder}}, we can turn it into a well-posed definition by finding a random set satisfying~\eqref{detectors} in an averaged sense. That is, a random set $\BB(\omega)\subset\VV$ is chosen according to some law with expectation $\E$ so that 
\begin{equation}\label{detectors2}
\E\Bigl[E_x\bigl[T_{\BB(\omega)}\bigr]\Bigr] \text{             is independent of  } x.
\end{equation}  
It then raises the question: does it exist such a random subset?

Our set of roots $\RR(\Phi_q)$ might be a good candidate since as previously noticed, the determinantality of the roots set stated in Theorem~\ref{detRoots} implies negative correlations suggesting that the roots in $\Phi_q$ tend to spread far apart from each other irrespective of the given network structure.
It turns out that indeed $\RR(\Phi_q)$ gives a positive answer to this question for any $q$. Further, as the next statement shows, the same is true when conditioning on having exactly $m$ energy plants, that is, with random subsets of prescribed size:
\begin{theo}[{\bf Well distributed roots}] \label{maschera}
For all $x \in \VV$ and all positive integer $m \leq n$ it holds
    \begin{equation*}
        \E\Bigl[E_x\bigl[T_{\RR(\Phi_q)}\bigr]\Bigr]
        = {1 \over q} \left(
            1 - \prod_{j > 0} {\lambda_j \over q + \lambda_j}
        \right)
\end{equation*}
and
\begin{equation}\label{fissom}
  \E\Bigl[E_x\bigl[T_{\RR(\Phi_q)}\bigr] \Bigm| |\RR(\Phi_q)| = m\Bigr]
        = {a_{m + 1} \over a_m}\,,
\end{equation}
where $a_{n + 1} = 0$, and $a_k$ denotes the coefficient of order $k$ of the characteristic polynomial of $L$.   
\end{theo}
This statement corresponds to Thm.~1 in~\cite{AG1}. When conditioning in Equation~\eqref{fissom} with either $m=1$ or $m=n-1$, it turns out that the property in~\eqref{detectors2} actually characterizes the law of $\BB(\omega)$.   
In particular, for $m=1$, the Markov tree-theorem, see e.g.~\cite{AT}, ensures that
conditionnally on  $\acc{|\RR(\Phi_q)| = 1}$, $\RR(\Phi_q)$ coincides with a point sampled according to the equilibrium measure $\mu$ of $X$. In this case,  Equation~\eqref{fissom}
was already known in the literature and often referred to as the {\bf \emph{ random target lemma}}, cf. e.g. Lemma 10.8 in~\cite{LPW08}. Our theorem is therefore a natural extension 
of this random target to subsets of arbitrary sizes. 
To get some insights in the corresponding proof, it is worth mentioning that the expected hitting time of a given deterministic set $\BB$ in Equation~\eqref{detectors} admits the following characterization due to Freidlin and Wentzell  in terms of forests, cf. Lemma 3.1 in~\cite{AG1} and Lemma~3.3 in~\cite{[FW]}: 
\begin{equation*} E_x\bigl[T_\BB\bigr]
    = \sum_{z \not\in \BB} G_\BB(x, z)
    = {1 \over Z_\BB(0)} \sum_{z \not\in \BB} \sum_{\scriptstyle \phi: \RR(\phi) = \BB \cup \{z\} \atop \scriptstyle \RR(\tau_x(\phi)) = z} w(\phi),
\end{equation*}
with $G_\BB$ denoting the Green's function in Equation~\eqref{Green} stopped at the hitting time $T=T_\BB$, and $\tau_x(\phi)$ standing for the tree in $\phi\in\FF$ containing $x\in\VV$. 

\begin{figure}
\centering
\includegraphics[height=10cm]{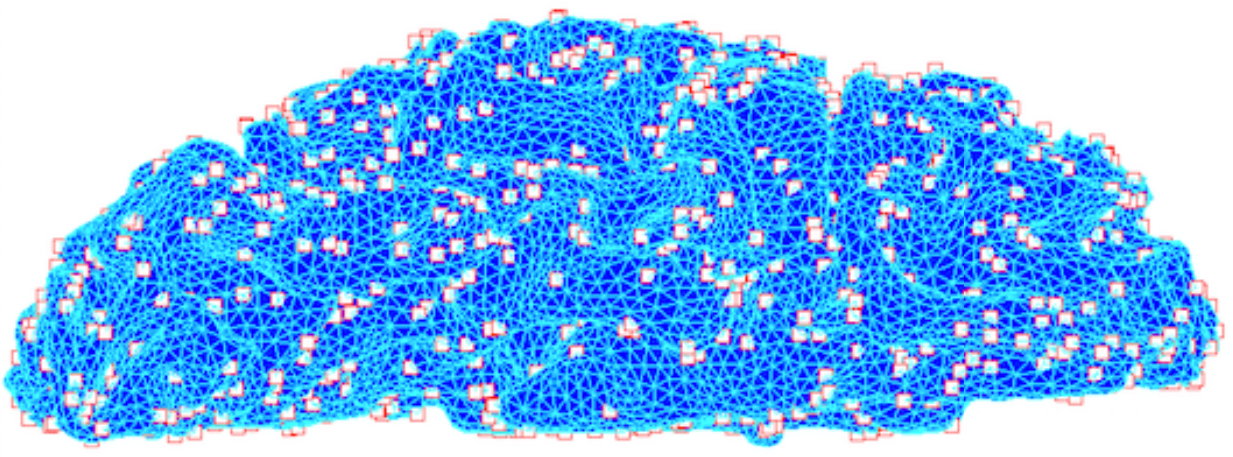}
\includegraphics[height=10cm]{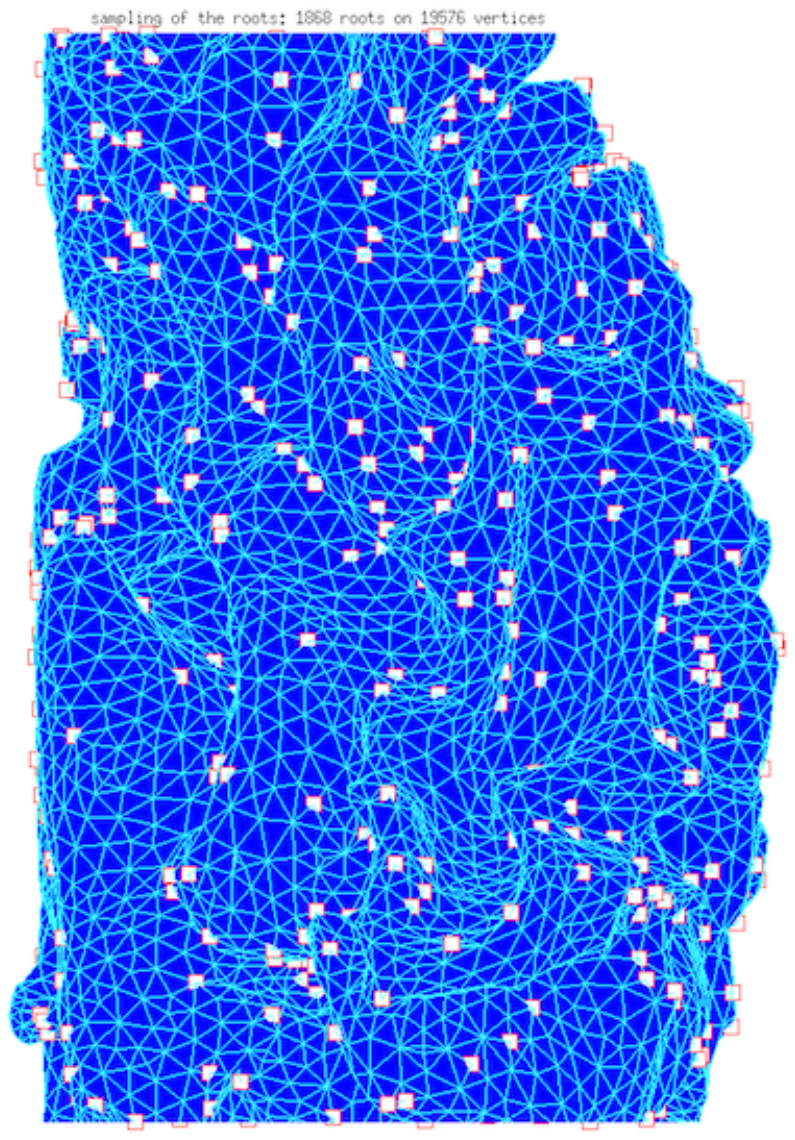}
\caption{Illustration of the roots in a non-trivial geometrical setup. 
The network is given by a triangulation of a brain cortex model due to J. Lef\`evre (LSIS, Marseille, France) with $1868$ roots denoted by white squares sampled from $(\Phi_q)$  out of $19576$ vertices.
}\label{brain}
\end{figure}

Concrete applications of this result will be given in the next sections to build up suitable subsampling procedures in the context of networks reduction and signal processing.
 
\section{\large{Applications: Network coarse graining via intertwining}}\label{app2} 
The goal of this section is to propose a random network coarse graining procedure which exploits the rich and flexible structure of the random spanning forest.
The problem can be formulated as follows:
\begin{center} {\bf (P1)} Given a network $\GG$ on $n$ nodes, find a smaller network $\bar{\GG}$ on $m<n$ nodes ``mimicking'' the original network $\GG$.
\end{center} 
For the sake of simplicity, we will restrict ourselves to {\bf networks in a reversible setup}, that is, when $\mu(x)w(x,y)=\mu(y)w(y,x)$ for any $x,y\in\VV$, with $\mu$ being the invariant measure of the Markov process X.
Now, it is a priori not clear what ``mimicking'' means and any meaningful answer would strongly depend on the implemented method and on the specific applications one has in mind.

Our approach will be process-driven. Namely, we saw that the structure of the starting network $\GG$ is encoded in the weighted graph Laplacian in Equation~\eqref{gen} 
which characterizes the Markov process $X$ with state space $\VV$. In view of the one-to-one correspondence between $\GG$ and the process $X$, we can look for a $\bar{\GG}$ in correspondence to another process $\bar{X}$  on a state space $\bar \VV$ with $m=|\bar \VV|<|\VV|=n$ points being some sort of coarse grained version of the process $X$. 
In this way, we are shifting perspective from graphs to Markov processes and, within this context, there is an interesting {\bf {\em duality notion called Markov chain intertwining}} which will be our lighthouse while addressing {\bf (P1)}. 
Before discussing this duality, we anticipate that concerning the {\bf{\em applications}}, our main motivation stems from two different problems:
\begin{itemize} 
\item
{\bf {\em metastability studies outside any asymptotic regimes}} (like low temperature or large volume limits), 
\item
  {\bf {\em signal processing on arbitrary networks}}. 
\end{itemize}
We will hence propose in sections~\ref{metagol} and~\ref{multigol} two different answers to {\bf (P1)}, namely, two network reduction procedures of general interest based on intertwining and random forests, inspired by and well-suited for the above problems. 
Our reduction schemes might also be useful for other applications, hence we will
start by discussing where the difficulties lie when looking at network reduction problems
through intertwining equations.
We start by introducing and discussing the intertwining duality.  

\subsection{\large{Intertwining and squeezing.}}\label{interdef}
Given two Markov chains with transition matrices $P$ and $\bar{P}$, and state spaces $\VV$ and $\bar{\VV}$, and given 
a rectangular stochastic matrix ${\Lambda}: \bar{\VV} \times \VV \rightarrow [0, 1]$,  
we say that {\bf the two chains are intertwined w.r.t. $\Lambda$ if }:
\begin{equation}\label{inter}
\Lambda P= \bar{P} \Lambda.
\end{equation}

Denoting by $\{\nu_{\bar x} = \Lambda(\bar x, \cdot) : \bar x\in \bar{\VV} \}$ the family of probability measures on $\VV$ identified by $\Lambda$, we see that 
this algebraic relation among matrices can be rewritten as 

\begin{equation} \label{banzai}
	\nu_{\bar x}P = \Lambda P(\bar x, \cdot) = \bar P \Lambda(\bar x, \cdot) = \sum_{\bar y \in \bar{\VV} } \bar P(\bar x, \bar y) \nu_{\bar y}
	, \quad \text{for all } \bar x \in \bar \VV,
\end{equation}
which says that: 

\begin{center}{\bf {\em the one-step-evolution of the $\nu_{\bar x}$'s according to $P$ 
 remains in their convex hull.}}
\end{center}

This duality notion can be equivalently formulated for continuous-time Markov processes by saying that 
 {\bf two Markov processes $X$ and $\bar X$ with generators $L$ and $\bar L$ and state spaces $\VV$ and $\bar \VV$ are intertwined  w.r.t. $\Lambda$ if }
\begin{equation}\label{intergen}\Lambda L= \bar{L} \Lambda.\end{equation}

By associating to the Markov process $X$ with generator $L$ in~\eqref{gen} the {\bf discrete-time skeleton chain} $\hat X$ as in~\eqref{skeleton}, 
we see that~\eqref{intergen} is equivalent to~\eqref{inter} if $\bar{P}= \frac{\bar{L}}{w_{max}}+\Id_{\bar{\VV}}$, and Equation~\eqref{banzai} reads as
: 
\begin{equation} \label{banzai2}
	\nu_{\bar x}L = \sum_{\bar y\in \bar{\VV} \setminus\{\bar{x} \}} \bar L(\bar x, \bar y) [\nu_{\bar y}-\nu_{\bar x}],
\end{equation}
which says again that, for each $\bar x\in \bar\VV$, by evolving the distribution $\nu_{\bar x}$ according to $L$, the process under consideration,  with rate  
$\bar L(\bar x, \bar y)$ is distributed according to $\nu_{\bar y}$. 

\subsubsection{\large{Intertwining in the literature.}}

Intertwining relations appeared in the context of diffusion processes in a paper by Rogers and Pitman~\cite{RP} as 
a tool to state identities in laws when the measures $\nu_{\bar x}$'s have disjoint supports. This method was later successfully applied 
to many other examples (see for instance~\cite{CPY},~\cite{MY}). In the context of Markov chains, intertwining was used by Diaconis and Fill~\cite{DF}
without the disjoint support restriction to build strong stationary times and to control convergence rates to equilibrium. 
At the time being, applications of intertwining include random matrices~\cite{DMDY}, particle systems~\cite{War}...

\subsubsection{\large{Solutions to intertwining equations, overlap and Heisenberg principle.}}
In the above references, intertwining relations have often been considered  with $|\bar \VV| $ being (much) larger than
or equal to $| \VV |$.
To address {\bf (P1)}, we will instead be naturally interested in the complementary case $|\bar \VV| < |\VV|$ and with the coarse grained process (or network identified by) $\bar P$ being irreducible. In this setup it is not difficult to show the existence of solutions $(\Lambda, \bar P)$ to Equation~\eqref{inter}, how they are related to the spectrum of $P$ , and how to construct some of them.
We refer the interested reader to section 2.2 in~\cite{ACGM1} and the statements therein. 
Let us simply note here
\begin{itemize}
\item that the intertwining equations generally have many solutions,
    including the trivial ones when, for any $\bar P$,
    all the $\nu_{\bar x}$ are equal to the invariant measure $\mu$ of $P$,
\item that the stability of the convex hull of the $\nu_{\bar x}$
    implies the stability of the vector space they span and the fact that the $\nu_{\bar x}$
    have to be linear combinations of at most $|\bar {\cal V}|$ left eigenvectors of $P$.
\end{itemize}
Looking at the eigenvectors of the generator $L$
as the analogue of the usual Fourier basis,
which diagonalizes the Laplacian operator,
we will say that
{\em the solutions $\nu_{\bar x}$ of the intertwining equations
have to be frequency localized}.
But we will see that the solutions
we will be interested in
for building a coarse-grained network 
should also be ``space localized''
or ``little overlapping''.
We will soon make precise 
what is meant here.
Let us simply stress for now that 
having both frequency and space localization
would contradict a (unfortunately not well established
for arbitrary graphs) Heisenberg principle.
To overcome this difficulty we will look at 
{\em approximate solutions} of intertwining equations
and we will focus on their {\em squeezing},
which is a measure of their joint overlap
or space localization,
that we can now introduce.

\subsubsection{\large{The squeezing functional.}}\label{compromise}
On the space of probability measures on $\VV$ with $|\VV|=n$, let us denote by $\langle \cdot, \cdot \rangle$ {\bf the scalar product} defined as
$$
    \langle \nu_{1}, \nu_{2}\rangle
    =\sum_{x\in\VV}{\nu_1(x) \over \mu(x)} {\nu_2(x) \over \mu(x)} \mu(x)
    =\sum_{x\in\VV}{\nu_1(x)\nu_2(x) \over \mu(x)}\,,
$$
for arbitrary probability measures $\nu_1,\nu_2$ on $\VV$, and let $\nor{\cdot}$ be {\bf the corresponding norm}. 
 
Given a family $\{\nu_{\bar x} = \Lambda(\bar x, \cdot) : \bar x\in \bar{\VV} \}$, since these measures form acute angles between them ($\langle \nu_{\bar x}, \nu_{\bar y}\rangle \geq 0$ for 
all $\bar x$ and $\bar y$ in $\bar{\VV}$) and have disjoint supports if and only if they are orthogonal,
one could use the volume of the parallelepiped they form to measure their ``joint overlap''.
The square of this volume is given by the determinant of {\bf the Gram matrix}:
$$
	{\rm Vol}(\Lambda) = \sqrt{\det\Gamma},
$$
with $\Gamma$ the square matrix on $\bar \VV$ with entries
$\Gamma(\bx, \by) = \langle \nu_{\bx},\nu_{\by}\rangle$, that is
\begin{equation}
\label{Gamma.def}
\Gamma= {\Lambda} D(1/\mu) {\Lambda}^t \, ,
\end{equation}
where $D(1/\mu)$ is the diagonal matrix with entries given by $(1/\mu(x), x \in \VV)$, and 
${\Lambda}^t$ is the transpose of ${\Lambda}$. 
Loosely speaking, {\bf {\em the less overlap, the largest the volume}}.

We will instead use the {\bf squeezing of $\Lambda$}, that we defined by
\begin{equation}
\label{Flat.eq}
\SS({\Lambda}):=\left\{\begin{array}{ccc}
&+\infty&\mbox{ if }\det(\Gamma)=0,\\
&\sqrt{\Tr\big(\Gamma^{-1}\big)} \in\,]0,+\infty[&\mbox{ otherwise,}
\end{array}\right. 
\end{equation} 
to measure this ``joint overlap''.
We call it ``squeezing'' not only because the $\nu_{\bar x}$ and the parallelepiped they form
is squeezed when ${\mathcal S}(\Lambda)$ is large,
but also because ${\mathcal S}(\Lambda)$ is  half the diameter of the rectangular parallelepiped
that circumscribes the ellipsoid defined by the Gram matrix $\Gamma$ :
this ellipsoid is squeezed too when ${\mathcal S}(\Lambda)$ is large.
We note finally that our squeezing controls the volume of $\Lambda$.
Indeed, by the comparison between harmonic
and geometric mean applied to the eigenvalue of the Gram matrix, {\bf {\em small squeezing implies large volume}}:
${\rm Vol}(\Lambda)^{1 / n} {\mathcal S}(\Lambda) \geq \sqrt n$.

The next statement, corresponding to Prop. 1 in~\cite{ACGM1}, gives bounds on this squeezing functional suitable for our approach.
\begin{prop}[{\bf Bounds on the squeezing functional}]\label{flatness.prop}
 Let $\{\nu_{\bar x}  : \bar x\in \bar{\VV} \}$ be a collection of $m=|\bar{\VV}|$ probability measures on $\VV$. 
\begin{itemize}
\item We have
\begin{equation} 
\label{Trace.ineq}
\SS({\Lambda}) \geq \sqrt{ \sum_{\bx \in \bV} \frac{1}{\nor{\nu_{\bx}}^{2}}} \, . 
\end{equation} 
Equality holds if and only if  the $\nu_{\bar x}$'s are pairwise orthogonal.
\item Assume that $\mu$ is a convex combination of the $\{\nu_{\bar x}  : \bar x\in \bar{\VV} \}$. 
Then,  
\[ \SS({\Lambda}) \geq 1 \, . 
\]
Equality holds if and only if  the $\nu_{\bar x}$'s are pairwise orthogonal.
\end{itemize}   
\end{prop}
\noindent
We notice in particular that $\SS({\Lambda})$ is maximal when the measures $\{\nu_{\bar x}  : \bar x\in \bar{\VV} \}$ are linearly dependent, and
minimal when they are orthogonal. Moreover, we know the minimal value of $\SS({\Lambda})$, 
when $\mu$ is a convex combination of the $\nu_{\bx}$'s.  Note that this is necessarily the case
if the convex hull of the $\nu_{\bar x}$ is stable under~$P$,
i.e. when~\eqref{inter} holds for some stochastic $\bar P$.
Indeed it is then stable under $e^{t{\mathcal L}}$ for any $t > 0$
and the rows of $\Lambda e^{t{\mathcal L}}$ converge to $\mu$ when $t$ goes to infinity.
We are now in shape to move to the applications.

\subsection{\large{Intertwining and metastability without asymptotics}}
A classical problem in meta\-sta\-bility studies can be described as follows. 
Associated with Markovian models
one is interested in making a coarse grained picture
of a dynamics which evolves on a large,
possibly very large,
configuration space,
according to some generator $L$.
A metastable state can be thought 
as a stationary distribution on this large configuration space
up to some random exponential time
that triggers a transition to a different metastable state,
possibly a more stable one.
By {\sl different} metastable states
we mean ``concentrated on different parts of the large configuration space.''
This is usually addressed in some asymptotic regime such as low temperature or large volume limits.
A natural and mathematically rigorous way (cf. Theorem~\ref{poireau} below)
to perform such a coarse graining, and avoiding limiting procedures,   
would be to provide solutions $(\Lambda, \bar L)$ to~\eqref{intergen} with
linking measures $\{\nu_{\bar x}=\Lambda(\bar x , \cdot) : {\bar x} \in {\bar \VV} \}$ having minimal squeezing $\SS({\Lambda})$. 
In this intertwining $\bar L$ would be the generator
of the coarse grained Markovian dynamics,
and the rows of $\Lambda$ would describe
the different metastable states.

\subsubsection{\large{A canonical example: the kinetic Ising model.}}\label{Ising} 
\begin{figure}[h]
$$
\includegraphics[width = 2.5 cm]{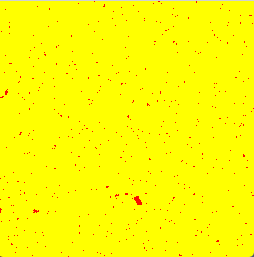}
\qquad
\includegraphics[width = 2.5 cm]{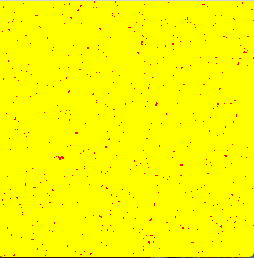}
\qquad
\includegraphics[width = 2.5 cm]{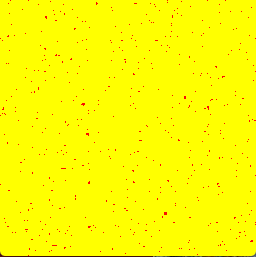}
$$
$$
\includegraphics[width = 2.5 cm]{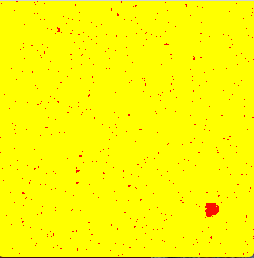}
\qquad
\includegraphics[width = 2.5 cm]{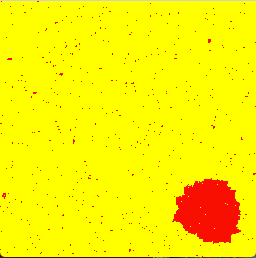}
\qquad
\includegraphics[width = 2.5 cm]{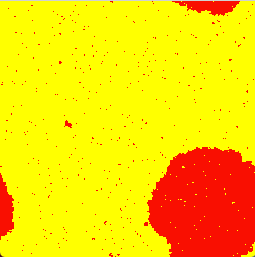}
\qquad
\includegraphics[width = 2.5 cm]{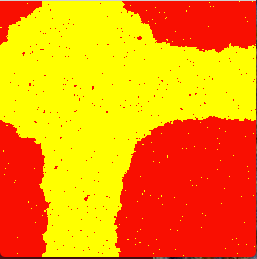}
\qquad
\includegraphics[width = 2.5 cm]{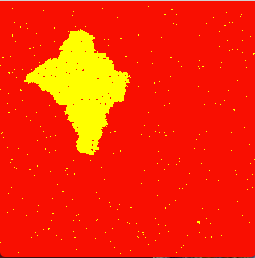}
$$
$$
\includegraphics[width = 2.5 cm]{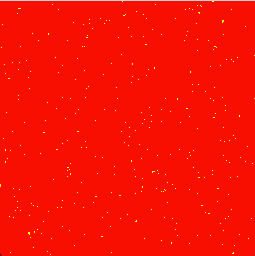}
\qquad
\includegraphics[width = 2.5 cm]{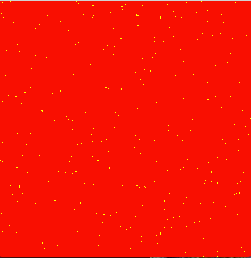}
\qquad
\includegraphics[width = 2.5 cm]{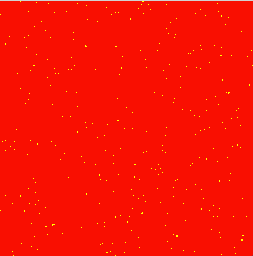}
$$

\caption{Snapshot of a kinetic Ising model
at times 471, 7482 and 13403 for the first line,
14674, 15194, 15432, 15892 and 16558 for the second line,
17328, 23645 and 40048 for the last line.}\label{pianta}
\end{figure}

In the following pictures (cf. Figure~\ref{pianta}) we consider
a Metropolis-Glauber kinetic Ising model
for a spin system started from aligned minus spins (yellow pixels on the pictures,
red pixels standing for plus spins)
and evolving under a small magnetic field $h = 0.14$,
at subcritical temperature $T = 1.5$,
in a $n \times n$ rectangular box $\BB_n$ with periodic boundary conditions and $n=256$ (thus here $\VV=\{+,-\}^{\BB_{256}}$).
The first three pictures can be thought 
as samples of a metastable state, 
which is concentrated on the minus phase of the system,
that is stationary up to the appearance (nucleation) of a supercritical droplet. This droplet
triggers the relaxation to a stable state,
which is concentrated on the plus phase,
and samples of which are given by the last three pictures.
In this case, $\bar L$ in \eqref{intergen} would be described by a $2 \times 2$ matrix.

Since the nucleation time is ``long''
(in the simulation the time needed for the supercritical droplet
to invade the whole box
is short with respect to the time needed
for its appearance),
solving \eqref{intergen} with little overlapping $\nu_{\bar x}$
would lead to a very small $\bar L$.
A natural alternative approach would be to provide,
with the same kind of link $\Lambda$,
an intertwining relation
between Markovian kernels (rather than generators) of the form:
\begin{equation}\label{intermeta} \Lambda K_{q'} = \bar P \Lambda,\end{equation}
with $K_{q'}(x, \cdot)$ as in~\eqref{kernel2}. The measure
$K_{q'}(x, \cdot)$ is indeed the distribution of the original process,
with generator $L$, started at a configuration $x\in \VV$ and looked
at an exponential time $T_{q'}$ with parameter $q'$.
This parameter should be
of the same order as the nucleation rate.

Unfortunately for this specific Ising example, we do not know how to write down
solutions of such intertwining equations with $\Lambda$ having small squeezing.
However, while the main results in metastability studies are usually 
written in some asymptotic regime
(e.g. for the Glauber dynamics illustrated in Figure~\ref{pianta}: $h \ll 1$ with an associated diverging box side length),
we stress that {\bf {\em looking at metastable dynamics
through intertwining equations makes possible
to deal with such dynamics outside of any asymptotic regime}}.
This is made clear through the following statement which we state in discrete-time since~\eqref{intermeta} is as in~\eqref{inter} with $P=K_{q'}$.

\begin{theo}[{\bf Metastability through intertwining}]\label{poireau}
Let $\hat X$ be a Markov chain with state space $\VV$ and transition kernel $P$. 
Assume Equation~\eqref{inter} holds for some $(\Lambda, \bar P)$. 
Then for each $\bar x$ in $\bar{\VV}$, 
there exists a stopping time $T_{\bar x}$ for $\hat X$ and a random variable
$\bar Y_{\bar x}$ with values in $\bar{\VV} \setminus \{\bar x\}$ such that
\begin{itemize}
\item $T_{\bar x}$ is distributed as a geometric variable with parameter $1 - \bar P(\bar x, \bar x)$;
\item $\nu_{\bar x}$ is stationary up to time $T_{\bar x}$, i.e., for all $t \geq 0$,
	\begin{equation}\label{mucca}
		P_{\nu_{\bar x}}\left(\hat X(t) = \cdot \bigm| t < T_{\bar x}\right) = \nu_{\bar x}
		\,;
	\end{equation}
\item $P_{\nu_{\bar x}}\left(\bar Y_{\bar x} = \bar y\right) = \frac{\bar P(\bar x, \bar y)}{1 - \bar P(\bar x, \bar x)}$
	for all $\bar y$ in $\bar{\VV}$;
\item $P_{\nu_{\bar x}}\left(\hat X(T_{\bar x}) = \cdot \bigm| \bar Y_{\bar x} = \bar y\right) = \nu_{\bar y}(\cdot)$;
\item $\left(\bar Y_{\bar x}, \hat X(T_{\bar x})\right)$ and $T_{\bar x}$ are independent.
\end{itemize}
\end{theo}
\noindent 
This statement corresponds to Prop. 6 in~\cite{ACGM1} and it can be seen as a partial rewriting of section~2.4 of~\cite{DF} in the spirit of~\cite{[Mi]}.

\subsubsection{\large{
A coarse-graining algorithm for metastability: from processes to measures}}\label{metagol}
To make use in practice of the result in Theorem~\ref{poireau}, as motivated so far, 
one should find explicit solutions $(\Lambda, \bar P)$ to~\eqref{intermeta} with corresponding $\nu_{\bar x}$'s having small joint overlap, which for relevant non-trivial examples could be too difficult if not unfeasible.    
We introduce here a deterministic algorithm depending on some tuning parameters to circumvent this problem. For a given Markov process $X$ with a big state space (or equivalently the given associated network $\GG$), the goal is to build {\bf {\em a measure-valued process on a small state space ``mimicking'' dynamical aspects of the original process}} in the sense of Theorem~\ref{poireau}. 
We will accomplish our task if the generator of the resulting measure-valued process is close to be intertwined with the original one, and if  the associated $\Lambda$ has small squeezing.
To guarantee these properties, we will then randomize the procedure through the random forests and give appropriate estimates for the tuning of the involved parameters.     

\begin{center}{\bf {\em Deterministic algorithm based on partitioning}}\end{center}
Given an irreducible and reversible network $\GG$ on $n$ vertices:
\begin{itemize}
\item[a.] Pick $m\leq n$ and consider a partition of the graph $\PP(\GG)=[\AA_1,\ldots, \AA_m]$ into $m$ blocks.
\item[b.] Set $\bar\VV
:=\{1,\ldots,m\}$ for the {\em new vertex set};
\item[c.] Set the {\em duality linking matrix} $\Lambda$ as follows: 
 for any $\bx \in \bar \VV$, 
 \begin{equation}\label{metanu}
 \Lambda(\bx, \cdot)=\nu_{\bx}(\cdot) := \mu_{\AA_{\bx}}(\cdot),
 \end{equation} 
 with $\mu_\AA$ being the probability $\mu$ conditioned to $\AA\subset \VV$, 
 i.e. $\mu_\AA = \mu(\cdot | \AA)$;
 \item[d.] The {\em new process} is given by  
 \begin{equation}
 \label{metaPbar}
 \bP_{q'}(\bx,\by):=  P_{\nu_{\bx}}\cro{X(T_{q'}) \in \AA(\by)},
 \end{equation} 
 for any $\bx,\by \in \bar{\VV}$ and with $T_{q'}$ being
an independent exponential random variable of parameter $q'>0$.
\end{itemize}

A few remarks: 

\begin{itemize}
\item Proposition~\ref{flatness.prop} guarantees that the linking measures in~\eqref{metanu} have minimal  squeezing  $\SS(\Lambda)$ equal 
 to one. 
 \item Provided that this pair $(\Lambda,{\bar P}_{q'})$ is close to be a solution to the intertwining Equation~\eqref{intermeta}, the network $\bar \GG=\bar \GG(q')$ identified by ${\bar P}_{q'}$ can be seen as a reduced measure-valued description of the original network $\GG$ on time-scale $T_{q'}$, hence a possible answer to our original problem {\bf (P1)}. In particular, by construction, inherited from $X$,  ${\bar P}_{q'}$ is again irreducible and reversible w.r.t. $\mu$.
\item For any choice of the parameter $q'>0$, in view of step d. and the irreducibility assumption, the resulting network $\bar \GG$ will be given by a complete graph with non-homogenous weights identified by~\eqref{metaPbar}. 
\end{itemize}

It remains to clarify how to choose the initial partition $\PP(\GG)$ in step a. above 
and how to guarantee that $(\Lambda,{\bar P}_{q'})$ is an approximate solution 
to the intertwining Equation~\eqref{intermeta}. 
To this end, we can exploit the nice properties of the random forests in~\eqref{nu} and 
randomize the deterministic algorithm above. 

\centerline{{\bf {\em Randomization through forests:}} }

Set
\begin{equation}\label{randomize} m=m(q):=|\RR(\Phi_{q})|\leq n \quad \text{           and        } \quad {\mathcal P}(\GG):={\mathcal P}(\Phi_{q})=[\AA_1,\ldots, \AA_{m(q)}]\end{equation} 
for $q>0$.
The choice of the random partition ${\mathcal P}(\Phi_{q})$ is motivated by Theorems~\ref{maude} and~\ref{maschera} and the fact that by Wilson's algorithms, two nodes 
$x,y\in\VV$ tend to be in the same block of ${\mathcal P}(\Phi_{q})$ 
if on scale $T_{q'}$, the process $X$ is likely to walk from $x$ to $y$ or vice-versa.
The following theorem quantifies in which sense the pair $(\Lambda, {\bar P}_{q'})$ is an 
approximate solution to~\eqref{intermeta} and can serve as a guideline to tune the involved 
parameters $(q,q')$. Since for each $\bx\in\bV$, $\Lambda K_{q'} (\bx,\cdot)$ and 
${\bar P}_{q'} \Lambda (\bx,\cdot)$ are probability measures on $\VV$, we will use the total 
variation distance defined in~\eqref{TV}.

\begin{figure}
    \centerline{
        \includegraphics[width=3.75cm]{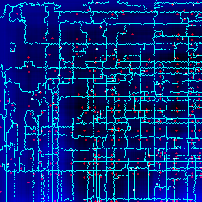}
        \hfill
        \includegraphics[width=3.75cm]{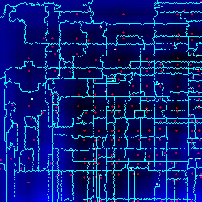}
        \hfill
        \includegraphics[width=3.75cm]{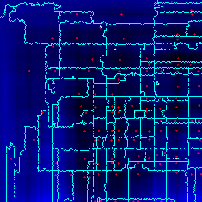}
        \hfill
        \includegraphics[width=3.75cm]{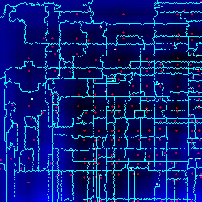}
    }
    \bigskip
    \centerline{
        \includegraphics[width=3.75cm]{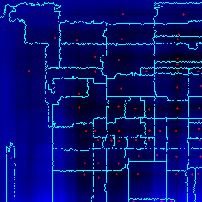}
        \hfill
        \includegraphics[width=3.75cm]{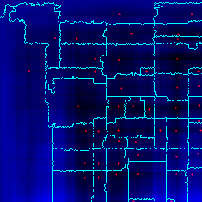}
        \hfill
        \includegraphics[width=3.75cm]{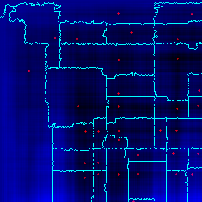}
        \hfill
        \includegraphics[width=3.75cm]{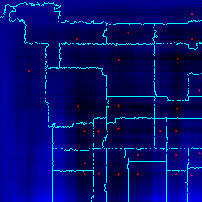}
    }
    \caption{
        Coupled partitions ${\cal P}(\Phi_q)$
        for a Metropolis random walk in Brownian sheet potential
        at inverse temperature $\beta = 0.16$ and in a $512 \times 512$
        square grid. The cyan vertices are the leaves of the trees,
        the roots are at the center of red diamonds
        and the other vertices are colored according to their potential:
        the darker the blue, the deeper the potential.
        The values of $q$ are
        $1.22 \times 10^{-4}$,
        $3.05 \times 10^{-5}$,
        $7.63 \times 10^{-6}$,
        $1.91 \times 10^{-6}$,
        $4.77 \times 10^{-7}$,
        $1.19 \times 10^{-7}$,
        $2.98 \times 10^{-8}$ and
        $7.45 \times 10^{-9}$.
    }\label{chocolat}
\end{figure}

\begin{theo}[{\bf Control on intertwining error for metastability}]
\label{TVmeta.theo} Let $p \geq 1$, and $p^*$ its conjugate exponent, i.e., $\frac 1 p + \frac 1{p^*} = 1$. Fix positive parameters $(q,q')$, and let 
$(\Lambda, {\bar P}_{q'})$ be the pair given by~\eqref{metanu} and~\eqref{metaPbar} 
randomized through $\Phi_q$ as in~\eqref{randomize}. Then,   
\[ \E\cro{\sum_{\bx =1}^{|\RR(\Phi_{q})|} d_{TV}\big(\Lambda K_{q'}(\bx, \cdot), \bP_{q'} \Lambda(\bx, \cdot)\big)}
\leq \pare{\E\cro{\abs{\RR(\Phi_{q})}}}^{1/p}
\pare{\frac{q'}{q} \sum_{x \in \VV} \E\cro{|\Gamma^x _{q'}|}}^{1/p^*} 
,
\] 
where 
$|\Gamma^x _{q'}|$ in the r.h.s. denotes the length (i.e. the number of crossed edges) of the trajectory of a loop-erased random walk
on the original graph started from $x\in\VV$ and stopped at an exponential time $T_{q'}$ (recall Proposition~\ref{march}). 
\end{theo}
The proof of the above statement together with some insights on how to tune $(q,q')$ to guarantee the bound in the r.h.s. to be small, can be found in~\cite{ACGM1}, cf. Thm. 4 and related discussion therein. For the interested reader, Thm. 5 in~\cite{ACGM1} also says how to modify the intertwining matrix $\Lambda$ in~\eqref{metanu}, and make it into an exact solution to~\eqref{intermeta} for small enough $q'$. 
Figure~\ref{chocolat}
illustrates the coupled partitions
associated with different values of $q$
for a random walk in random potential.
We consider a Metropolis nearest-neighbour random walk
associated with Brownian sheet potential $V$
and inverse temperature $\beta$
on the square box $[0, 511]^2 \cap \mathbbm{Z}^2$.
This means that the rates $w(x,y)$
are given by
$w(x, y) = \exp\bigl\{-\beta[V(y) - V(x)]_+\bigr\}$
if $x$ and $y$ are nearest neighbours, and 0 if not.
In this picture, the vertical and horizontal axes
are oriented southward and eastward respectively,
so that $V$ is 0 on the left and top boundaries.
Since the value of $V$ at the bottom-right corner
is of order 500 (it is the sum of $510^2$ independent
normal random variables),
we already have a metastable situation 
for the case $\beta = 0.16$, illustrated in Figure~\ref{chocolat}:
on the time scales $1 / q$ corresponding to the different partitions
${\cal P}(\Phi_q)$, the random walk tends to be trapped
in each piece of the partition.

Let us finally observe that if such a coarse-graining scheme
can be practically implemented on very large networks,
the typical size of statistical mechanics configuration spaces
is too large for these algorithms to be run in such cases. 
For the Ising example in section~\ref{Ising},
we would need to deal with a network made of $|\VV|=2^{65536}$ nodes.
But it could be that a similar coarse-graining procedure 
applied on the graph on which the spin are defined
(the $256 \times 256$ two-dimensional torus in our case),
rather than on the configuration graph (with its $2^{65536}$ nodes
in our case), can be used to derive a coarse description
of the original dynamics. We believe that such an approach would deserve some investigation.
In section~\ref{multigol} we give another similar algorithm to address these slightly different questions.
This other algorithm will be concretely motivated by applications in signal processing rather than in metastability.
After describing it, we will turn back to our random walk on the Brownian sheet.

\subsection{\large{
Network coarse-graining from a signal processing point of view}}\label{multigol}
We present here a different network reduction scheme.
Our original motivation will be fully explained in the next section~\ref{app3} where we construct a pyramidal algorithm to \emph{process} signals (i.e. real valued functions) defined on the vertex set of a network.
In first approximation, the idea of these pyramidal algorithms is to build  from a given network and associated signal, a multiscale invertible reduction scheme where at each scale, after so-called \emph{downsampling} and \emph{filtering} procedures, one can define a network of smaller size (typically having a fraction of the original nodes) with associated a coarse grained approximation of the signal. We again have to address problem {\bf (P1)} from the beginning of this section, but the purpose is different. We will follow a similar strategy as for metastability, but looking at solutions $(\Lambda, \bar L)$ to Equation~\eqref{intergen} under a different perspective. In this case, as intertwining rectangular matrix of size $|\bV| \times |\VV|$, we use  
 \begin{equation}\label{lamq} 
 \Lambda=\Lambda_{q'}:=[K_{q'}]_{\bV, \VV}, \quad \text{  with parameter  } q'>0,
 \end{equation}
that is, the restriction of the kernel in~\eqref{kernel2}.   
Actually, based on the concrete problem under investigation, other suitable kernels can also be used, this is just a convenient choice for our application explained in more details in section~\ref{app3}.
We only anticipate that as for metastability, the measures identified by the $\Lambda$'s rows should concentrate on different regions of the underlying state space. They will play the role of \emph{low frequency filters} in the wavelets construction to build coarse-grained versions of the original signal through local averages. Anyhow, in this section, we focus only on the network reduction step, that is, on the $\bar L$  we want to consider, and how we can  guarantee that it is almost intertwined with the given $L$.

\subsubsection{\large{Graph reduction via roots subsampling and the trace process}}\label{multigol2}

The reduction algorithm presented here makes use of the notion of Schur complement of a matrix which we first recall.  
Let $M$ be a  matrix of size $n=p+r$ and let 
\[ M=  \begin{pmatrix} A & B \\ C & D \end{pmatrix} \]
be its block decomposition, $A$ being a square matrix of size $p$ and $D$ a square 
matrix of size $r$. If $D$ is
invertible, {\bf the Schur complement of $D$ in $M$} is the square matrix of size
 $p$ defined by
\[ S_M(D) := A - B D^{-1} C  \, .
\] 

\begin{center}{\bf {\em Deterministic algorithm based on representative vertices}}\end{center}

Given an irreducible and reversible network $\GG$ on $n$ vertices:

\begin{itemize}
\item[a.] Let $\bar\VV\subset \VV$ be a subset of $m=|\bV| \leq n$ {\em selected vertices} in $\GG$. 
\item[b.] For the {\em reduced process}, set 
\begin{equation}
\label{qtrace} \bar L \text{ to be the Schur complement of } [L]_{\brV} \text{ in } L
\end{equation} 
with $\brV=\VV\setminus\bar{\VV}$. 
\item[c.] The {\em coarse-grained network} $\bar \GG$ is the graph 
with vertex set $\bar\VV$, weights 
$\bar w(\bx,\by):= \bar L(\bx,\by)$ and edge set identified by positive weights among vertices.
\end{itemize}

Let us stress the main features: 

\begin{itemize}
\item The Schur complement is a classical electrical-network-like reduction, sometimes 
referred to as Kron's reduction, having nice properties and a clear probabilistic interpretation 
which we explain in Proposition~\ref{SchurProp} below. In particular, as for the algorithm in 
section~\ref{metagol}, the resulting process is still  irreducible and reversible, thus allowing for 
successive iterations of this coarse graining procedure (see example in Figure~\ref{nonsparse}).   
\item The resulting network $\bar \GG$ tends to be (depending on local bottlenecks and the 
locations of the selected vertices in the original graph) a complete graph with non-
homogenous weights, cf. Equation~\eqref{trace}. Depending on the specific problem, one 
would wish to obtain a sparser graph (e.g. when iterating the scheme, dealing with sparse 
matrices can reduce the algorithmic complexity). To this aim, a possibility is to redistribute part
 of the masses in~\eqref{trace} and disconnect edges with low weights. An example of such a 
 ``sparsification'' procedure will be briefly mentioned in section~\ref{minne} below. 
\item Unlike the algorithm presented in section~\ref{metagol}, for any $q'>0$, the measures 
identified by $\Lambda_{q'}$ do not have disjoint supports. Thus quantitative bounds on their 
squeezing are desirable.
\item As in the case of metastability, the pair $(\Lambda_{q'}, {\bar L})$ identified by~
\eqref{lamq} and~\eqref{qtrace} is not a solution to~\eqref{intergen}, but can be turned into an 
approximate solution.
\end{itemize}

In the following proposition, we collect relevant properties of the Schur complement and its probabilistic interpretation. Its proof can be found in~\cite{ACGM2}, see Lemma 13 in there. 

 \begin{prop}[{\bf Schur complement and trace process}]\label{SchurProp} 
 \label{LbarSchur.lem} Let $\bar \VV$ be any subset of $\VV$ and $L$ as in~\eqref{gen} reversible w.r.t. $\mu$. Set $\brV=\VV\setminus\bar{\VV}$.
 Then the restricted matrix $[L]_{\brV}$ is invertible, and the Schur complement $\bar{L}$ of $[L]_{\brV}$ in $L$ is an irreducible Markov generator on $\bV$ reversible w.r.t. $\mu$. This Markov process with generator $\bar{L}$ is often referred to as {\bf trace process} on the set $\bar\VV$.
 Further, the discrete-time kernel given by  $\bP = \Id + \bLL/w_{max}$ admits the following interpretation:
 \begin{equation}
 \label{trace}
 \bP(\bx,\by)=P_{\bx}\cro{\hat X(T_{\bV}^+)=\by}, \quad \text{ for all  }\bx,\by \in \bV, 
 \end{equation}
with $T_{\bV}^+$ being {\bf the first return time} of $X$ in $\bV$.
 \end{prop}

To sum up,  we link any $\bx, \by \in \bV$  with a weight (possibly zero) proportional to~\eqref{trace}, that is, the probability that the original walk $X$ starting from $\bx$ lands in $\by$ when hitting again the subset $\bV\subset \VV$.
We can next move to control squeezing and intertwining error for this choice of $(\Lambda_{q'}, {\bar L})$. In view of the well distributed property of the forests roots in Theorem~\ref{maude}, a natural way to select the reduced vertex set in step a. is to consider 
\begin{equation}\label{random2} \bV=\bV(q):= \RR(\Phi_q), \quad \text { for } q>0,\end{equation}
which, as for metastability, makes the pair $(\Lambda_{q'}, {\bar L})$ {\bf {\em randomized}} w.r.t. $\Phi_q$.        
We next move to control squeezing and intertwining error in this case for which we use the $\ell_{p}$-norm in~\eqref{pnorm}: 

\begin{theo}[{\bf Squeezing vs Intertwining}]
\label{TVmulti} Fix $q'>0$. Consider $\bV\subset \VV$, and set $\brV = \VV \setminus \bV$. Then 
the deterministic pair $(\Lambda_{q'}, {\bar L})$ given by~\eqref{lamq} and~\eqref{qtrace}  satisfies:
 \begin{equation}
 \label{intertwiningp.ineq}
 \nor{\pare{\bar{L} \Lambda_{q'} - \Lambda_{q'} L} f}_{p,\bV} 
 \leq 2 q'  \pare{\frac{w_{max} }{\beta}}^{1/p^*} \frac{1}{\mu(\bV)^{1/p}} \nor{f}_{p,\VV}
  \, 
  \end{equation}
for any $p \geq 1$, $p^*$ being its conjugate exponent, and any $f$ in $\ell_p(\VV,\mu)$, 
where \begin{equation}\label{beta}
\frac{1}{\beta} = \max_{\bx \in \bV} \sum_{z \in \VV} P(\bx, z) E_{z}\cro{T_{\bV}} \text{   and } P  \text{ as in }\eqref{skeleton}.
 \end{equation}

Further, by randomizing $(\Lambda_{q'}, {\bar L})$ as in Equation~\eqref{random2} for some $q>0$,  we get the squeezing control:
\begin{equation*}
	\E\cro{ \SS\left(\Lambda_{q'}\right) \bigm | |\RR(\Phi_q)|=m}
	\leq \frac{
		\min\left\{
			\sqrt{1 + \sqrt{\frac{T_n}{S_n}}} \; \exp\pare{\sqrt{S_n T_n}- V_n} ;
			\sqrt{1 + T_n} \exp\pare{\frac{\left(1 + S_n T_n)\right)}{2} - V_n}
		\right\}
	}{\P\cro{|\RR(\Phi_q)| = m}}
\end{equation*}
for any $m \leq n$ , with 
\[ S_n = \sum_{j=1}^{n-1} p_j(q')^2 \left(1-p_j(q)\right)^2 \, \, , \, \, 
T_n= \sum_{j=1}^{n-1} \frac{p_j(q)^2}{p_j(q')^2 } \, \, , \, \,  
V_n =  \sum_{j=1}^{n-1} p_j(q) \left(1-p_j(q)\right) 
\,, 
\]
and $p_j(\cdot)$ as in~\eqref{pj} with $\BB=\emptyset$.
 \end{theo}
 
This result corresponds to part of Theorem 3 in~\cite{ACGM1} and Proposition 16 in~\cite{ACGM2}.
Note that the upper bound on the squeezing depends on $L$ through its spectrum only. 
For general remarks and insights on how to tune $(q,q')$ based on this result, we refer the reader to the comment after Theorem 3 in~\cite{ACGM1} and to 
section 6 in~\cite{ACGM2}. In particular, section 6.2 gives estimates on $\beta$ defined in~\eqref{beta}. 
We conclude this part on coarse-graining algorithms by giving a concrete example on a classical road network benchmark.

\subsubsection{\large{An experiment: reduction of Minnesota road network }}\label{minne}

\begin{figure}[H]
\includegraphics[height=6cm,width=16cm]{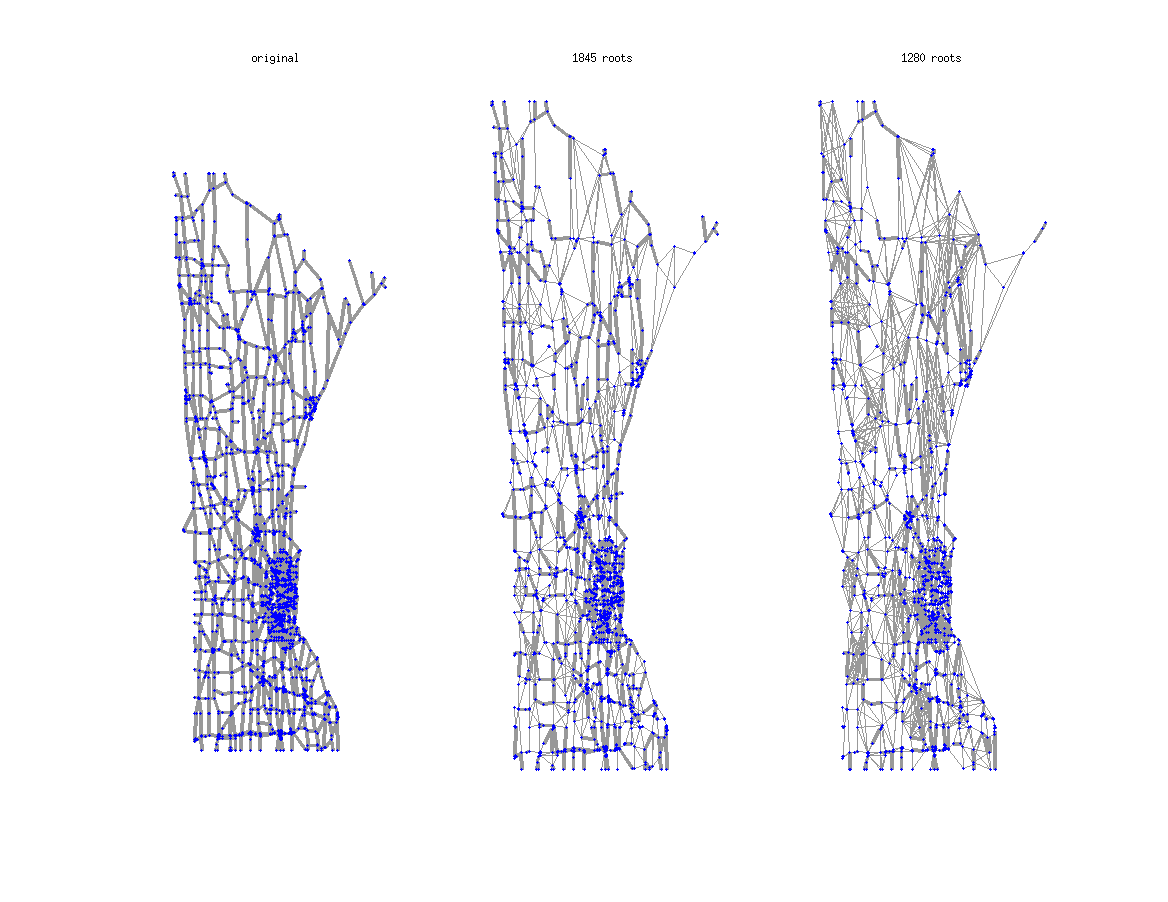}
\caption{Two successive reductions of the Minnesota's graph using the algorithm
 in section~\ref{multigol}. The wider the edge, the larger its weight.}
 \label{nonsparse}
\end{figure}

\begin{figure}[H] 
   \includegraphics[height=6cm,width=16cm]{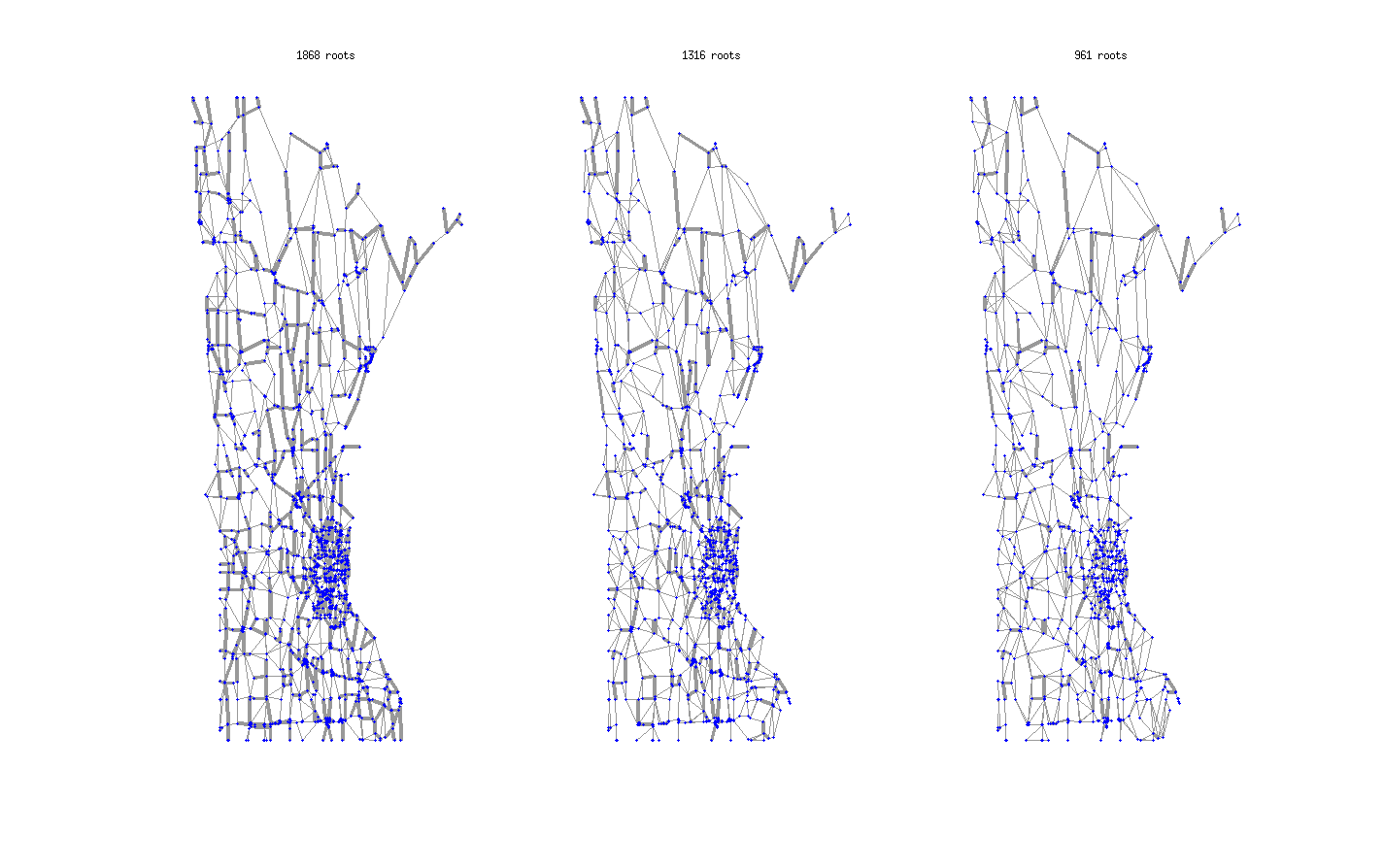}
\caption{Two successive reductions of the Minnesota's graph using the algorithm in section~\ref{multigol} with sparsification.}
 \label{sparse}
\end{figure}

Figures~\ref{nonsparse} and~\ref{sparse} illustrate 
our recursive coarse-graining algorithm
with and without sparsification procedure.
We start from the Minnesota road network
with unitary weights and choose our downsampling parameter
as explained in section~\ref{app3}.

Figure~\ref{nonsparse} shows that, as mentioned previously,
our original sparse graph is rapidly changed into a dense one
after a few Kron's reductions.
The intertwining error,
i.e., the approximation error in the intertwining equation,
is then our guideline to build a sparsification procedure.
We set down to 0 a family of small weights $\bar w(\bar x, \bar y)$
obtained by computing the Schur complement $\bar L$,
to get a sparser $\bar L_s$.
These are chosen in such a way that for each $\bar x$
the new error 
$\|\nu_{\bar x} - \bar L_s\Lambda(\bar x, \cdot)\|_\infty$
does not exceed the original error
$\|\nu_{\bar x} - \bar L\Lambda(\bar x, \cdot)\|_\infty$
by more than a fraction of it.
We refer to section~7.2 in~\cite{ACGM2} for the details.
Figure~\ref{sparse} illustrates the result. 

\subsubsection{Back to the random walk in random potential in Figure~\ref{chocolat}}
\begin{figure}
    \includegraphics[width=5.4cm]{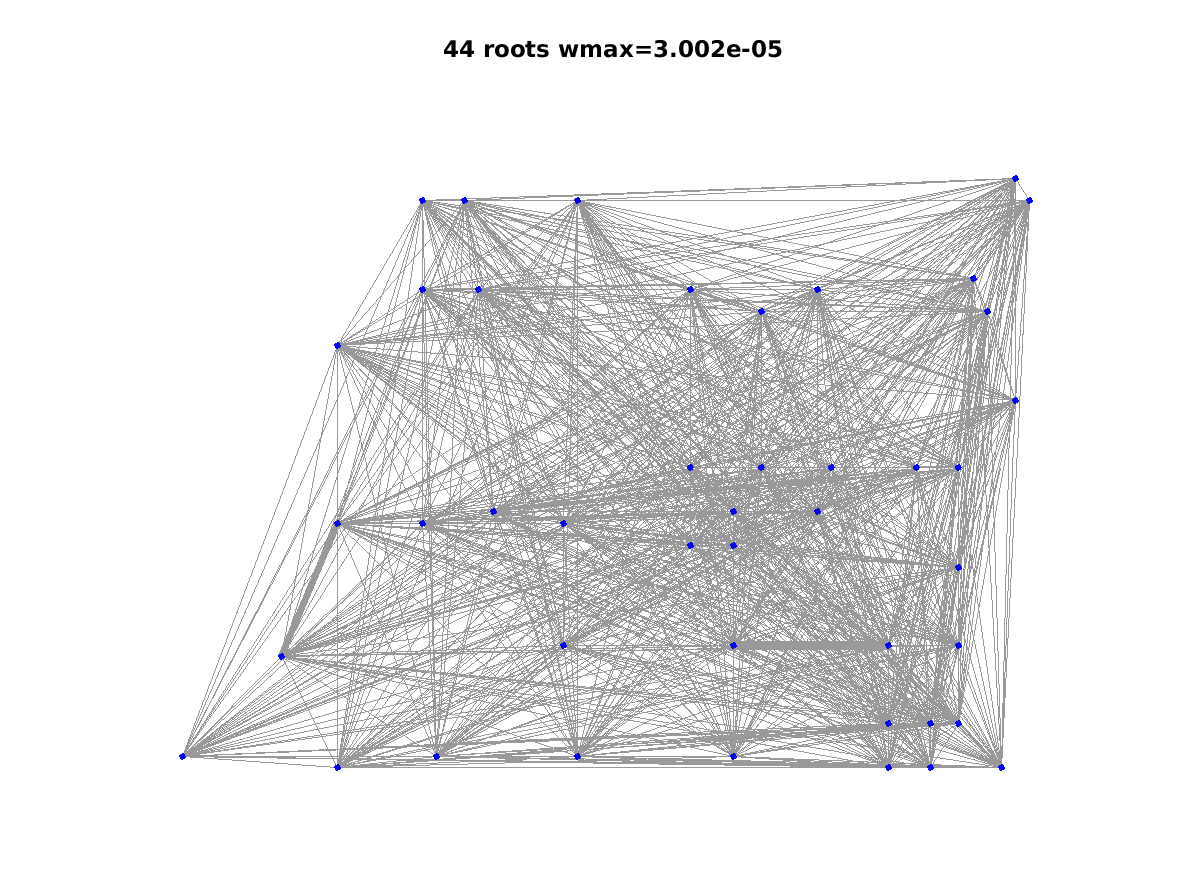}
    \includegraphics[width=5.4cm]{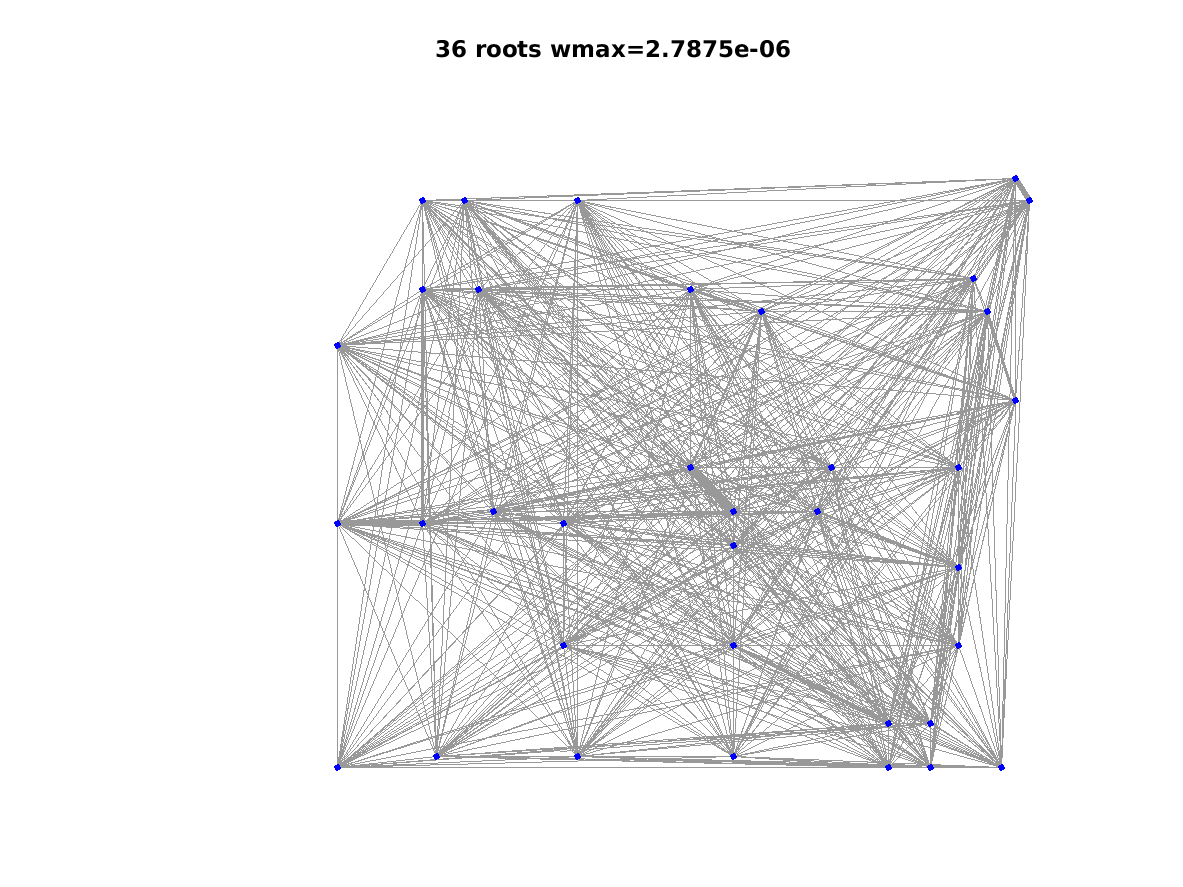}
    \includegraphics[width=5.4cm]{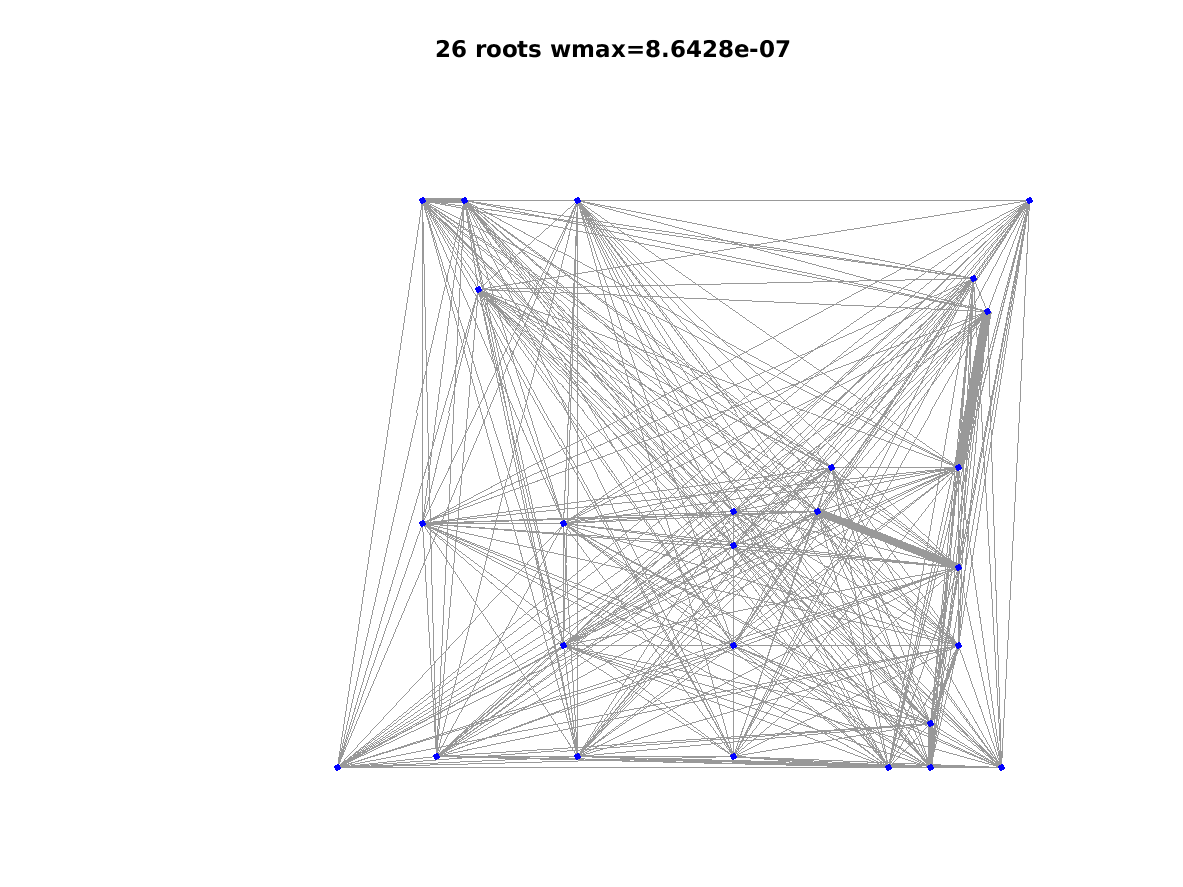}
    \includegraphics[width=5.4cm]{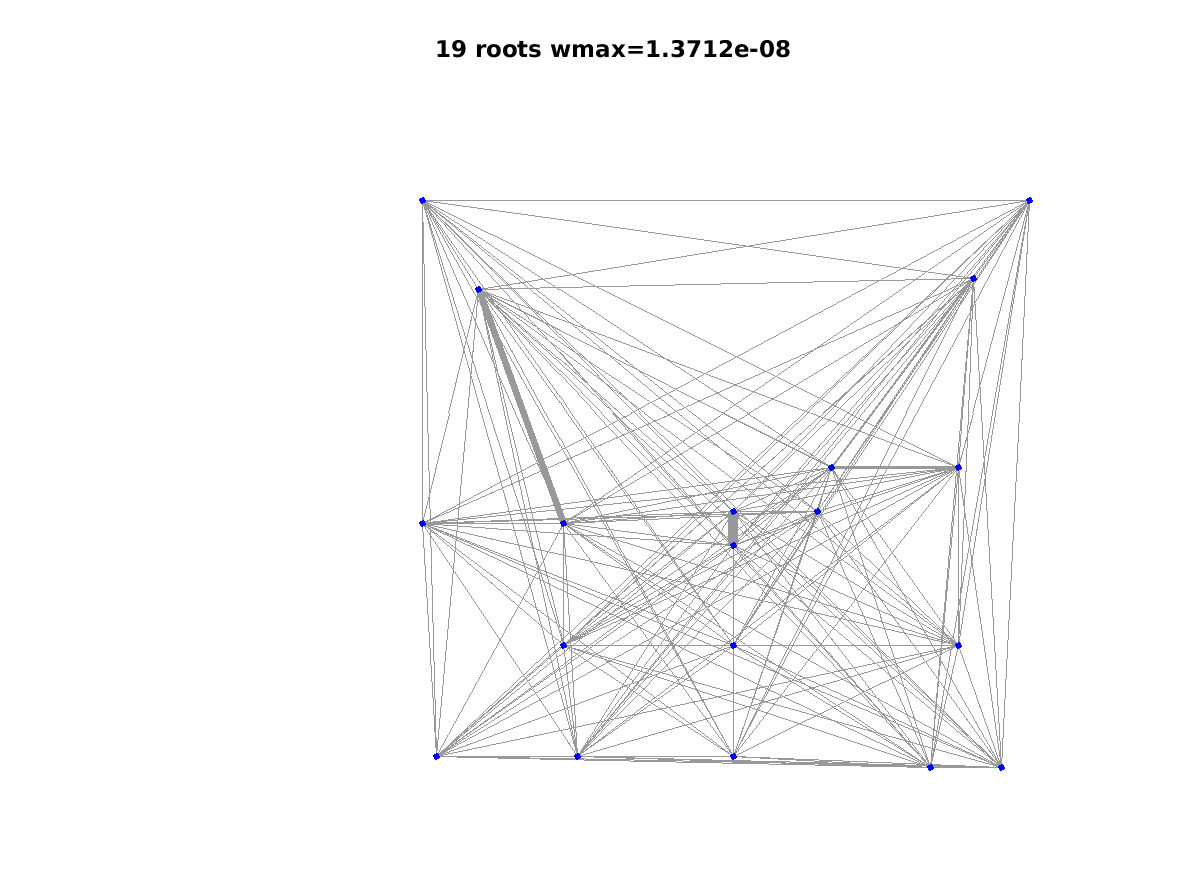}
    \includegraphics[width=5.4cm]{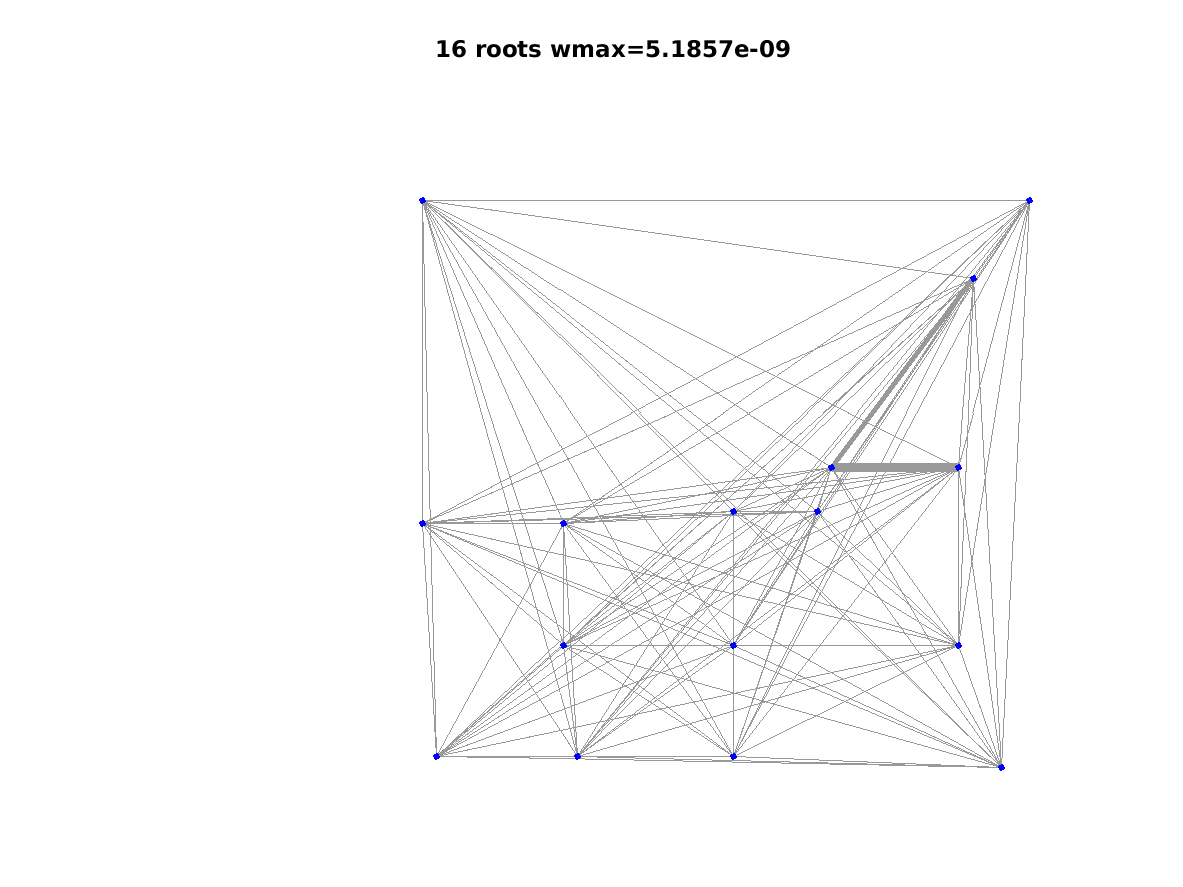}
    \includegraphics[width=5.4cm]{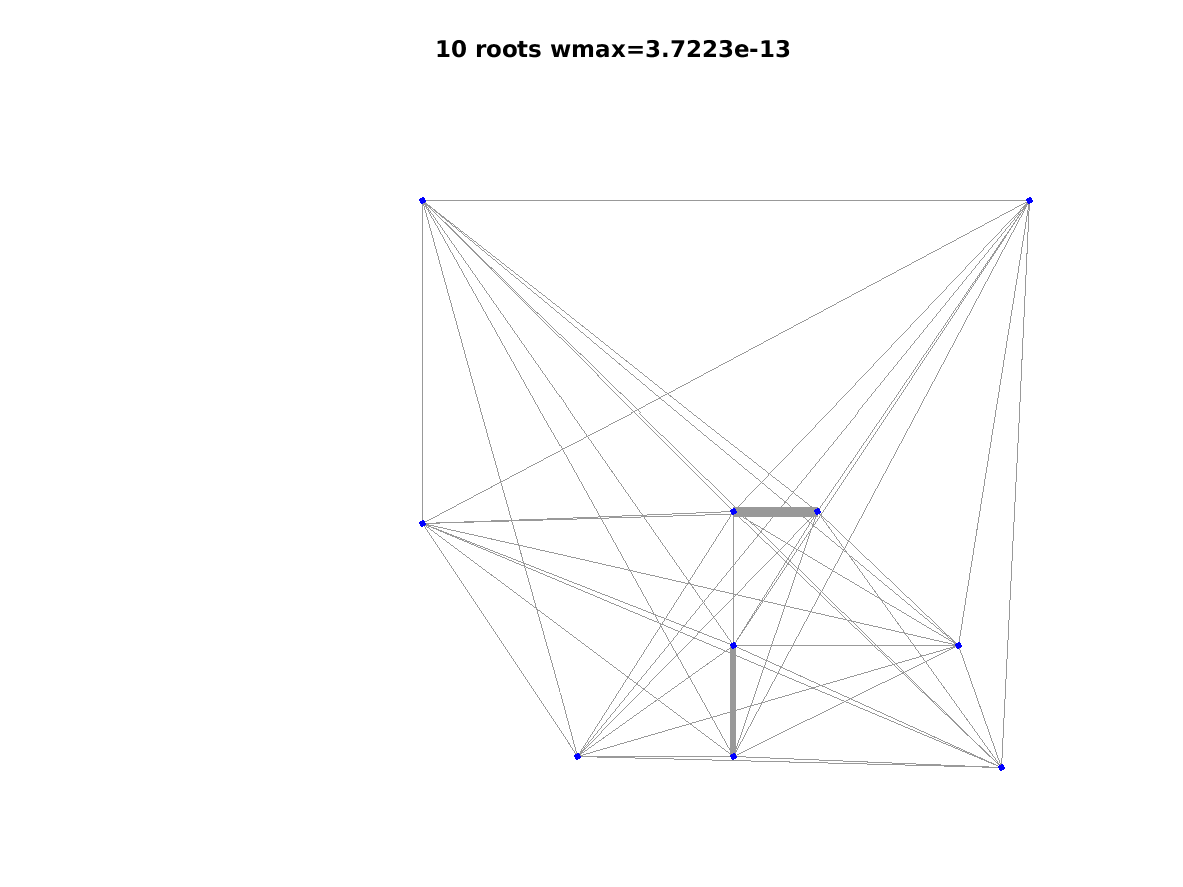}
    \includegraphics[width=5.4cm]{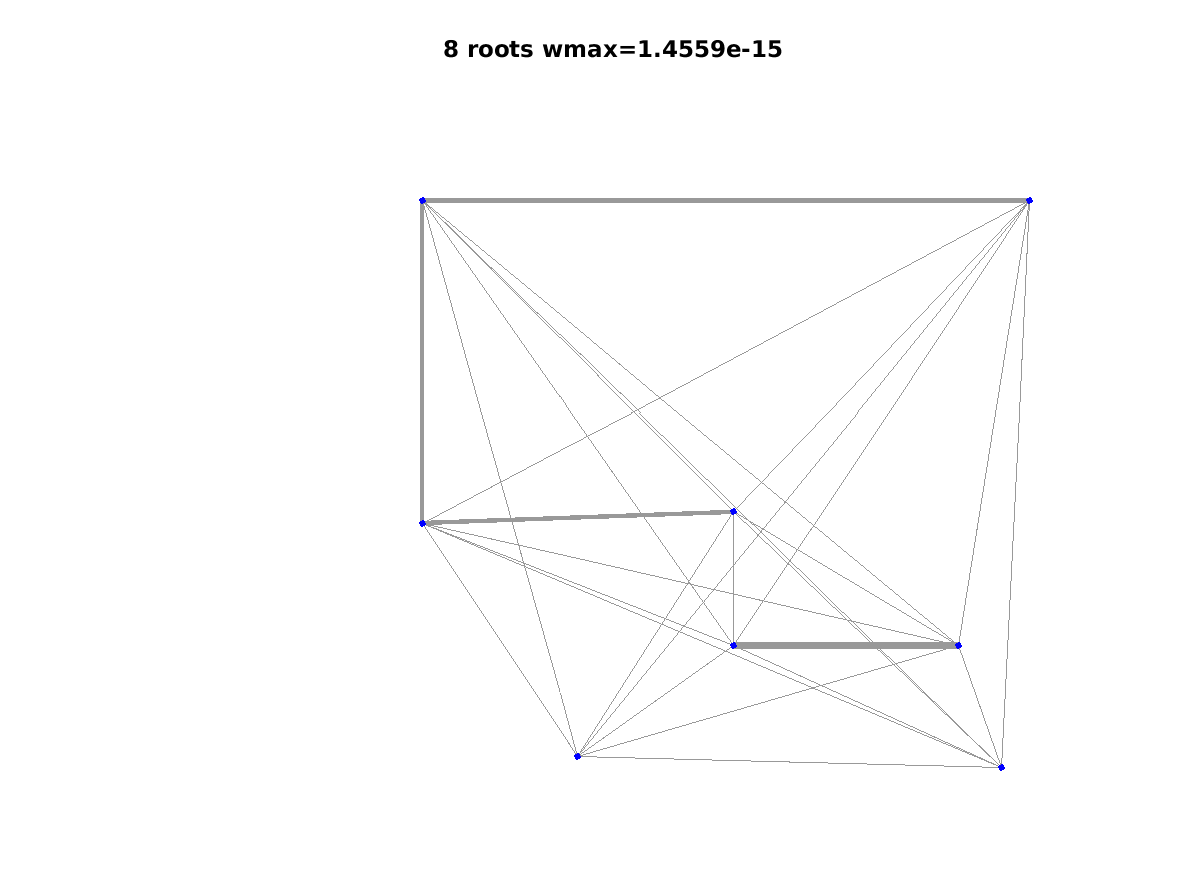}
    \includegraphics[width=5.4cm]{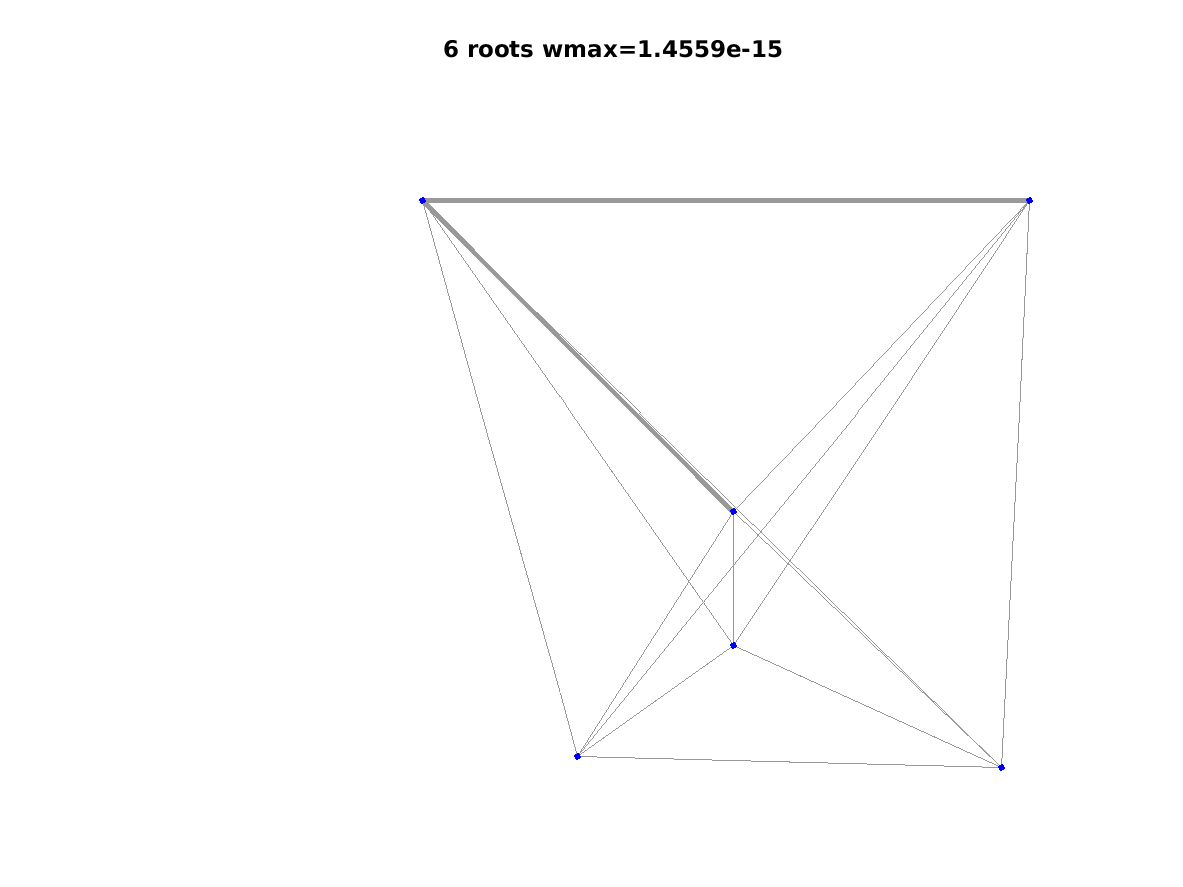}
    \includegraphics[width=5.4cm]{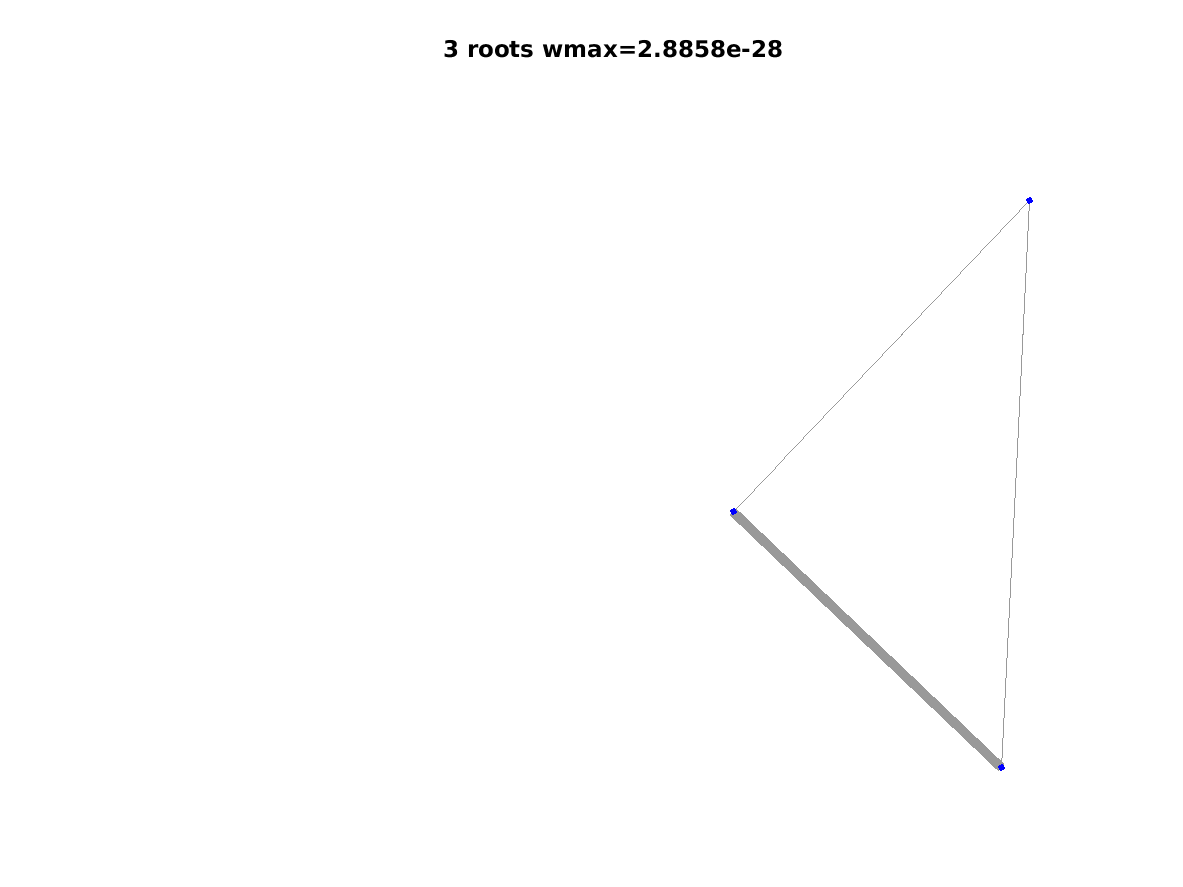}
    \includegraphics[width=5.4cm]{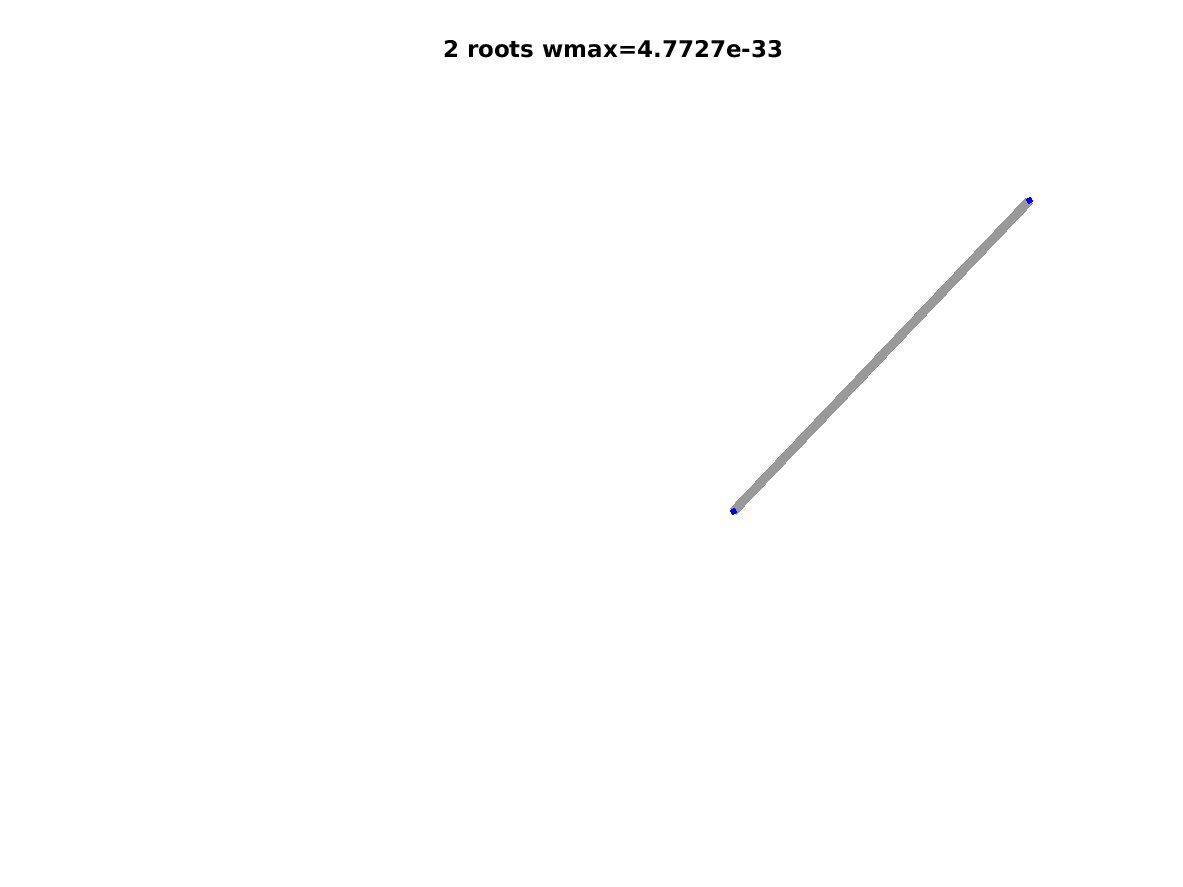}
    \caption{
        Successive reduced networks for a random walk
        in Brownian sheet potential.
        The vertex number as well as the largest jumping rate
        are reported on each picture.
    }\label{clo}
\end{figure}

The very nature of metastable dynamics
makes impossible in practice
to run Wilson's algorithm to get a partition ${\cal P}(\Phi_q)$
with a few pieces only, even for relatively small networks. 
The parameter $q$ should be so small that we would get
a huge running time of the algorithm.
But the recursive procedure we described,
with at each step a downsampling parameter $q$
of the same order as the current maximal jumping rate $w_{max}$,
makes it possible to describe the long time dynamics through the computation
of the trace process on a very small set of points.
Figure~\ref{clo} shows some of these reduced graphs 
for the random walk in Brownian potential on a smaller grid
($64 \times 64$) but with larger inverse temperature
$\beta = 2.56$, so that we face the same kind of difficulty.
This reduction allows us to describe the dynamics up to convergence
to equilibrium, 
which occurs on a time scale of order $10^{33}$.
This is similar in spirit to the renormalization
procedure introduced in~\cite{Sco}.
The main difference here
is that each vertex of the reduced graphs
is associated through the linking matrices
with a probability measure on the original graph. 
Our reduced process is a measure-valued process
rather than a Markov process on some local minima
of an energy landscape.

\section{\large{Applications: Intertwining Wavelets, multiresolution for graph signals}}\label{app3}
Weighted graphs provide a flexible representation of geometric structures of irregular domains, and they are now a commonly used tool to encode and analyze data sets in numerous 
disciplines including neurosciences, social sciences, biology,  transport, communications... 
Edges between vertices represents interactions between them, and 
 weights on edges quantify the strength of the interaction. In data modeling, it is often the case
  that a {\bf graph signal } comes along with the network structure $\GG$, this is simply a real-
  valued function on the vertex set of the network: $$f:\VV\rightarrow \R.$$
Functional magnetic resonance images measuring brain activity in distinct functional regions 
are classical example of such a graph signal. 
Signal processing is the discipline devoted to develop tools and theory to process and analyze signals. Depending on concrete instances, {\bf {\em processing means e.g. classifying, removing noise, compressing or visualizing.}}
In case of signals on regular domains very robust tools and algorithms have been developed over more than half a century mainly based on Fourier analysis. 
For signals on arbitrary networks studies are less advanced and only in recent years significant efforts from different communities have been dedicated to develop suitable efficient methods. We refer  to~\cite{SCHU} for a recent review on this growing investigation line.  
We present here a new and rather general {\bf {\em multiresolution scheme}} we introduced in~\cite{ACGM2}. 
Multiresolution scheme is a generic name for several multiscale algorithms allowing to decompose and process signals.
We will start with a quick recap of classical multiresolution schemes on regular domains. In section~\ref{multiforest} we then present our new forest-algorithmic-scheme and in section~\ref{interwave} we collect the main theorems providing a solid theoretical framework to our method and guidelines for the practice. Illustrative 
numerical experiments of various nature are given in the concluding section~\ref{action}. 

\subsection{\large{Classical multiresolution: wavelets and pyramidal algorithms on regular grids}}\label{classic} 
Let us consider a discrete periodic function  $f:\Z_n=\Z/n\Z\rightarrow \R$, viewed as a vector in $\R^n$.  
The multiresolution analysis of $f$  is based on {\em wavelet analysis} which roughly   
amounts to compute an ``approximation'' $\baf \in \R^{n/2}$ and a ``detail'' component $\brf \in\R^{n/2}$ through classical operations in signal processing such as ``filtering'' and ``downsampling''. 
The idea is that {\bf {\em the approximation gives the main trends}} present in $f$ whereas {\bf {\em the detail contains more refined information}}. This is done by splitting the frequency content of $f$ into two components:  the approximation $\baf$ focuses on the low frequency part of $f$ whereas the high frequencies in $f$ are contained in $\brf$. 
\begin{itemize}\item
{\bf {\em Filtering}} is the operation allowing to perform such frequency splittings and it consists of computing a convolution $f\star k$ for some well chosen kernel $k$: 
$K_{l}(f)=f \star k_{l}$ yields $K_{l}(f)$ as a low frequency component of $f$ and 
$K_{h}(f)=f \star k_{h}$ yields $K_{h}(f)$ as a high frequency version of $f$.
\item {\bf {\em Downsampling:}} The  vectors $\baf$ and $\brf$ are ``downsampled'' versions of $K_l(f)$ and $K_h(f)$ by a factor of 2, which means that one keeps one coordinate of $K_l(f)$ and  $K_h(f)$ out of two, to build $\baf$ and $\brf$ respectively. Thus the total length of the concatenation of the two vectors $[\baf,\brf] $ is exactly $n$, hence the length of $f$. 
\end{itemize}

To sum up we have
\begin{eqnarray*}
&\baf(\overline{x})=K_l(f)(\overline{x})=<\varphi_{\overline{x}},f>= \text{ Approximation }\\
& = \text{ downsampled low-frequency components of } f,
\end{eqnarray*}
and
\begin{eqnarray*}
&\brf(\overline{x})=K_h(f)(\overline{x})=<\psi_{\overline{x}},f>= \text{ Detail }\\
& = \text{ downsampled high-frequency components of } f,
\end{eqnarray*}

where 

\begin{itemize}
\item $\overline{x}$ belongs to the set of downsamples $\overline{\Z_n} $ isomorphic to $\Z/\frac{n}{2}\Z $.
\item $\{\varphi_{\overline{x}},\overline{x}\in \overline{\Z_n}\}$ is the set of functions such that the equality between linear forms $<\varphi_{\overline{x}},.>=K_l(\cdot)(\overline{x}) $ holds for all $\overline{x}\in\overline{\Z_n}$.
\item In the same way, $\{\psi_{\overline{x}},\overline{x}\in \overline{\Z_n}\}$ is such that  $<\psi_{\overline{x}},.>=K_h(\cdot)(\overline{x}) $ holds for all $\overline{x}\in\overline{\Z_n}$.
\end{itemize}

The choice of $k_l$ and $k_h$ is clearly crucial and done in a way that perfect reconstruction of $f$ from $\baf$ and $\brf$ is possible with no loss of information in the representation $[\baf,\brf]$. 
By denoting $f_0=f$, $f_1=\baf$ and $g_1=\brf$, we see that this splitting scheme can be successively iterated starting from $f_1$ to obtained a sequence $f_{N}\in\R^{n/2^{N}}, g_{N}\in\R^{n/2^{N}},\ldots,g_{1}\in\R^{n/2} $, for any integer $N$ where the total length of the concatenated vectors $[f_{N},g_{N},\ldots,g_{1}]$ is exactly $n$. 
This leads to \begin{center}{\bf the  multiresolution scheme:}\end{center}

\begin{equation}\label{multischeme} \begin{matrix} f_0 = f & \rightarrow & f_1=\baf & \rightarrow & f_2 & \cdots &  \rightarrow & f_N 
\\ & \searrow &  & \searrow &   & & \searrow & 
\\[-.1cm]
&                 & g_1=\brf &         &g_2 &              &  & g_N 
\end{matrix}
\end{equation}
with $f_i, g_i \in \R^{n/2^{i}}, i=0 \ldots, N.$

We remark that the perfect reconstruction condition amounts to have $\{\varphi_{\overline{x}},\psi_{\overline{x}},\overline{x}\in\overline{\Z_n}\} $ a basis for the signals 
$f$ on $\Z_n$. A famous construction by Ingrid Daubechies~\cite{DAUB} derives several families of orthonormal compactly supported such basis. 
It is worth mentioning that these families combine {\bf {\em localization in space }} around the point $\overline{x}$ and  {\bf {\em localization properties in frequency }} due to the filtering step they have been built from. 
Using this space-frequency localization one can derive key properties of the wavelet analysis of a signal which rely on the deep links between the local regularity properties of $f$ and the behavior and decay properties of detail coefficients. We refer the interested reader to one of the numerous books on wavelet methods and their applications such as~\cite{DAUB} or~\cite{MAL}, in all these methods the word $wavelets$ denotes the family $\{\psi_{\overline{x}},\overline{x}\in \overline{\Z_n}\}$ spanning for the high-frequencies components.

\subsection{\large{Forest-multiresolution-scheme on arbitrary networks}}\label{multiforest} 
When considering signals $f$ on irregular networks $\GG$, it is not clear how to reproduce the classical multiresolution scheme described above.
In other words, in a non-regular network, there are no canonical neither obvious answers to the following questions:

%
%

\begin{enumerate}
\item[{\bf (Q1)}] What kind of downsampling should one use? What is the meaning of ``every other point''?
\item[{\bf (Q2)}] On which weighted graph should the approximation $f_1$ be defined to iterate the procedure?
\item[{\bf (Q3)}] Which kind of filters should one use? What is a good notion of low frequency components of a signal?
\end{enumerate}  

In light of the properties and applications of the random spanning forests described in the previous sections, we do have natural answers to the first two questions. 
In fact, Theorem~\ref{maschera} suggests that the set of roots is a random downsample tailored to any network providing a possible answer to {\bf (Q1)}\footnote{This proposal has already received some attention within the signal processing community~\cite{Trembley}}. And the network coarse-graining algorithm presented in 
section~\ref{multigol} is a good candidate to address {\bf (Q2)}.
It remains to make sense of filtering in {\bf (Q3)}, that is, how to capture the low frequencies components or the main trend of a graph signal $f$.   
But if we use the algorithm in 
section~\ref{multigol}, then it is natural to define the approximation component $\baf(\bx)$ as a local average around the downsampled $\bx\in\bV$
w.r.t. the measures identified by the intertwining matrix $\Lambda=\Lambda_{q}$ in~\eqref{lamq}.
That is, for each downsampled $\bx\in\bV$:  
 \[ \baf(\bx) := \nu_{\bx}(f) = \sum_{x \in \VV} \Lambda(\bx,x) f(x) \,.
  \]
Before proceeding with our construction, let us give some remarks on one of the main problems in defining good filters in signal processing.

\subsubsection{\large{Filtering : fighting against Heisenberg.}}\label{filters} 
To be of any practical use in signal processing, the filters $\{\nu_{\bx}, \bx \in \bV\}$ have to be 
well localized both in space and frequency. This is violating Heisenberg uncertainty principle, a delicate problem which needs a proper compromise. 
In the graph context, frequency localization means that the filters belong to an eigenspace of the graph laplacian $L$. Hence, 
we are interested in solutions to Equation~\eqref{intergen} such that the measures  $\{\nu_{\bx}, \bx \in \bV\}$ (our proposed filters) 
are linearly independent measures tending to be non-overlapping (space localization), and contained in eigenspaces of $L$.  
We already observed that  saying that $(\Lambda, \bLL)$ is an exact solution to~\eqref{intergen} implies that the linear space spanned by the measures $\{\nu_{\bx}, \bx \in \bV\}$ 
is stable by $L$, and is therefore a direct sum of eigenspaces of $L$, so that these measures provide filters which are frequency localized. 
Hence {\bf {\em the error in the intertwining relation is a measure of frequency localization:}} the smaller the intertwining error, the better the frequency localization. 
And through Theorem~\ref{TVmulti} we can control such frequency localization. 

Concerning space localization, we then want small squeezing on $\Lambda$. Notice further that our $\Lambda$ in~\eqref{lamq}  
is just the restriction of the Green's kernel $K_{q'}$ which is very sensitive to localization by tuning the parameter $q'$. In fact:   
\begin{itemize}
\item when $q'$ goes to $0$,  for any $\bx \in \bV$, 
$K_{q'}(\bx,y)$ goes to $\mu(y)$ so that~\eqref{intergen} is trivially satisfied. Since $\mu$ is   
the left-eigenvector of $L$ corresponding to the eigenvalue $0$, the $K_{q'}(\bx,y)$ are well 
localized in frequency. However, the vectors
$\{K_{q'}(\bx, \cdot) \mid \bx \in \bV\}$ become linearly dependent and very badly localized in space. 
\item On the other extreme, when $q'$ goes to $\infty$, $K_{q'}(\bx, \cdot)$ goes to $\delta_{\bx}$. Hence, the space localization is perfect. However, the frequency localization is lost, and  the error in~\eqref{intergen} tends to grow.
\end{itemize}
A compromise has to be made for the choice of $q'$. Since in our method there is also the 
parameter $q$ controlling the fraction of downsampled points, we will need a suitable joint 
tuning of the pair $(q,q')$ to optimize localization properties. As explained in section~\ref{tuneqs}, our tuning choice is indeed 
guided by Theorem~\ref{TVmulti} but also on stability results of the proposed method which we state in section~\ref{interwave} below. For a detail discussion on the actual choice of the parameters, we refer to section 6 in~\cite{ACGM2}.

\subsubsection{\large{Approximation, detail and the full forest-multiresolution}}
Let us summarize our proposed forest-multiresolution-scheme and present the corresponding basis construction. For arbitrary real-functions $f,g$ on $\VV$, we will denote by 
\begin{equation}\label{muscalar}\left\langle f,g \right\rangle_{\mu}= \sum_{x \in \VV} f(x) g(x) \mu(x) \end{equation}
the {\bf scalar product w.r.t. $\mu$}.

\begin{center}{\bf {\em Intertwining Wavelets }}\end{center}

Given an irreducible and reversible network $\GG=(\VV, w)$ on $n$ vertices and a signal $f: \VV \rightarrow \R$.
\begin{itemize}
\item[a.] {\bf {\em Forest-downsampling:}} Choose a fix $q>0$, let $\bar{\GG}=\bar{\GG}(q)=(\bV,\bw)$ be the randomized coarse-grained (irreducible and reversible) network given by the algorithm in section~\ref{multigol2} with $\bV=\bV(q)$ as in~\eqref{random2}, and set $\brV =\VV\setminus\bV$.

\item[b.] {\bf {\em Forest-filtering:}} Fix $q'>0$ and let $\Lambda=\Lambda_{q'}$ be as in~\eqref{lamq}. Define the {\bf approximation component of $f$} as 
the function $\baf$ defined on $\bV$ by 
\begin{equation}
\label{barf.def}
  \baf(\bx) := \Lambda f(\bx)= K_{q'}f(\bx) \, , \quad  \forall \bx \in \bV \, , 
\end{equation}
and the {\bf detail component of $f$} as the function  $\brf$ defined on $\brV$ by 
\begin{equation}
\label{brevef.def}
 \brf(\brx) := (K_{q'}-\Id_{\brV})f(\brx) \, , \quad \forall \brx \in \brV \, .  
\end{equation}
\end{itemize}

\begin{theo}[{\bf Basis and wavelets}]\label{basis} Fix a parameter $q'>0$. For each $\bx \in \bV\subset \VV$ and $\brx \in \brV,$ respectively, define on $\VV$ the densities functions of 
 the measures $\Lambda(\bx,\cdot)$ w.r.t. $\mu$:
\begin{equation}
\label{nononde}
\phi_{\bx}(\cdot) := \frac{\Lambda(\bx,\cdot)}{\mu(\cdot)}=\frac{ K_{q'}(\bx,\cdot)}{\mu(\cdot)}  \, \, ,
\end{equation}
and the functions
\begin{equation}
\label{onde}
\psi_{\brx}(\cdot) := \frac{ (K_{q'}-\Id_{\brV})(\brx, \cdot)}{ \mu(\cdot)} \, ,
\end{equation}
and abbreviate 
$ \{\xi_x\}_{x \in \VV}= \{ (\phi_{\bx},\psi_{\brx}) \mid \, \bx \in \bV, \brx \in \brV \}.$
Then the family $\{\xi_x\}_{x \in \VV}$ is a basis of $\ell_2(\VV,\mu)$.
In particular,  for any  $\brx \in \brV$,  $ \bra{\psi_{\brx},\ind} = 0$.
\end{theo}
The statement above corresponds to Lemma 9 in~\cite{ACGM2}. As in classical multiresolution analysis, the functions $\{\psi_{\brx}, \brx \in \brV\}$ represent our {\bf wavelets}. 
The basis functions given by $\{\xi_x\}_{x \in \VV}$ are not pairwise orthogonal w.r.t. $\bra{\cdot,\cdot}$ but, by considering the corresponding {\bf dual basis in $\ell_2(\VV,\mu)$}, that is, the family  $\{\tilde{\xi}_x\}_{x \in \VV}$ defined through $$\bra{\tilde{\xi}_x,\xi_y} = \delta_{xy}, \quad x,y \in\VV, $$ for any $f\in\ell_2(\VV,\mu)$, we get the following  representation
\begin{equation}\label{expansion} 
f = \sum_{x \in \VV} \bra{\xi_x,f} \tilde{\xi}_x=\sum_{\bx \in \bV} \bra{\phi_{\bx},f} \tilde{\xi}_{\bx} + \sum_{\brx \in \brV} \bra{\psi_{\brx},f} \tilde{\xi}_{\brx}\,  
\end{equation}
identifying our ``split into low and high frequency'' components. We call {\bf analysis operator } $U=U_{q'}$ the operator 
 \begin{equation}\label{anoperator} U: \begin{array}[t]{lcl} \ell_2(\VV,\mu) & \rightarrow & \R^{\VV}= \R^{\bV} \times \R^{\brV}
 	\\ f & \mapsto & \{\bra{\xi_x,f}\}_{x\in \VV}=\{\bra{\phi_{\bx},f}, \bra{\psi_{\brx},f}\}_{\bx \in \bV, \brx \in \brV}=[\baf,\brf]
	\end{array} \,  
	\end{equation}
assigning to $f\in\ell_2(\VV,\mu)$ its coefficients in the expansion in~\eqref{expansion}. As explained in the following theorem, the reconstruction of $f$ from the knowledge of its coefficients $U(f)$ can be made operative:

\begin{theo}[{\bf Reconstruction formula}]
\label{reconstruction}
Fix $q'>0$. For any $f \in \ell_2(\VV,\mu)$, consider $\baf \in \ell_2(\bV,\mu)$ and $\brf \in \ell_2(\brV,\mu)$ respectively given by 
\[ \baf(\bx) = K_{q'}f(\bx)=U(f)(\bx) \, , \quad  \bx \in \bV \, , \]
\[ \brf(\brx) =( K_{q'}-\Id_{\brV})f(\brx)=U(f)(\brx) \, , \quad \brx \in \brV. 
\]

Then, 
\[ f = \bR_{q'} \baf + \brR_{q'} \brf  \, ,  
\]
where
\begin{equation}\label{apoperator}
 \bR_{q'} = \begin{pmatrix} \Id_{\bV} - \frac{1}{q'} \bLL 
 \\[.2cm]
  [-L]_{\brV}^{-1} [L]_{\brV}
\end{pmatrix}  = \text{ {\bf approximation  operator ,}}    
\end{equation} 
and 
\begin{equation}\label{deoperator} 
\brR_{q'} =\begin{pmatrix}  [L]_{\bV} [-L]_{\brV}^{-1} \\  - \Id_{\brV} - q' [-L]_{\brV}^{-1} 
\end{pmatrix} = \mbox{ {\bf detail  operator}}
 \, . 
\end{equation}
\end{theo}

This last statement correspond to Prop. 10 in~\cite{ACGM2} and fully describes our multiresolution scheme. In fact, in view of the properties of $\bar{\GG}$, we can simply iterate the procedure with $(\bV,\bLL,\mu_{\bV})$ in place of $(\VV, L,\mu)$ resulting into a pyramidal algorithm as described in Equation~\eqref{multischeme}.

Still, we have not enough motivated the choice of the filter bank given in~\eqref{barf.def} and~\eqref{brevef.def} for our ``high and low frequencies'' components in relation to regularity properties of  signals. To clarify this point, let us first emphasize that 
in any reasonable wavelets construction, as mentioned in section~\ref{filters}, good space and 
frequency localization properties are necessary ingredients. And in our approach, this issue is partially addressed via Theorem~\ref{TVmulti}. Notice in particular that space localization is achieved using the determinantality of $\bV$, and the fact that $\bV$ is well spread on 
$\VV$. Since $\brV= \VV \setminus \bV$ is also a determinantal process, with kernel $ \Id-K_q$, this suggested the detail definition~\eqref{brevef.def} (the sign convention for the $\psi_{\breve x}$ is chosen to have a self-adjoint analysis operator $U: f \mapsto [\baf, \brf]$ in $\ell_2(\VV,\mu)$). 
On the other hand, another fundamental ingredient is that if the signal $f$ is ``regular enough'', then the corresponding detail component $\brf$ should be small. For instance, if the original function is constant it should not contain any high frequency component, that is, the corresponding $\brf$ being identically zero. As the last remark in Theorem~\ref{basis} states, this is in fact the case in our intertwining wavelets. However, more generally, a way to capture and guarantee that the size of the details is small for non--constant but regular functions is desirable. In the next section we give bounds on the involved operators as a function of $q'$ and we make sense of this latter regularity issue through the norm of $Lf$ in Theorems~\ref{brf.prop} and~\ref{jack}.

\subsection{\large{Quality and stability of intertwining wavelets}}\label{interwave} 
We collect here our results to guarantee numerical stability when using the forest-multiresolution-scheme and to guide in the choice of the parameters.  
We stress that the scheme presented in the previous section can be implemented when the downsample $\bV\subset\VV$ is chosen according to any (even deterministic) rule. 
For this reason, the following statements controlling norms and sizes of the main involved objects are stated for arbitrary $\bV\subset\VV$ and will not depend on $q$.
We start by the control on the approximation and detail operators in~\eqref{apoperator} and~\eqref{apoperator}, respectively.

\begin{theo}[{\bf Bound on the norm of the approximation operator}] 
  \label{approxoperator}
Fix $q'>0$.  Let $\bV$ be any subset of $\VV$, $\brV = \VV \setminus \bV$, and let 
  $\bR_{q'}$ be the operator defined in~\eqref{apoperator}. 
  For any $p \geq 1$,  for any $f \in \ell_{p}(\bV, \mu_{\bV})$, 
  \begin{equation}
  \label{bR-norm.eq}
  \nor{\bR_{q'} f}_{p,\VV}
   \leq  
\cro{\pare{1+ 2 \frac{\bar{w}_{max}}{q'}}^p + \frac{w_{max}}{\beta}}^{1/p} \,  
\mu(\bV)^{1/p}  \nor{ f}_{p,\bV} 
  \,  , 
  \end{equation}
  with $w_{max}$  as in~\eqref{wmax}, $\bar{w}_{max}$ defined analogously w.r.t. the generator $\bLL$, and $\beta$ as in~\eqref{beta}.
 \end{theo}

 \begin{theo}[{\bf Bound on the norm of the detail operator}] 
  \label{detailoperator}
Fix $q'>0$.  Let $\bV$ be any subset of $\VV$, $\brV = \VV \setminus \bV$, and let 
  $\brR_{q'}$ be the operator defined in~\eqref{deoperator}. 
 For all $p \geq 1$,  for all $f \in \ell_{p}(\brV, \mu_{\brV})$, 
  \begin{equation}
    \label{brR-norm.eq}
 \nor{\brR_{q'} f}_{p,\VV}
   \leq  
\cro{\pare{\frac{w_{max}}{\beta}}^{p/p^*}
      + \pare{1+ \frac{q'}{\gamma}}^p}^{1/p} \,  
\mu(\brV)^{1/p}  \nor{ f}_{p,\brV} 
  \,  , 
  \end{equation}
   where  $p^*$ is the conjugate exponent of $p$, $\beta$ as in~\eqref{beta}, and  
   \begin{equation}\label{gamma2}\frac{1}{\gamma} = \max_{\brx \in \brV} E_{\brx}\cro{T_{\bV}}.\end{equation}  
   \end{theo} 
We refer to section 6.2 in~\cite{ACGM2} for estimates on $\gamma, \beta$ and $\bar{w}_{max}$ 
defined in~\eqref{gamma2}. 
We now move to the regularity issue mentioned above, that is, for ``regular'' signal we wish small details. By measuring the modulus of continuity of the original function through $\nor{L f}_{p, \VV}$, we get the following statement corresponding to Prop. 15 in~\cite{ACGM2}:   

\begin{theo}[{\bf Control on the size of details}]
\label{brf.prop}
For any $p \geq 1$ and any $f \in \ell_{p}(\VV,\mu)$, 
\[  \nor{\brf}_{p,\brV} =  \nor{(K_{q'}-\Id) f}_{p,\brV}
 \leq \frac{\max_{x\in\VV} K_{q'}(x,\brV)^{1/p}}{q' \mu(\brV)^{1/p}} \nor{L f}_{p, \VV} \, . 
\]
\end{theo}

The next result gives a control on the size of the coefficients at arbitrary scale $i\leq N$ when   
implementing the multiresolution scheme with $N\geq 1$ successive reductions. To this end, let us introduce suitable abbreviations for the objects at the different scales $i=1,\ldots,N$.
We have a given sequence of $N$ (non--empty) nested vertex sets 
\[ \VV_0 \supsetneq \VV_1 \supsetneq \cdots \supsetneq \VV_N , \quad \text{ starting from } \VV_0= \VV,\] 
with associated parameters $\{q'_i \mid i = 0, \cdots N-1\}$.

For each $i=0,\ldots, N-1$ set: \begin{itemize}
\item $\brV_i = \VV_i \setminus \VV_{i+1}$,
\item $L_0 = L$, and $L_{i+1}$  the Schur complement of  $[L_i]_{\brV_i}$  in $L_i$,
\item the kernels $K_i=q'_i(q'_i \Id_{\VV_i} -L_i)^{-1}$,
\item $w_{i}, \beta_i$ and $\gamma_i$ as in~\eqref{wmax},~\eqref{beta} and~\eqref{gamma2}, respectively, w.r.t. $L_i$. 
\item $f_0=f$, and the succesive approximation and detail components $f_{i+1}=K_i f_i$ and $g_{i+1}=(K_{i} - \Id_{\brV_i}) f_i$ as in~\eqref{barf.def} and  
\eqref{brevef.def}, respectively,  
so that  
\[ f_i = \bR_{i} f_{i+1} + \brR_{i} g_{i+1}  \, ,
\]
with
\[ \bR_i =  \begin{pmatrix} \Id_{\VV_{i+1}} - \frac{1}{q'_i} L_{i+1}
 \\ [- L_i]_{\brV_i}^{-1} [L_i]_{\brV_i \VV_{i+1}}
\end{pmatrix}  
\,  
\,\, 
\mbox{ and }
\quad \brR_i 
=\begin{pmatrix} [L_i]_{\VV_{i+1} \brV_i} [- L_i]_{\brV_i}^{-1} 
\\  q' _i [L_i]_{\brV_i}^{-1} - \Id_{\brV_i}
\end{pmatrix}
 \, . 
\]
\end{itemize}

We can now extend the analysis operator in~\eqref{anoperator} to arbitrary scale $N\geq 1$ by setting  
 
\[ U_N : \begin{array}[t]{lcl} \ell_p(\VV) & \rightarrow & 
	\ell_p(\VV_N,\mu_{\VV_N}) \times \ell_p(\brV_{N-1}, \mu_{\brV_{N-1}}) \times \cdots \times \ell_p(\brV_{0}, \mu_{\brV_{0}})
		\\
		f & \mapsto & [f_N,g_N,g_{N-1},\cdots, g_1] 
		\end{array}
		\]
where the space $\ell_p(\VV_N,\mu_{\VV_N}) \times \ell_p(\brV_{N-1}, \mu_{\brV_{N-1}}) \times \cdots \times \ell_p(\brV_{0}, \mu_{\brV_{0}})$
 is endowed with the norm 
\[ \nor{[f_N,g_N,g_{N-1},\cdots, g_1]}_p = 
\pare{ \mu(\VV_N) \nor{f_N}_{p,\VV_N}^p +  \sum_{i=1}^N \mu(\brV_{i-1}) \nor{g_i}_{p,\brV_{i-1}}^p}^{1/p}
\, . 
\]
Here is the control on the vectors identified by $U_N$ we derived in Prop. 17 in~\cite{ACGM2}.
   
\begin{theo}[{\bf Bound on the norm of the analysis operator}]
Let $p\geq 1$ and $p^*$ be its conjugate exponent. For any $N\geq 1$ and $f \in \ell_p(\VV,\mu)$,
\begin{equation}
\label{analyse.eq}
\nor{U_N(f)}_p \leq 2^{1/p^*} (1+N)^{1/p} \nor{f}_p 
\, . 
\end{equation}
\end{theo}

Our last result is a form of so-called Jackson's inequality. This result is in general crucial for numerical stability in approximation theory and multiresolution analysis (and in particular it plays an important role in our approach to tune the involved parameters, see section 6 in~\cite{ACGM2}). It guarantees small error for ``smooth'' functions, when reconstructing an approximating function of the original one after setting the details $g_i$'s to zero at all scales. This is clearly relevant if e.g. the aim of the multiresolution is to compress a signal. 
To formulate it, notice that by performing $N$ reduction steps in our multiresolution, from the coefficients $[f_N,g_N,g_{N-1},\cdots, g_1]$
we can reconstruct $f=f_0$ as follows:
\begin{align*}
 f=f_0 & = \bR_0 f_1 + \brR_0 g_1
 \\
	  &  = \bR_0 \bR_1 f_2 + \bR_0 \brR_1 g_2 + \brR_0 g_1
\\
& = \bR_0 \bR_1 \cdots \bR_{N-1} f_N
+ \sum_{j=0}^{N-1} (\bR_0 \cdots \bR_{j-1})   \brR_j g_{j+1}
 \, .
\end{align*}
The {\bf approximation of $f$ associated to scale $N$} is thus the function on $\VV$: 
\begin{equation}\label{approx}\tilde{f}(N)=\bR_0 \bR_1 \cdots \bR_{N-1} f_N,\end{equation} and 
we have the following Jackson's inequality measuring its quality.
\begin{theo}[{\bf Jackson's inequality: quality of the approximation of a signal}]
\label{jack}
For any $p \geq 1$ and any $f \in \ell_p(\VV,\mu)$, let $\tilde{f}(N)$ be the approximation associated to scale $N\geq 1$ in~\eqref{approx}. Then
\begin{align}
& \nor{f-\tilde{f}(N)}_{p,\VV}
\nonumber 
\\
& \leq  \acc{\sum_{j=0}^{N-1}  \prod_{i=0}^{j-1} \cro{\pare{1 + 2 \frac{w_{i+1}}{q'_i}}^p + \frac{w_{i}}{\beta_i}}^{1/p}  \cro{\pare{\frac{w_{j}}{\beta_j}}^{p/p^*} 
+ \pare{1 + \frac{q'_j}{\gamma_j}}^p}^{1/p}  \frac{1}{q'_j} } \nor{L f}_{p,\VV}
\nonumber 
\\ 
 \label{jackson.ineq}
 & + \acc{\sum_{j=0}^{N-1}  \prod_{i=0}^{j-1} \cro{\pare{1 + 2 \frac{w_{i+1}}{q'_i}}^p + \frac{w_{i}}{\beta_i}}^{1/p}  \cro{\pare{\frac{w_{j}}{\beta_j}}^{p/p^*} 
+ \pare{1 + \frac{q'_j}{\gamma_j}}^p}^{1/p}  \frac{1}{q'_j}
 \sum_{k=0}^{j-1} 2 q'_{k} \pare{\frac{w_{k}}{\beta_k}}^{1/p^*}} \nor{f}_{p,\VV}
\end{align}
\end{theo}
For the proof of this statement see Prop. 18 in~\cite{ACGM2}. It is worth stressing that the second summand involving $\nor{f}_{p,\VV}$ is due to the propagation of the error in the intertwining relation, it would vanish if the generators at all scales were perfectly intertwined.  

\subsection{Choosing the downsampling and concentration parameters}\label{tuneqs}
Intertwining error and localization property have to be our guidelines
to choose our downsampling parameter $q$ 
and our concentration parameter $q'$.
Beyond squeezing measurement,
another way to look at localization properties
is to focus on the reconstruction operator norm:
spread out wavelets $\phi_{\breve x}$ 
and scaling functions $\psi_{\bar x}$ will lead
to bad reconstruction properties,
i.e., to a large operator norm.
These norms are controlled by Equations \eqref{bR-norm.eq}
and \eqref{brR-norm.eq}.
Since the $\bar R$ operators only are composed together
in the reconstruction scheme, Equation~\eqref{bR-norm.eq}
is the crucial inequality and we look at it
for $p = + \infty$.
As far as the intertwining error is concerned
we look at Equation~\eqref{intertwiningp.ineq},
for $p = +\infty$ too.
These inequalities say that we want $q'/\beta$ 
and $\bar w_{max}/q'$ as small as possible.
This implies that we want $\bar w_{max} / \beta$
as small as possible.
Note that this is a random function of $q$.
It does not depend on $q'$, 
but it is very costly to estimate it by direct Monte-Carlo methods.

Now the same kind of algebra we used to prove that the mean hitting time
of the root set  does not depend on the starting point, offers a workaround.
We can estimate $\bar w_{max}$ with
$$
    \tilde w = \mathbbm{E}\left[
        {1 \over |{\cal R}(\Phi_q)|}
        \sum_{\bar x \in {\cal R}(\Phi_q)} -\bar L(\bar x, \bar x)
    \right]
$$
and $1 / \beta$ with 
$$
    {1 \over \tilde \beta} = \mathbbm{E}\left[
        {1 \over |{\cal R}(\Phi_q)|}
        \sum_{\bar x \in {\cal R}(\Phi_q)}
        \sum_{z \in {\cal V}} P(\bar x, z) E_z[T_{\cal R}(\Phi_q)]
    \right].
$$
It turns out that these expected values 
are respectively equal to 
$$
    \tilde w = q \mathbbm{E}\left[
        {|{\cal V} \setminus {\cal R}(\Phi_q)| \over 1 + |{\cal R}(\Phi_q)|}
    \right]
$$
and 
$$
    {1 \over \tilde \beta} = {1 \over w_{max}} \mathbbm{E}\left[
        {|{\cal V} \setminus {\cal R}(\Phi_q)| \over |{\cal R}(\Phi_q)|}
    \right].
$$
These expected values are easy
to estimate by Monte-Carlo simulations for $q$ between two bounds
of order $w_{max}$ 
---such a restriction is natural if we expect $|\bar {\cal V}|$
to be a fraction of $|{\cal V}|$---
since we have a practical algorithm to sample all the $\Phi_q$ together.
We can then choose $q$ by optimization between the two bounds.
We refer to section~6 in~\cite{ACGM2} for more details.

It remains to choose $q'$ once we have chosen $q$
and sampled $\Phi_q$.
It turns out that the previous estimations on $\bar w_{max}$
and $1 / \beta$ suggest that the norm of the composed $\bar R$
could be bounded by $|{\cal V}|$
by choosing $\bar w_{\max} / q'$ and $q' / \beta$ of the same order.
This is actually ensured by setting
$$
    q' = 2 w_{max}{|{\cal R}(\Phi_q)| \over |{\cal V} \setminus {\cal R}(\Phi_q)|}\,.
$$
While ensuring numerical stability of the algorithm, 
this, in turn, essentially amounts to make
intertwining error and localization property 
of the same importance.
We refer once again to section 6 in~\cite{ACGM2} for details
and we conclude this survey with
giving some numerical experiments based on these results.

\subsection{Experiments, intertwining wavelets in action}
\label{action}
\subsubsection{Comparison with classical wavelet algorithms on the torus}
\begin{figure}
    \includegraphics[width=7cm]{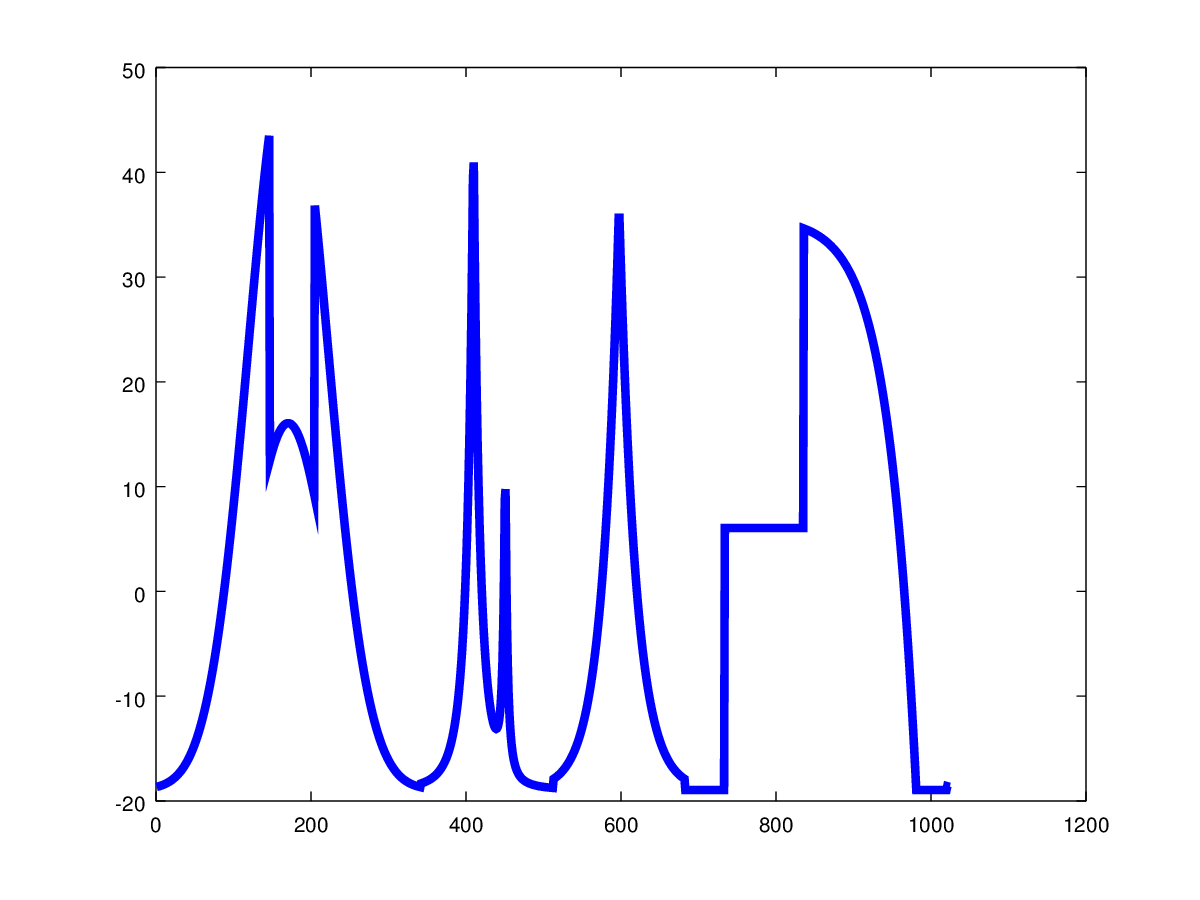}
    \caption{Original signal for analysis on the one dimensional torus.}
    \label{sole}
\end{figure}

\begin{figure}
  \includegraphics[scale=0.6]{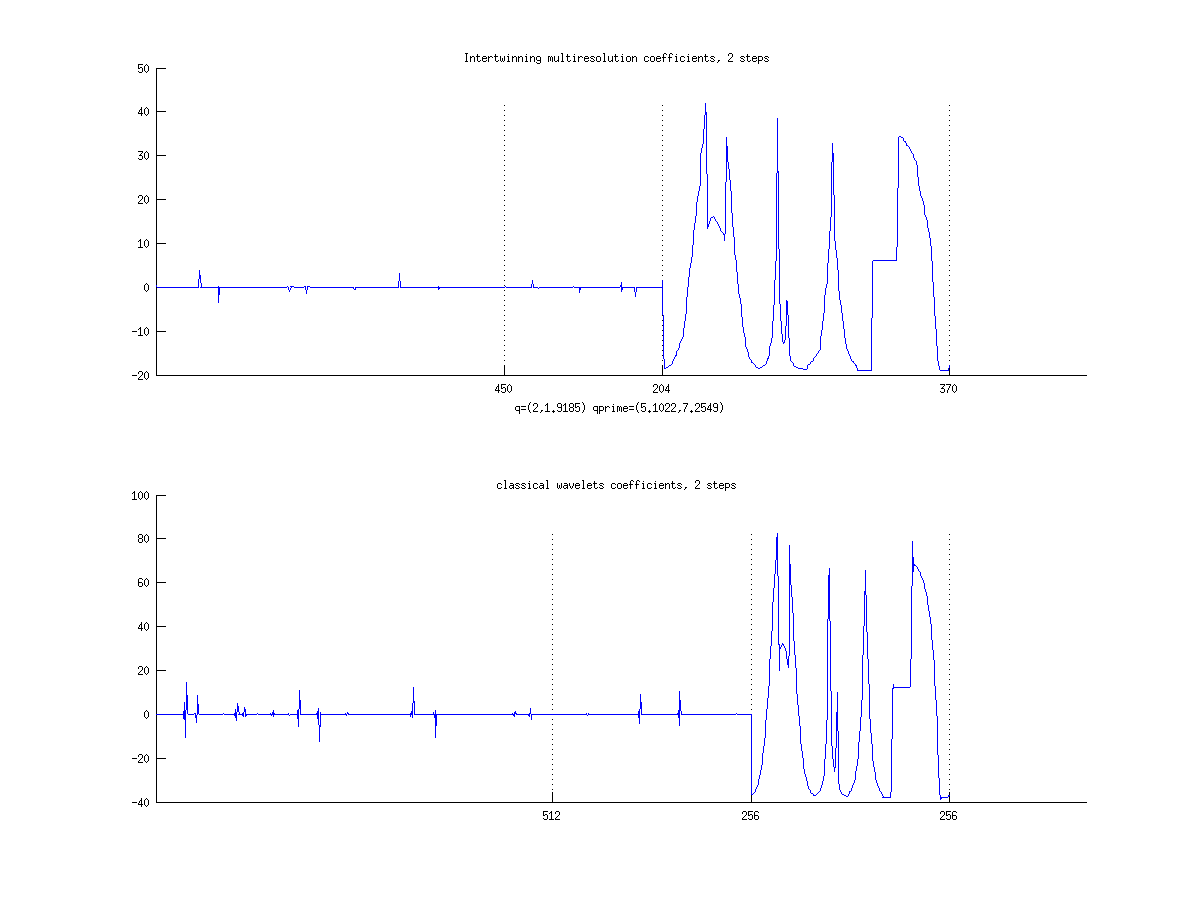}
  \caption{Two steps of the multiresolution for the intertwining wavelets (first graph) and Daubechies-12 wavelets. The three parts of a graph 
  give $(g_1,g_2,f_2)$.}
  \label{alba}
\end{figure}

\begin{figure}
  \includegraphics[width=7cm]{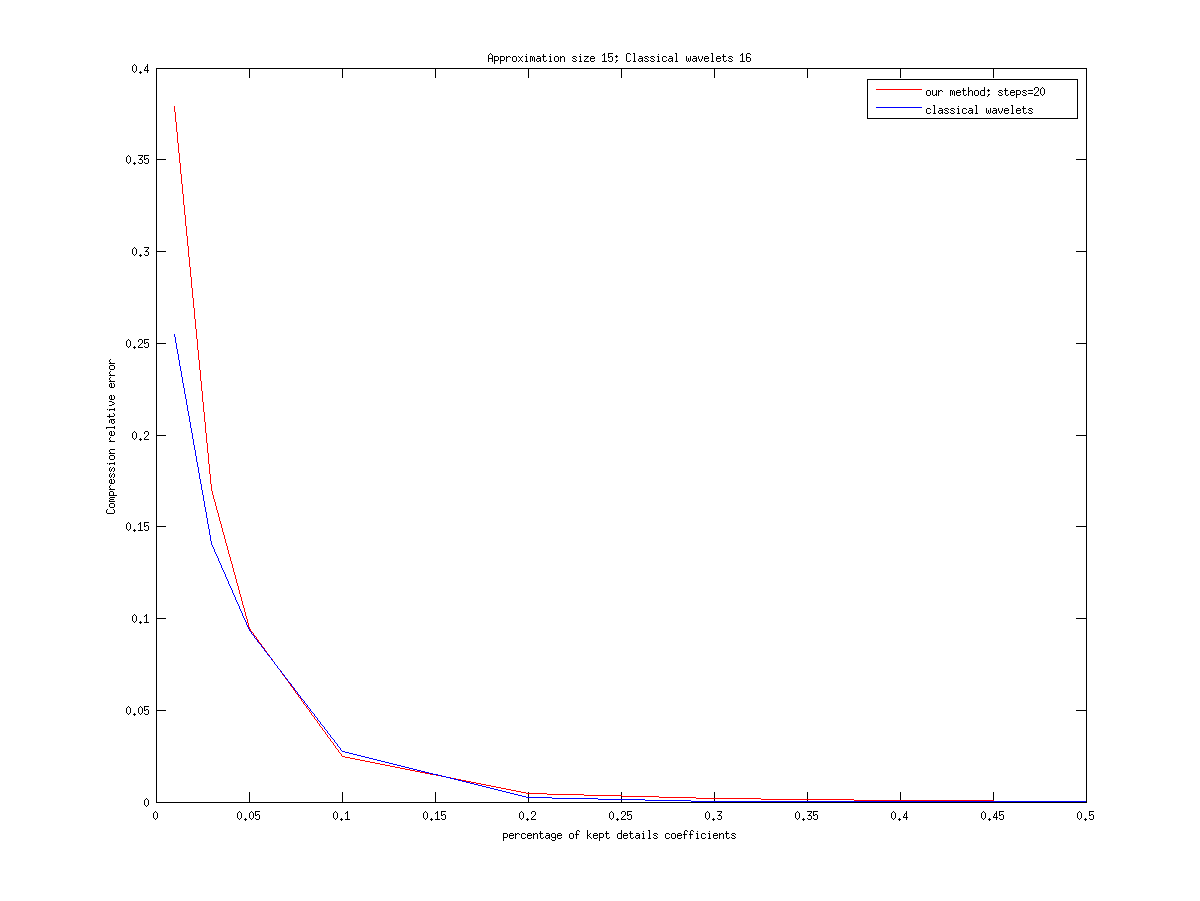}
  \caption{Relative compression error of signal in Figure~\ref{sole} in terms of the percentage of kept normalized detail coefficients. In red, the error using intertwining wavelets. In blue, error using Daubechies12 wavelets. }
  \label{domani}
\end{figure}

\begin{figure}
  \includegraphics[width=7cm]{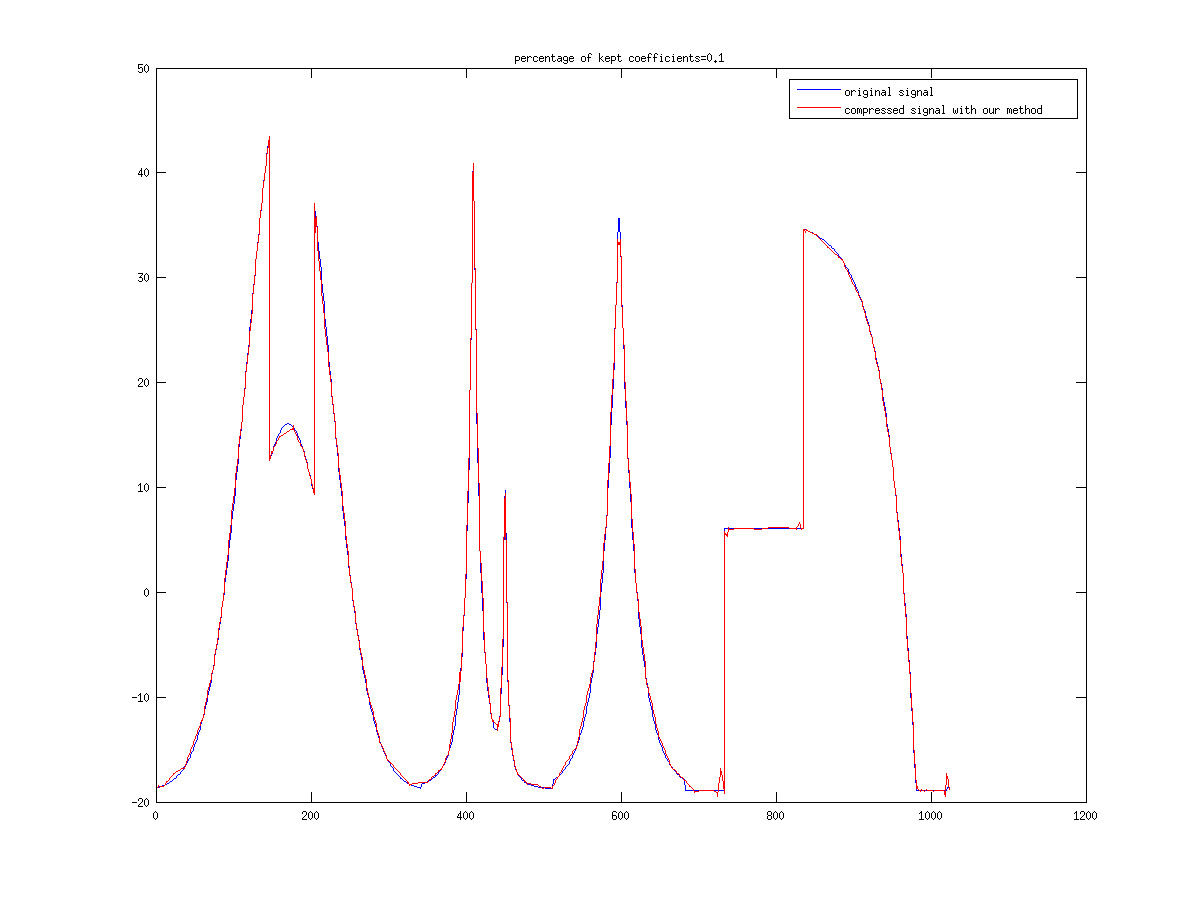}
  \caption{Original signal (blue line) and compressed one (red line) keeping 10\% of normalized detail coefficients.  }
  \label{luna}
\end{figure}

In Figure~\ref{alba}, we show two steps of the multiresolution analysis
for intertwining wavelets and Daubechies-12 wavelets applied to the signal
shown on Figure~\ref{sole}.
We refer to~\cite{ACGM2} for the description of the natural compression algorithm
that is associated with our multiresolution scheme.
Figure~\ref{domani} shows the relative compression error in terms of the percentage
of kept detail coefficients.
In Figure~\ref{luna}, we compare the original signal with
the compressed one we obtained by keeping 10\% of the detail coefficients.

\begin{figure}
    \vbox{%
        \hbox to \hsize{%
            \includegraphics[width=4cm]{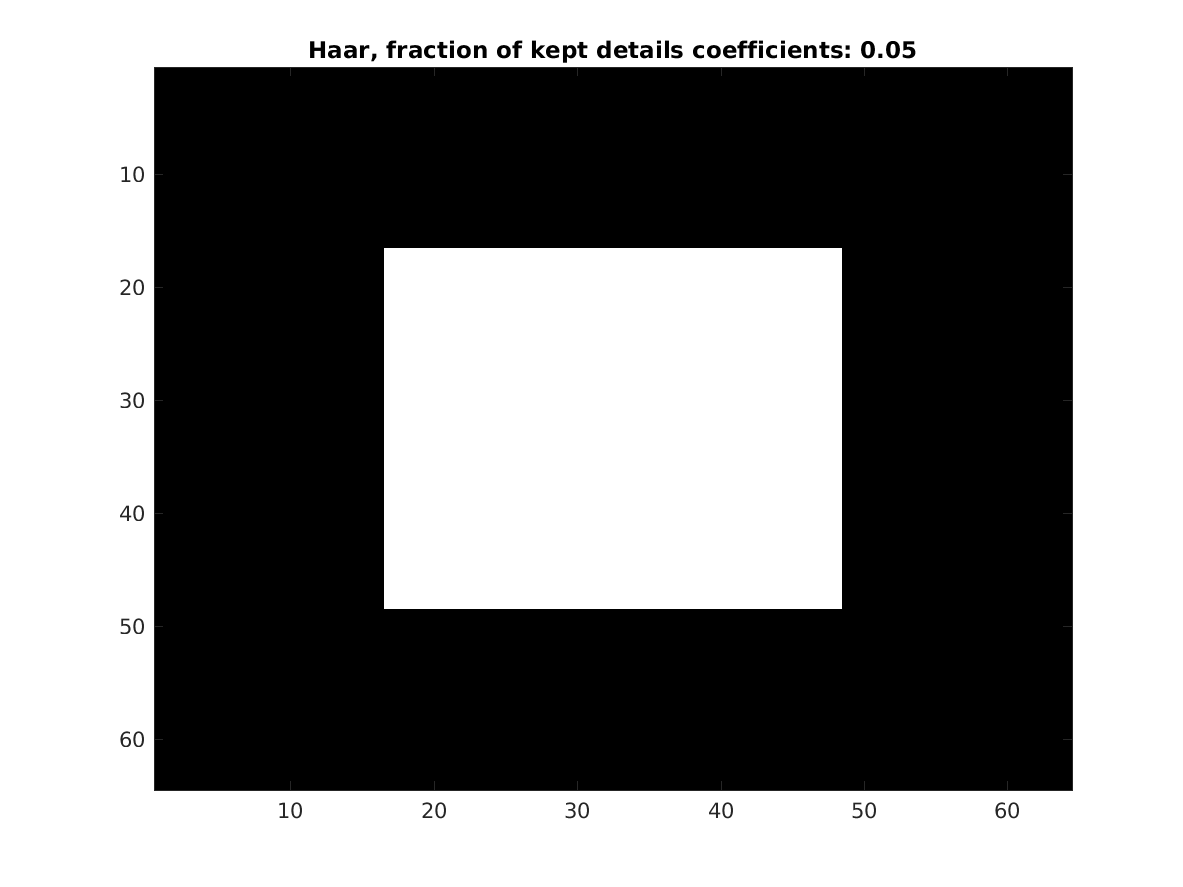}
            \hfill
            \includegraphics[width=4cm]{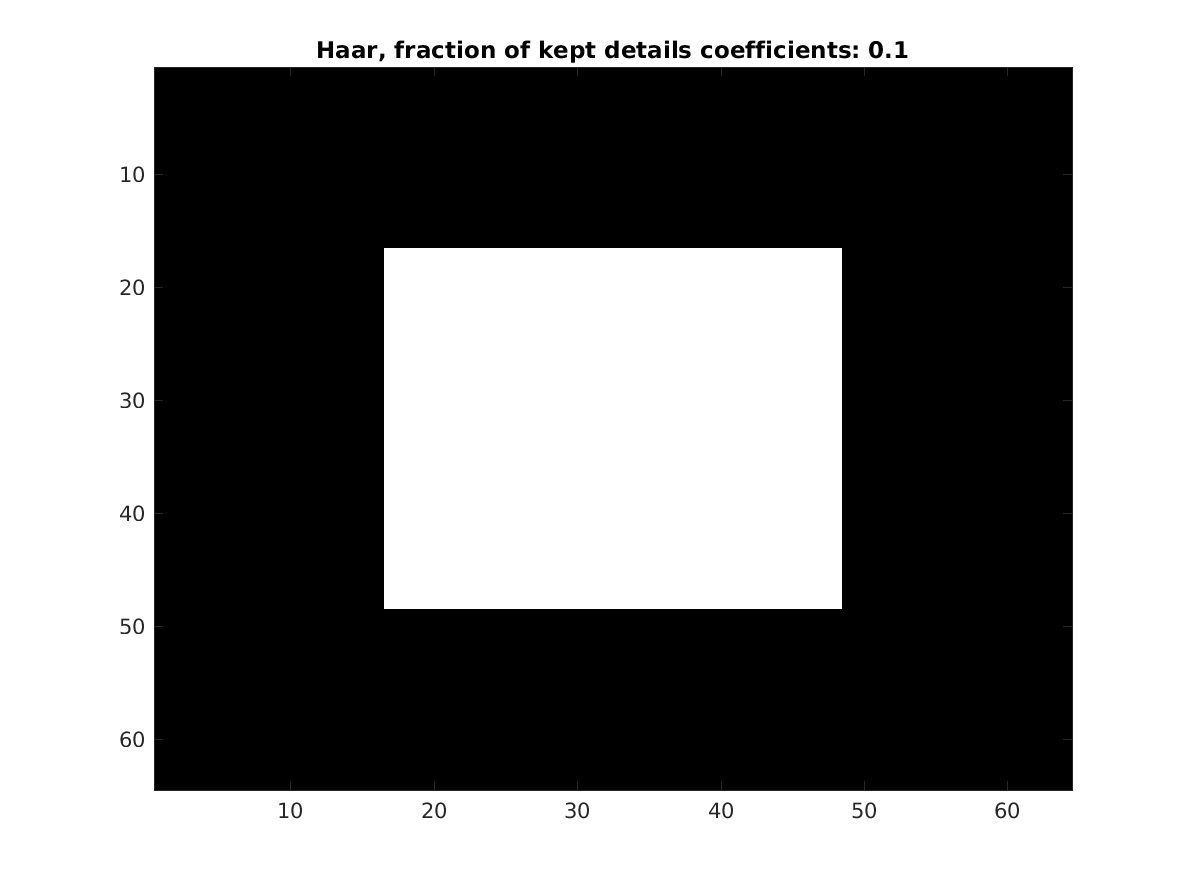}
            \hfill
            \includegraphics[width=4cm]{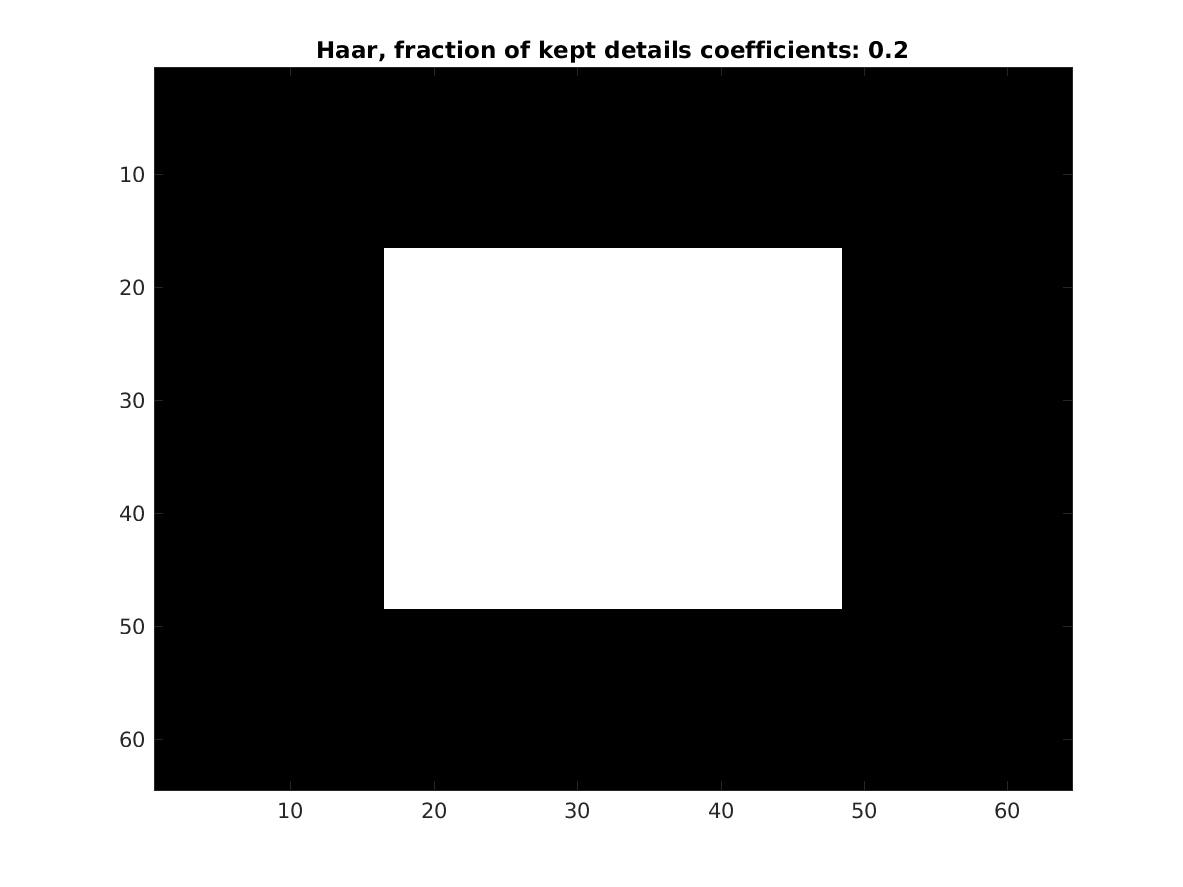}
            \hfill
            \includegraphics[width=4cm]{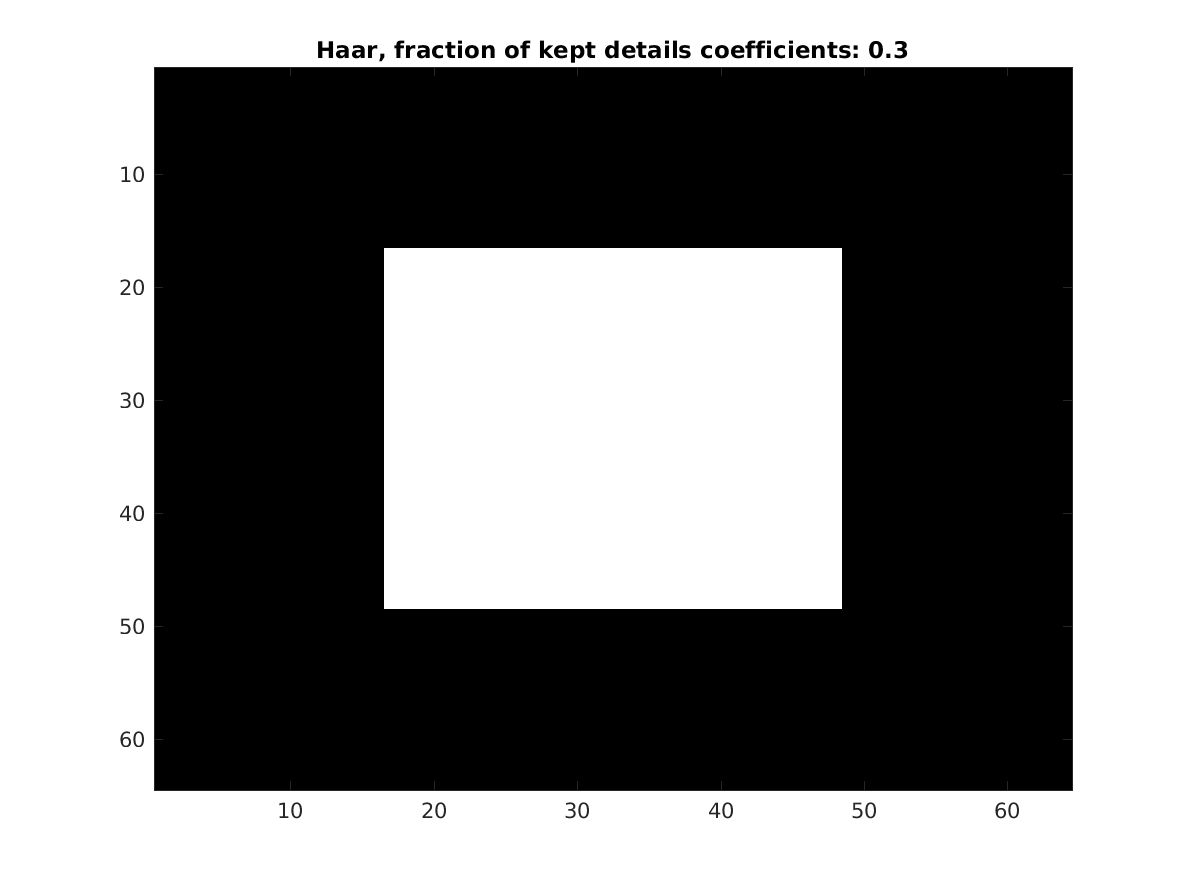}
        }           
        \medskip
        \hbox to \hsize{%
            \includegraphics[width=4cm]{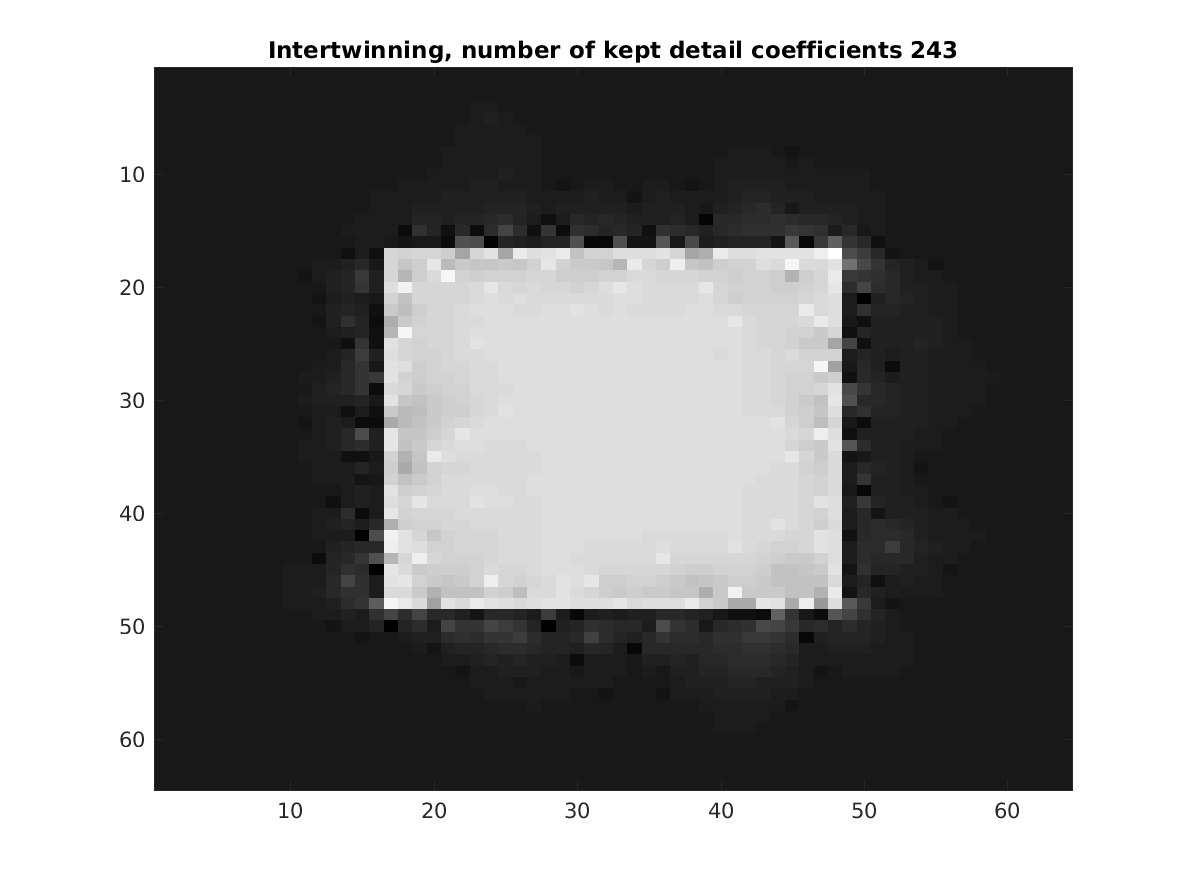}
            \hfill
            \includegraphics[width=4cm]{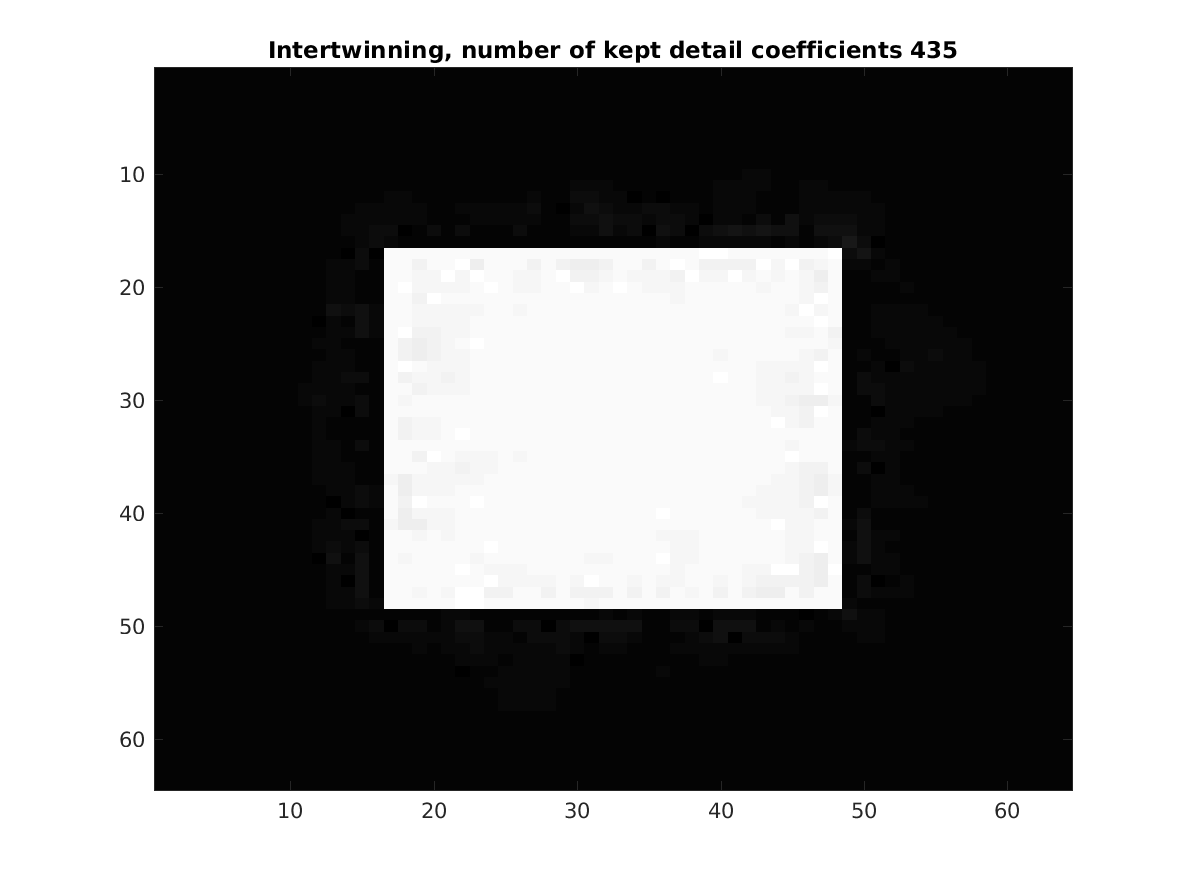}
            \hfill
            \includegraphics[width=4cm]{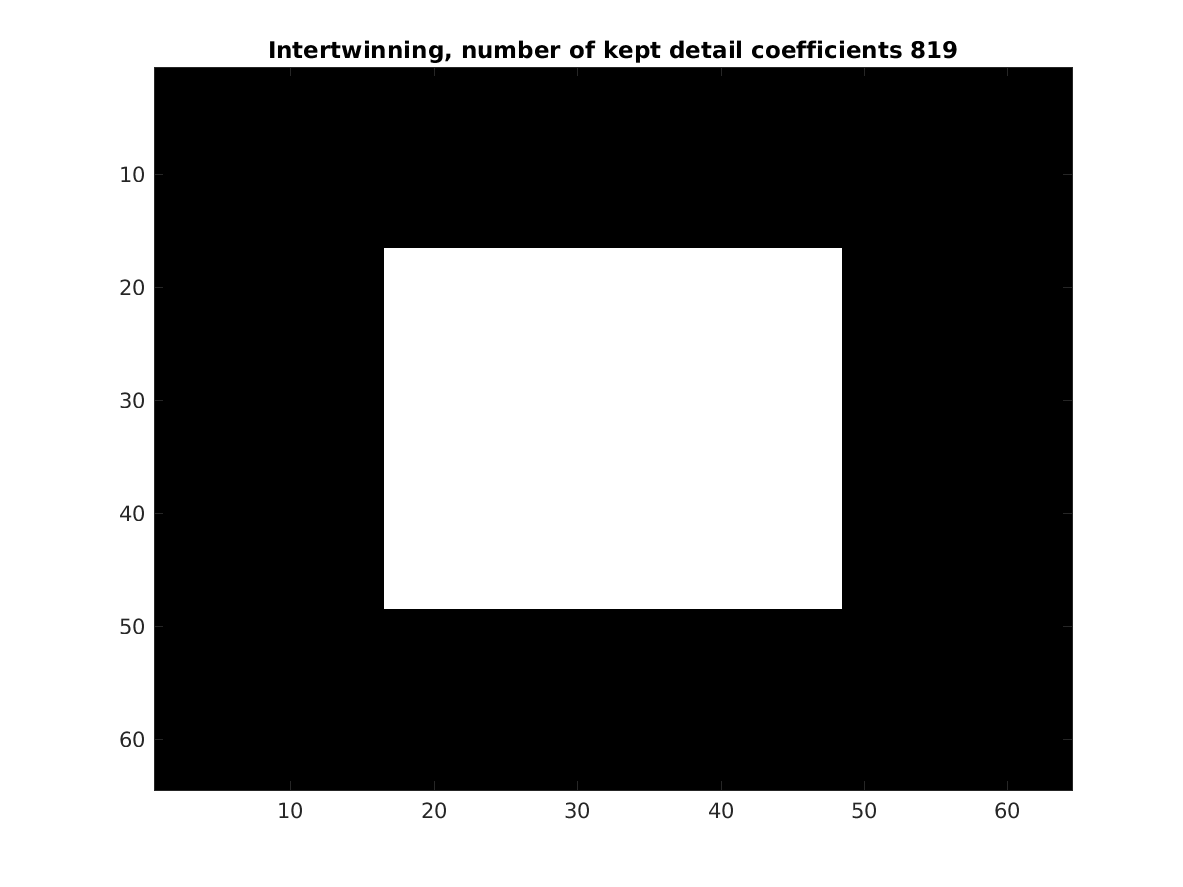}
            \hfill
            \includegraphics[width=4cm]{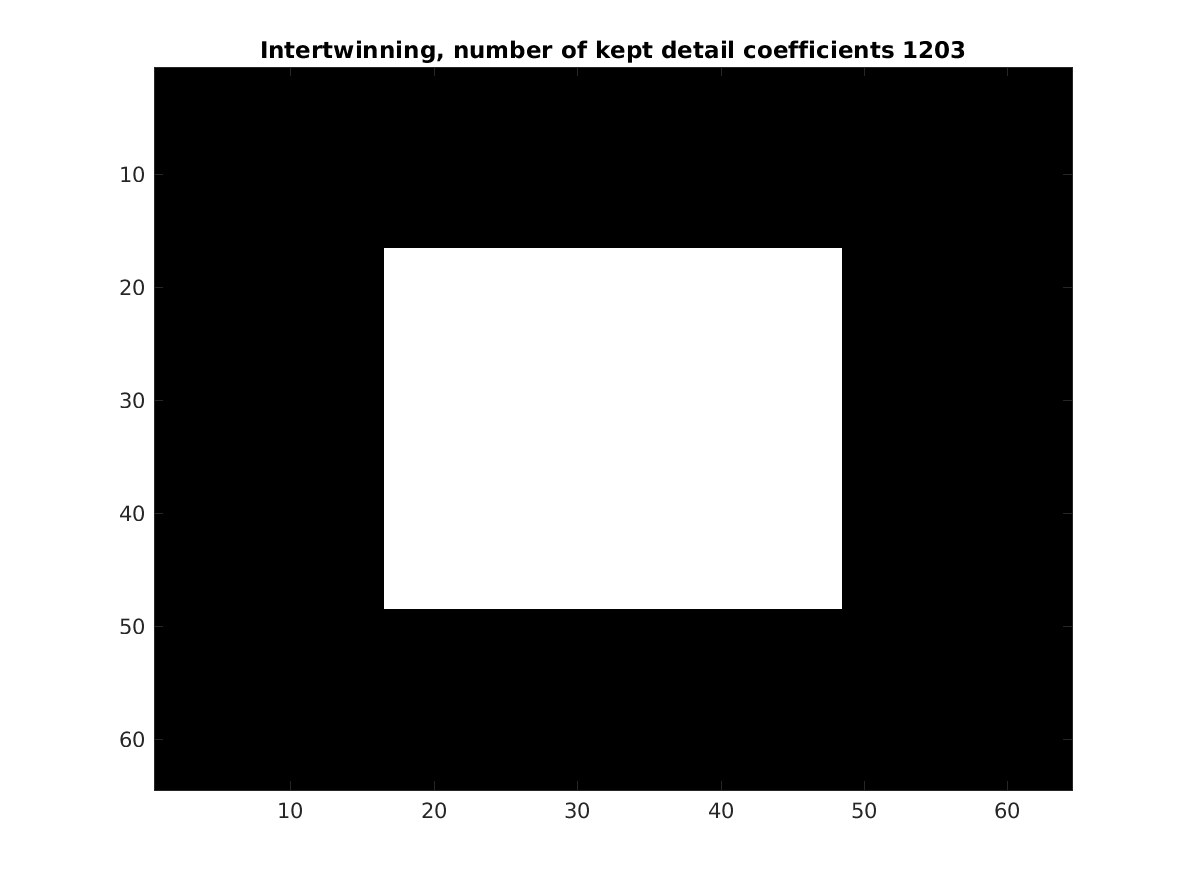}
        }           
    }
    \caption{
        A compressed white rectangle with Haar wavelets on the first line, for which the percentage
        of kept detail coefficients is reported in each picture.
        The same rectangle is compressed with intertwining wavelets on the second line,
        with the same number of kept coefficients.
    }
    \label{nuvola}
\end{figure}

\begin{figure}
    \includegraphics[width=4cm]{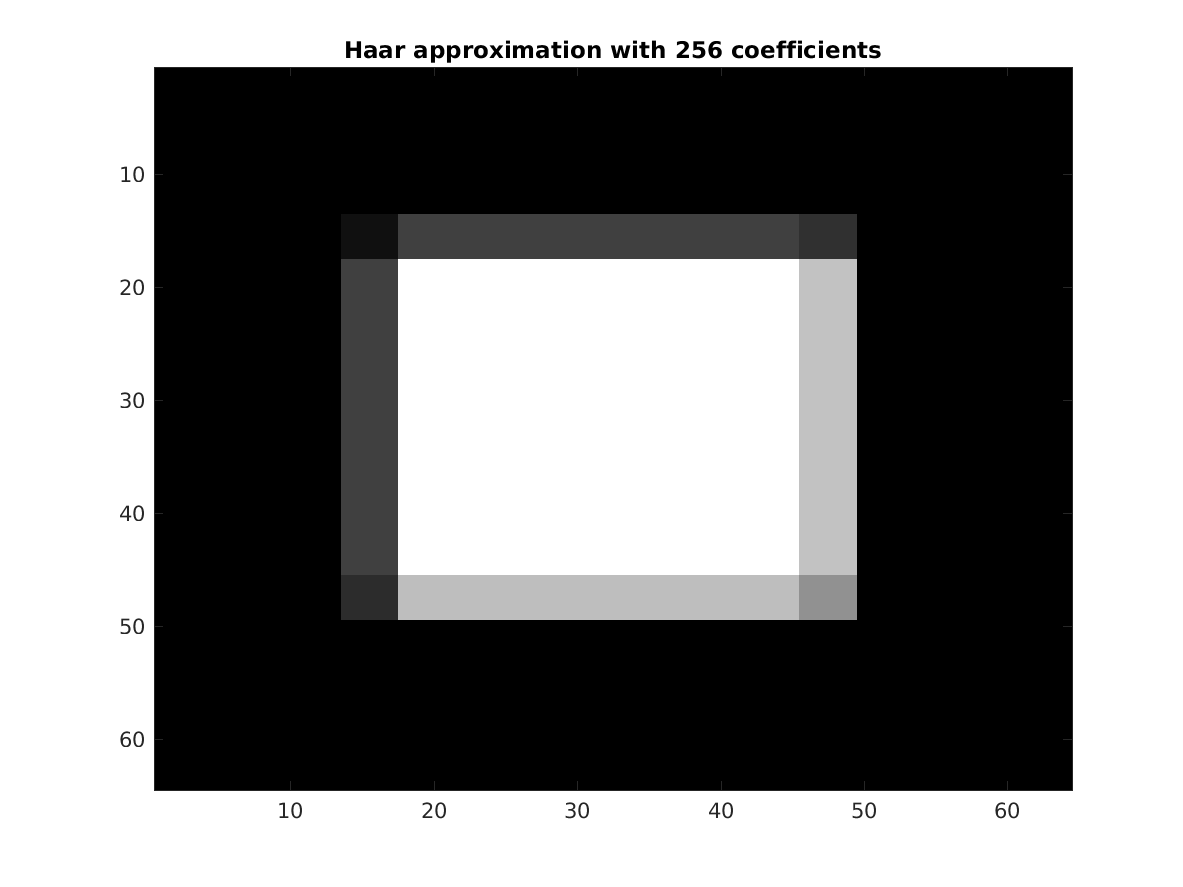}
    \includegraphics[width=4cm]{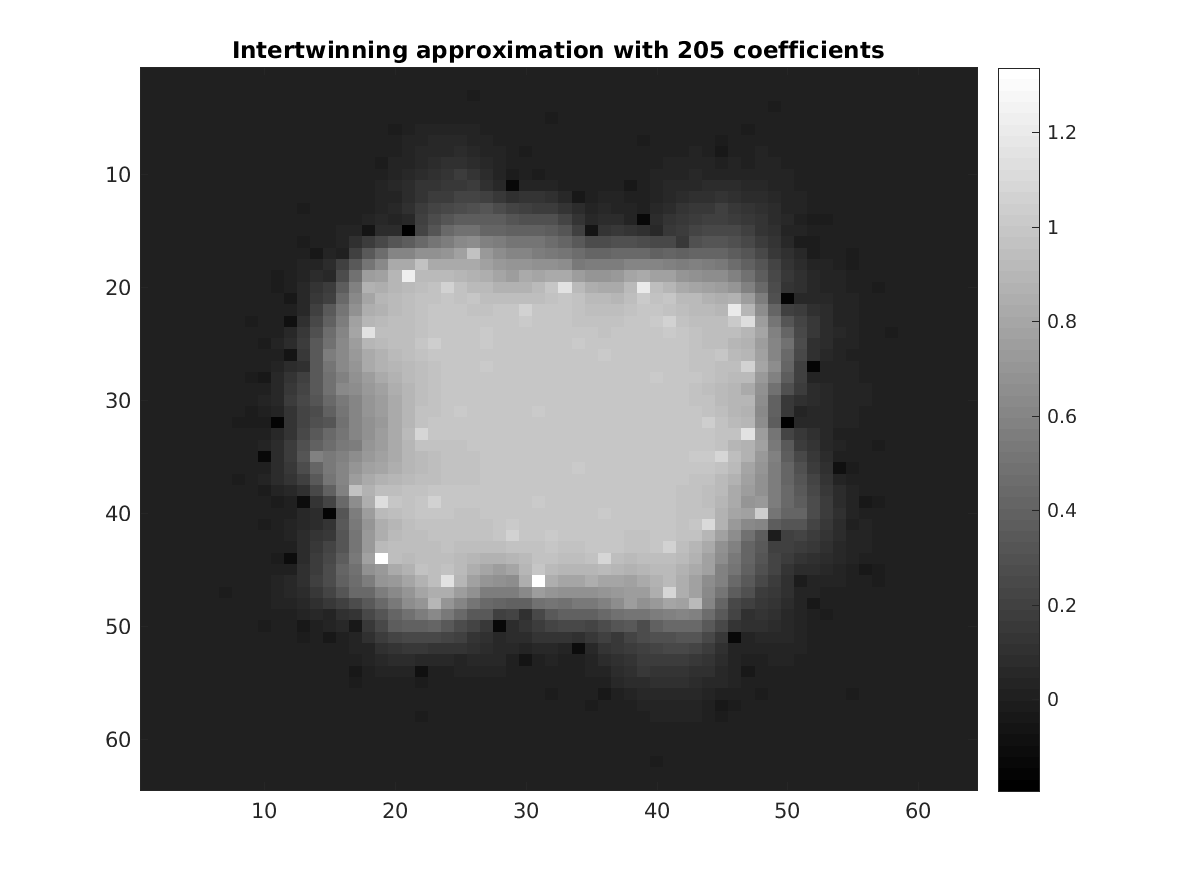}
    \caption{
        Approximations for Haar and intertwining wavelets,
        with 256 and 205 coefficients respectively.
    }\label{verde}
\end{figure}

\begin{figure}
    \hbox to \hsize{\hss
        \vbox{%
            \hbox to 11.5cm{%
                \includegraphics[width=5.2cm]{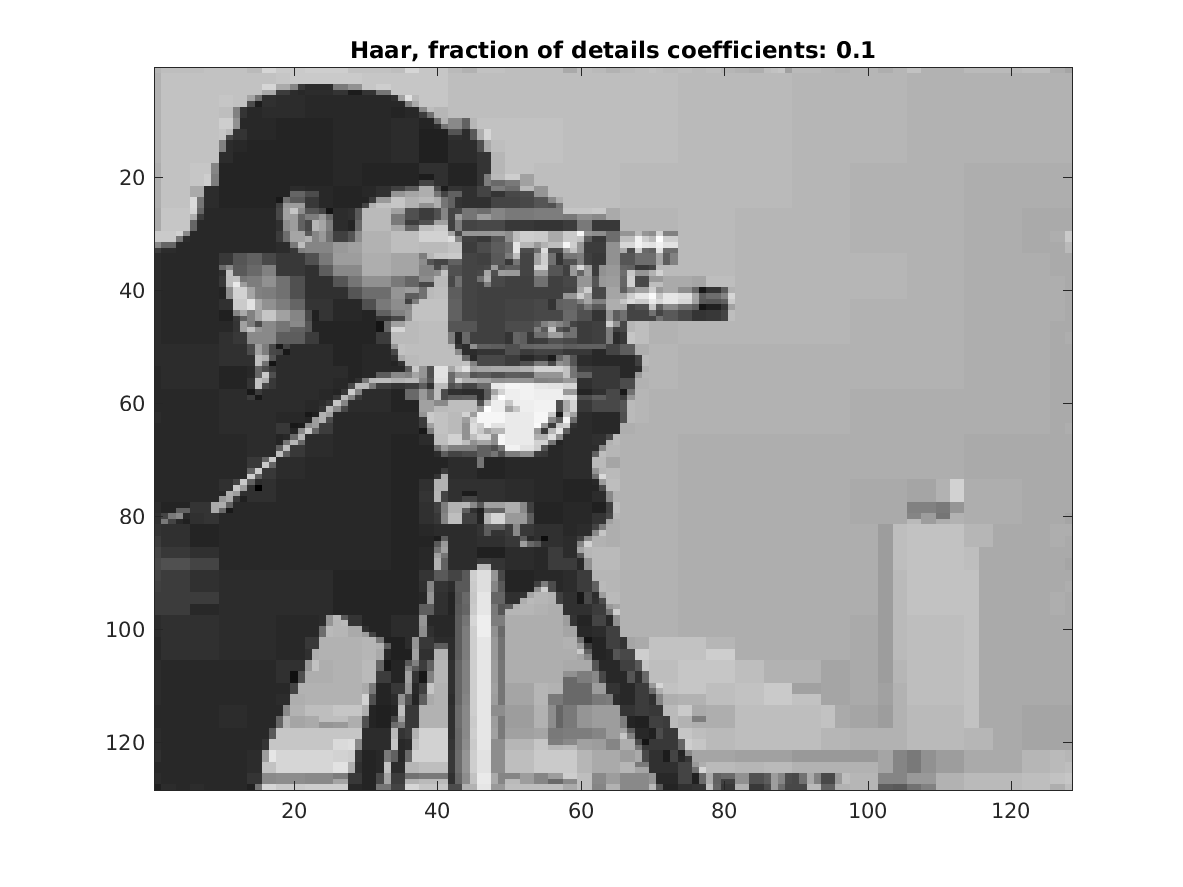}
                \hfill
                \includegraphics[width=5.2cm]{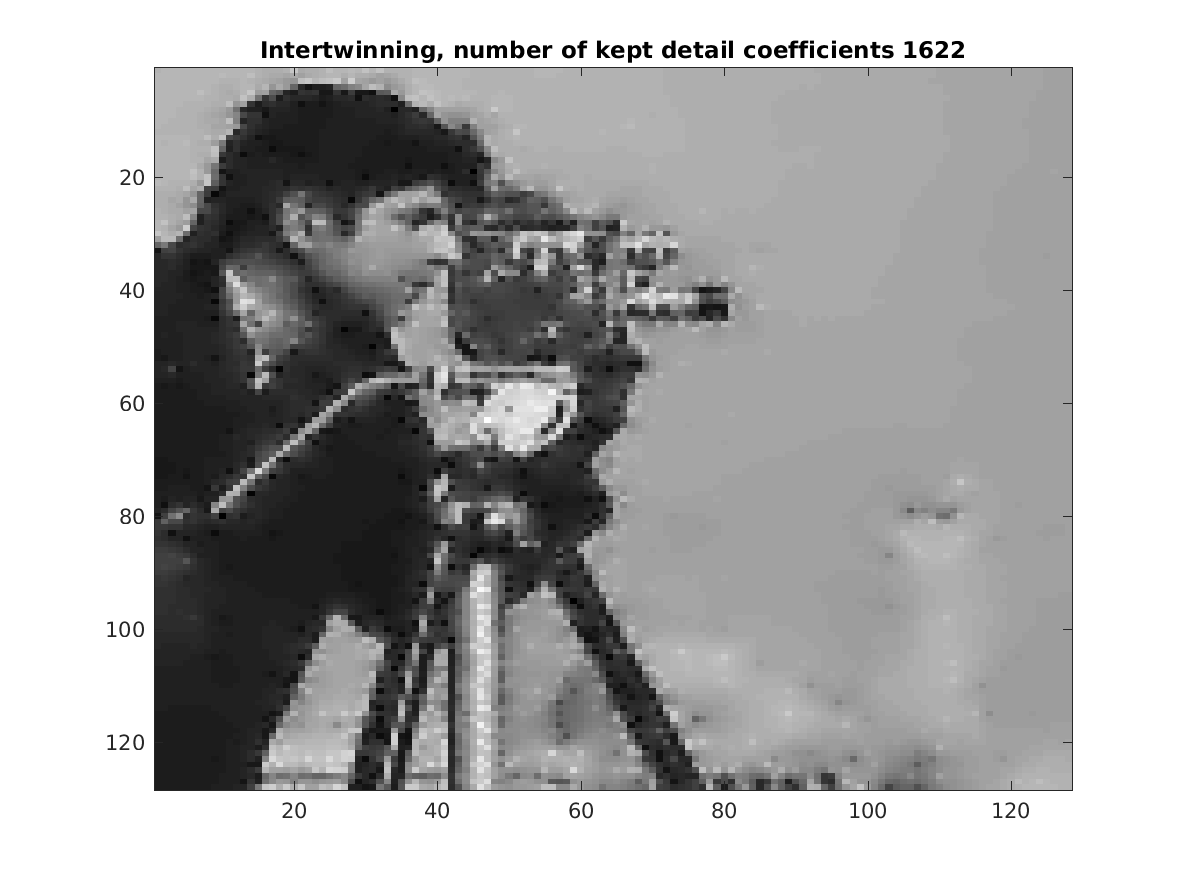}
            }
            \medskip
            \hbox to 11.5cm{%
                \includegraphics[width=5.2cm]{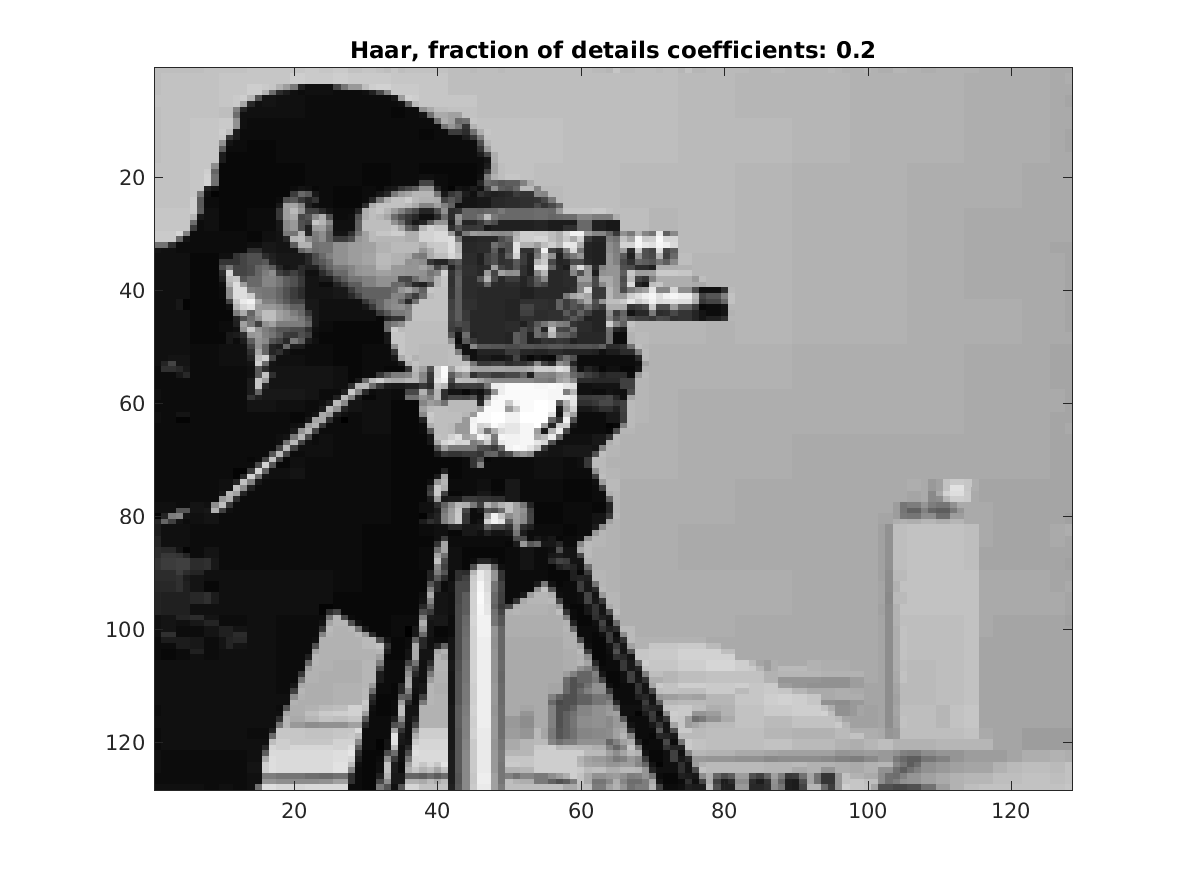}
                \hfill
                \includegraphics[width=5.2cm]{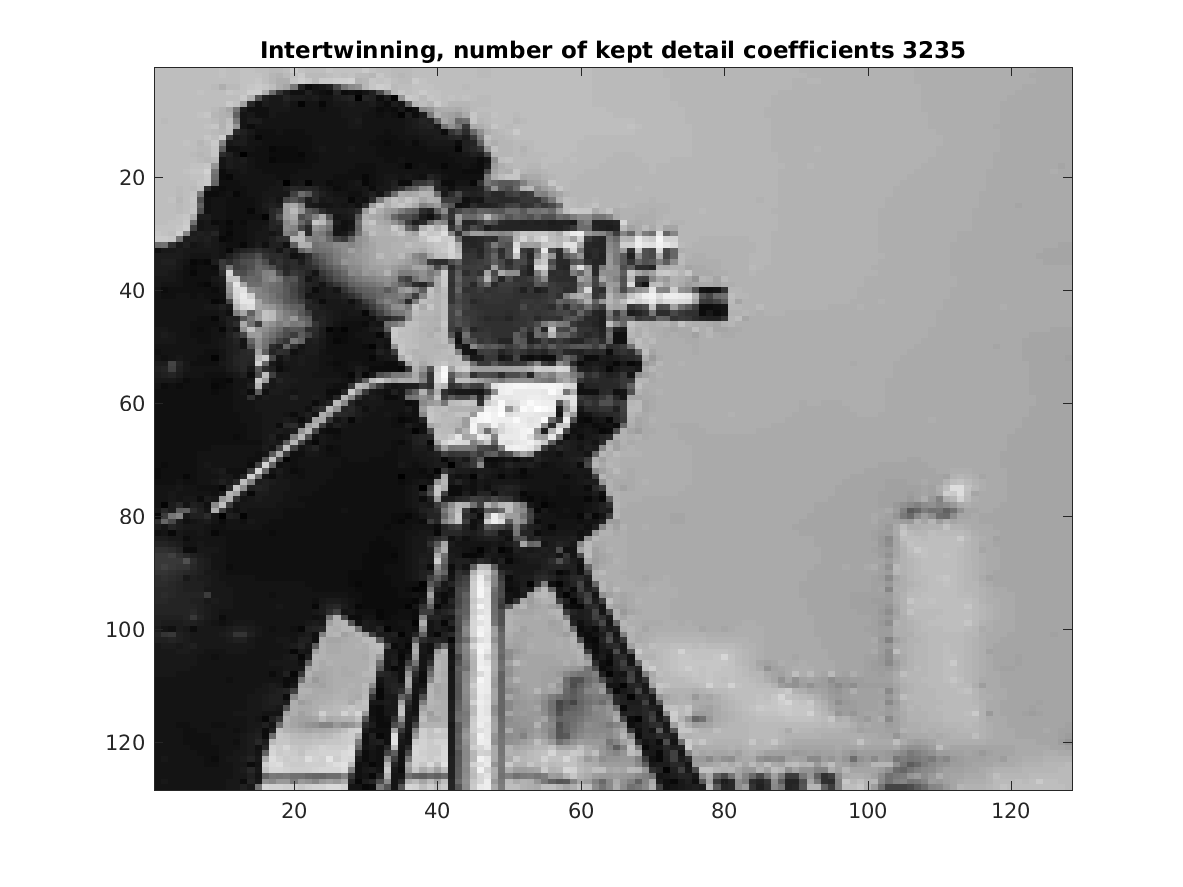}
            }
            \medskip
            \hbox to 11.5cm{%
                \includegraphics[width=5.2cm]{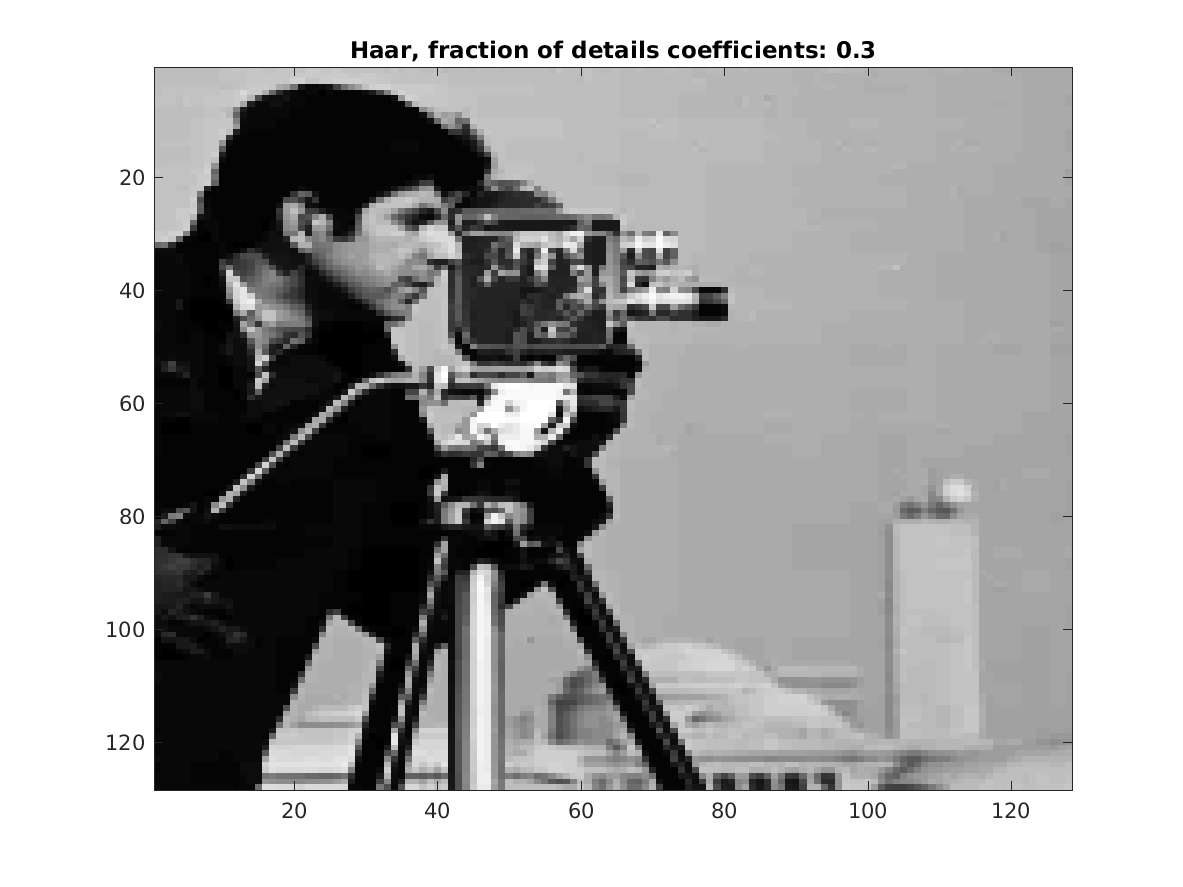}
                \hfill
                \includegraphics[width=5.2cm]{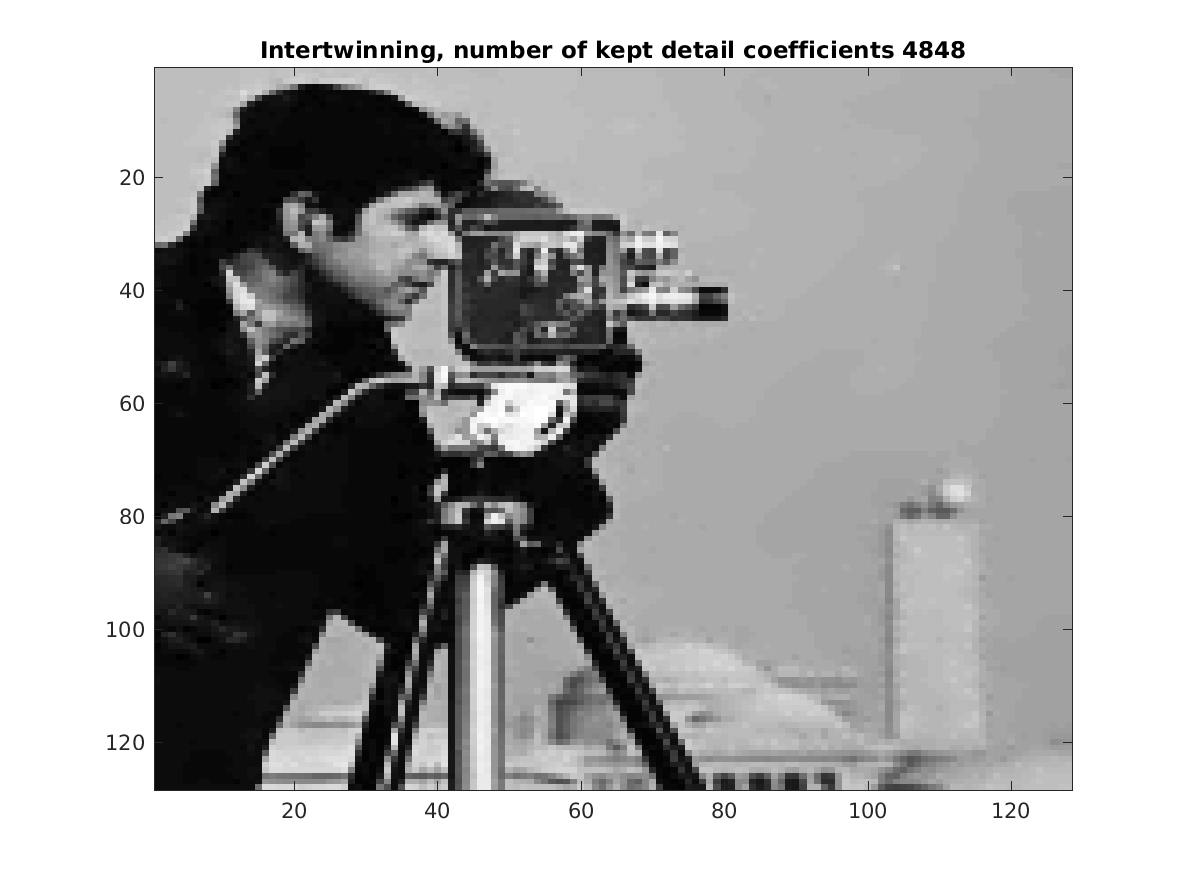}
            }
            \medskip
            \hbox to 11.5cm{%
                \includegraphics[width=5.2cm]{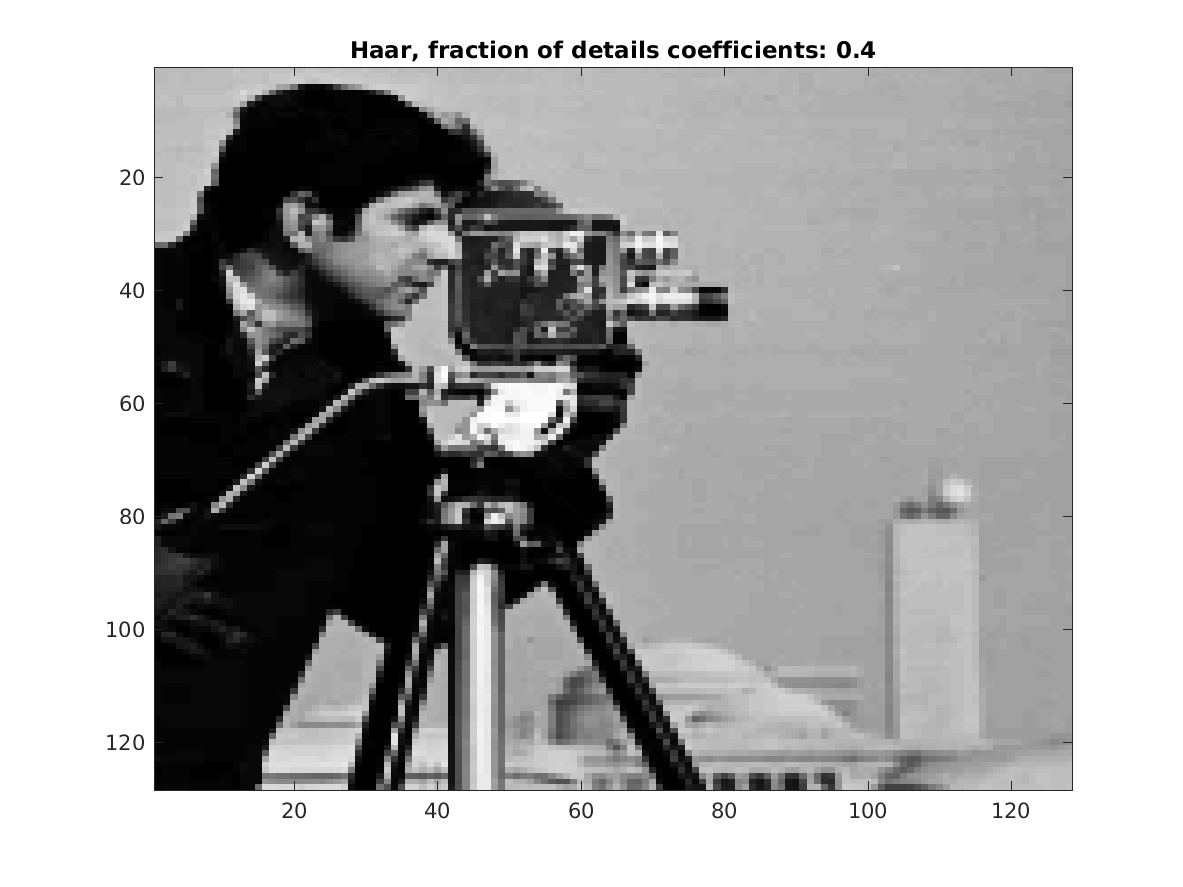}
                \hfill
                \includegraphics[width=5.2cm]{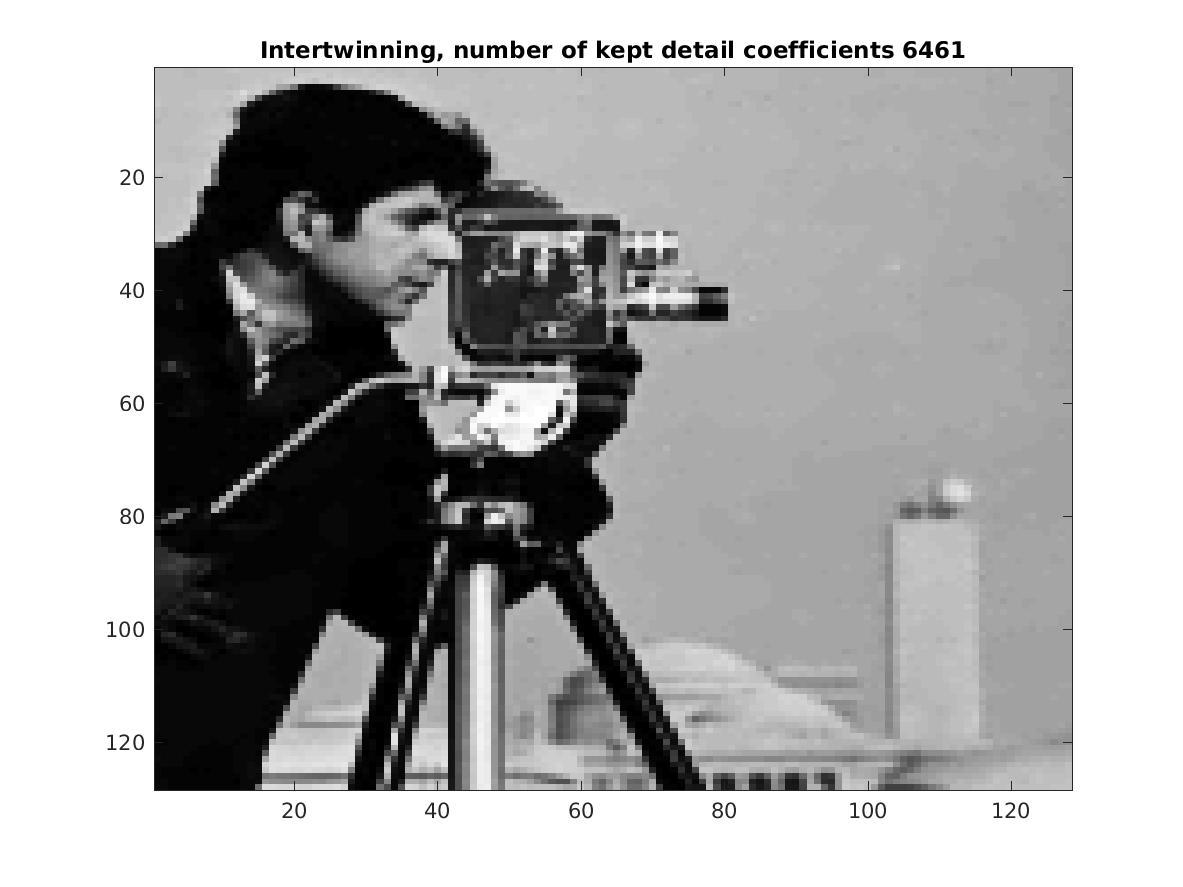}
            }
            \medskip
            \hbox to 11.5cm{%
                \includegraphics[width=5.2cm]{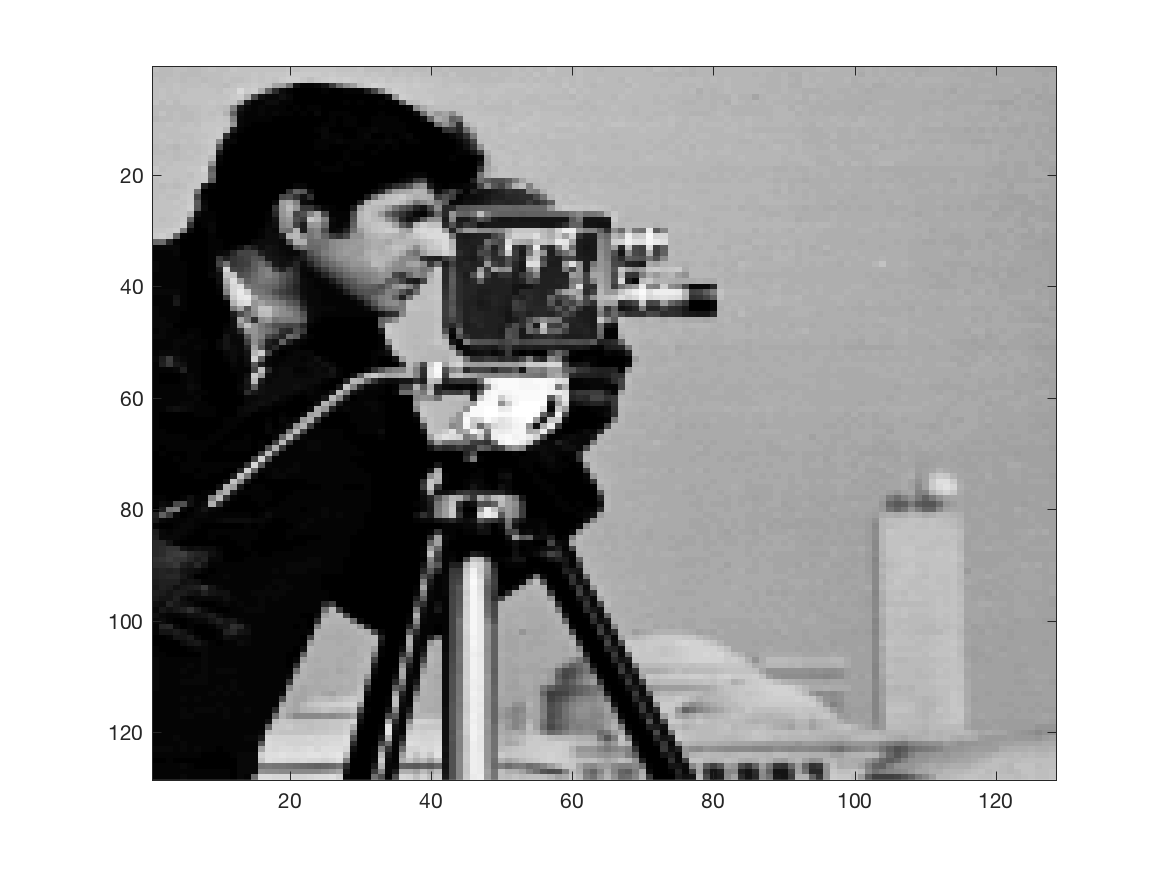}
                \hfill
                \includegraphics[width=5.2cm]{cameraman_original.png} 
            }
        }
    \hss}
    \caption{%
        Compression with Haar wavelets on the left,
        with intertwining wavelets on the right.
        In both case the last picture is the original signal.
    }\label{pioggia}
\end{figure}

\begin{figure}
    \includegraphics[width=4cm]{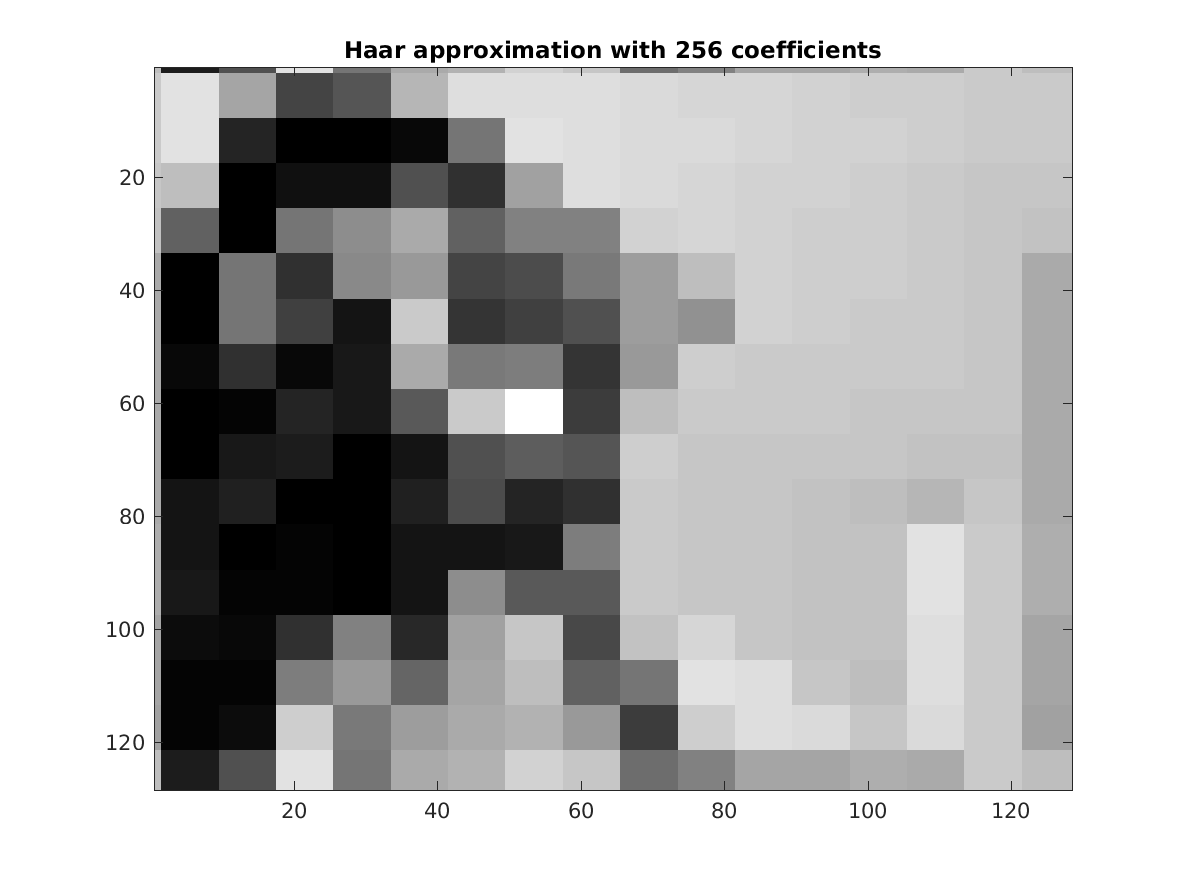}
    \includegraphics[width=4cm]{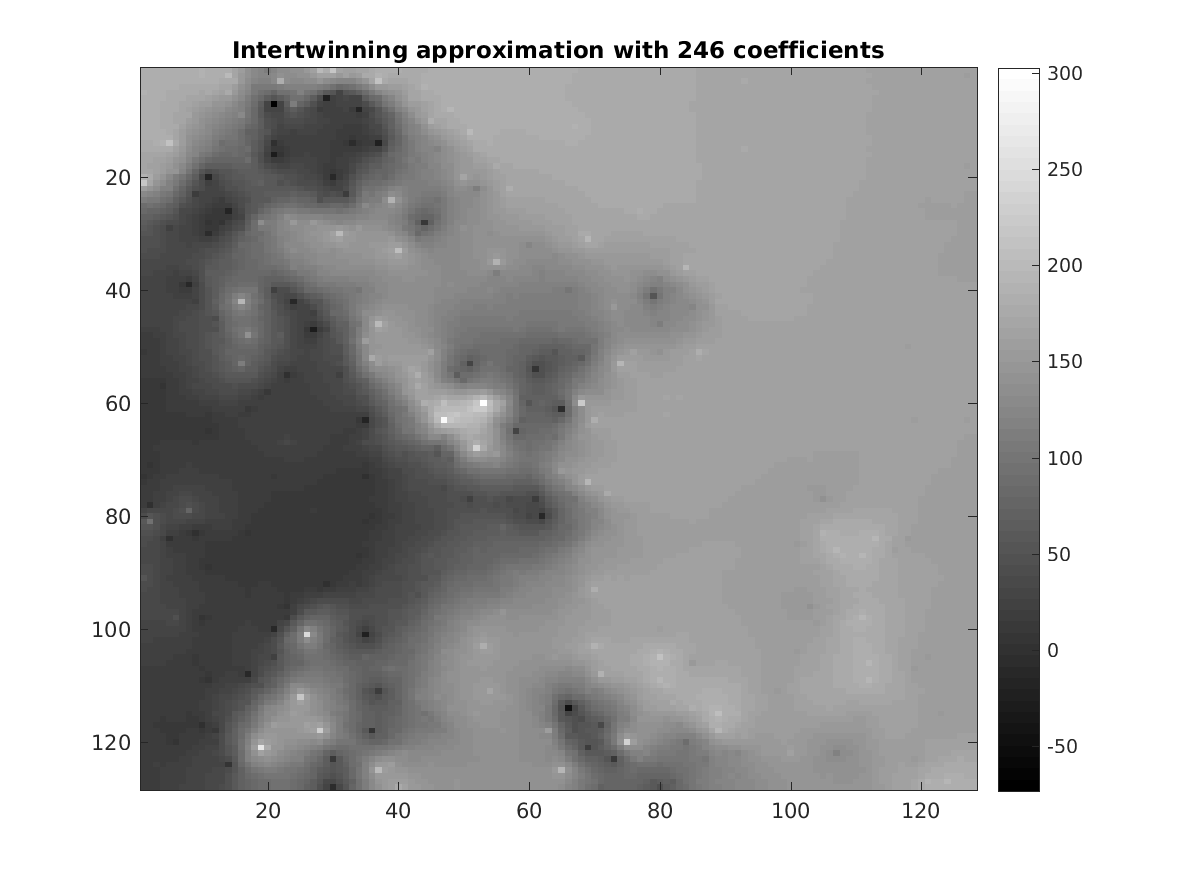}
    \caption{
        Approximations for Haar and intertwining wavelets,
        with 256 and 246 coefficient respectively.
    }\label{verde2}
\end{figure}

In Figures~\ref{nuvola} and~\ref{pioggia},
we show compressed images of a rectangle and a cameraman, respectively, obtained by keeping 
the same number of coefficients.
Figures~\ref{verde} and~\ref{verde2} show the upsampled approximations
(obtained by setting to zero all detail coefficients).

\subsubsection{A last example}
\begin{figure}
\begin{multicols}{2}
(a) Original signal.

\centerline{\includegraphics[width=6cm,height=4cm]{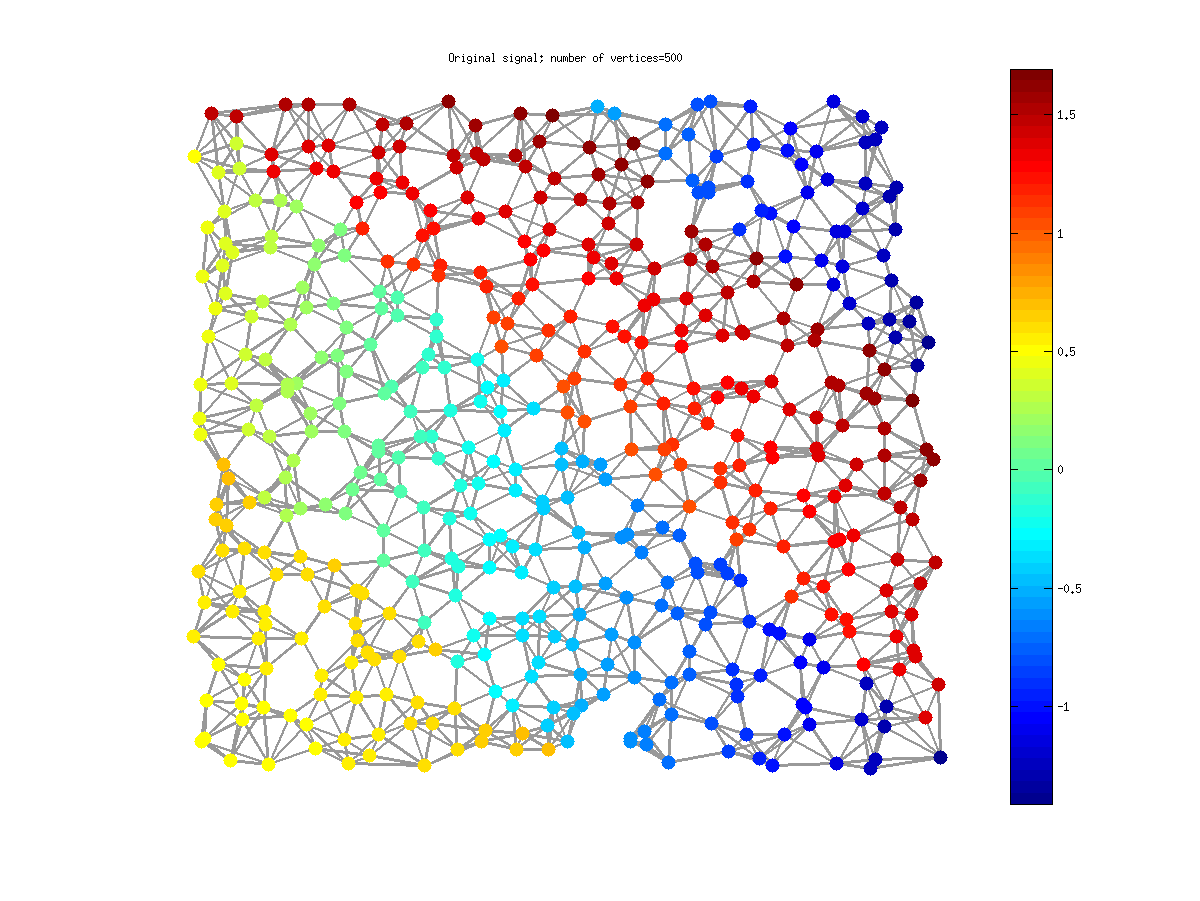}}

(b) Approximation. Size of $f_2$: 206.

\centerline{\includegraphics[width=6cm,height=4cm]{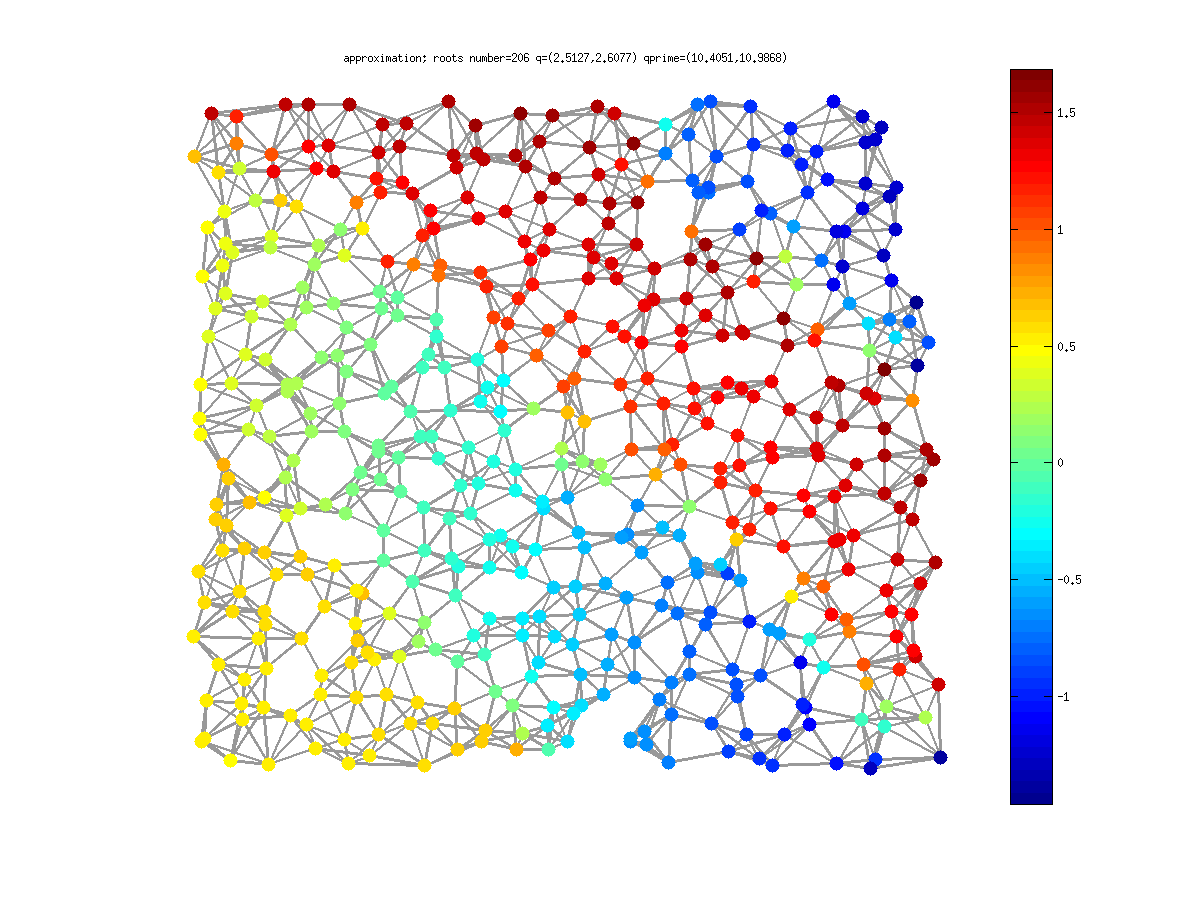}}
\columnbreak

(c) Detail at scale 1.   Size of $g_1$: 188.

\centerline{\includegraphics[width=6cm,height=4cm]{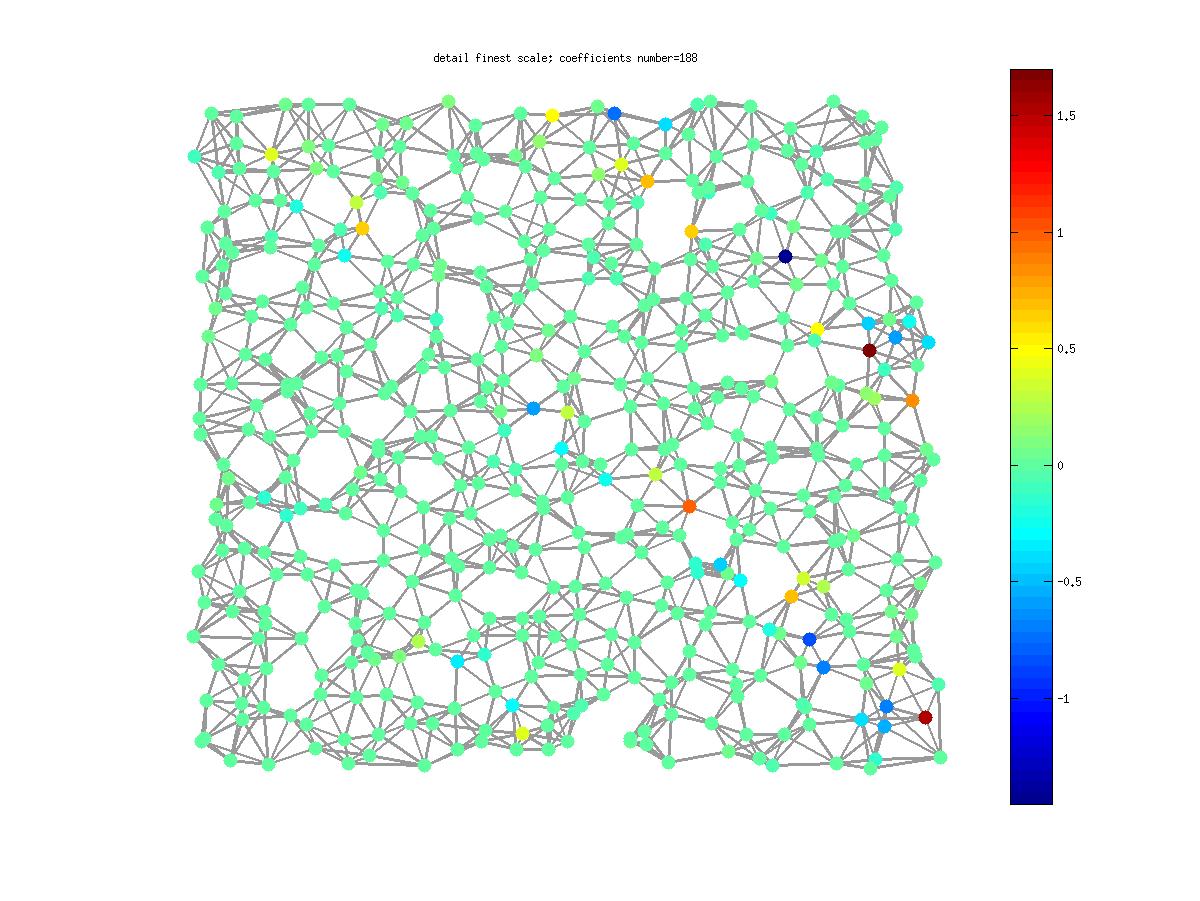}}

(d) Detail at scale 2.  Size of $g_2$: 106.

\centerline{\includegraphics[width=6cm,height=4cm]{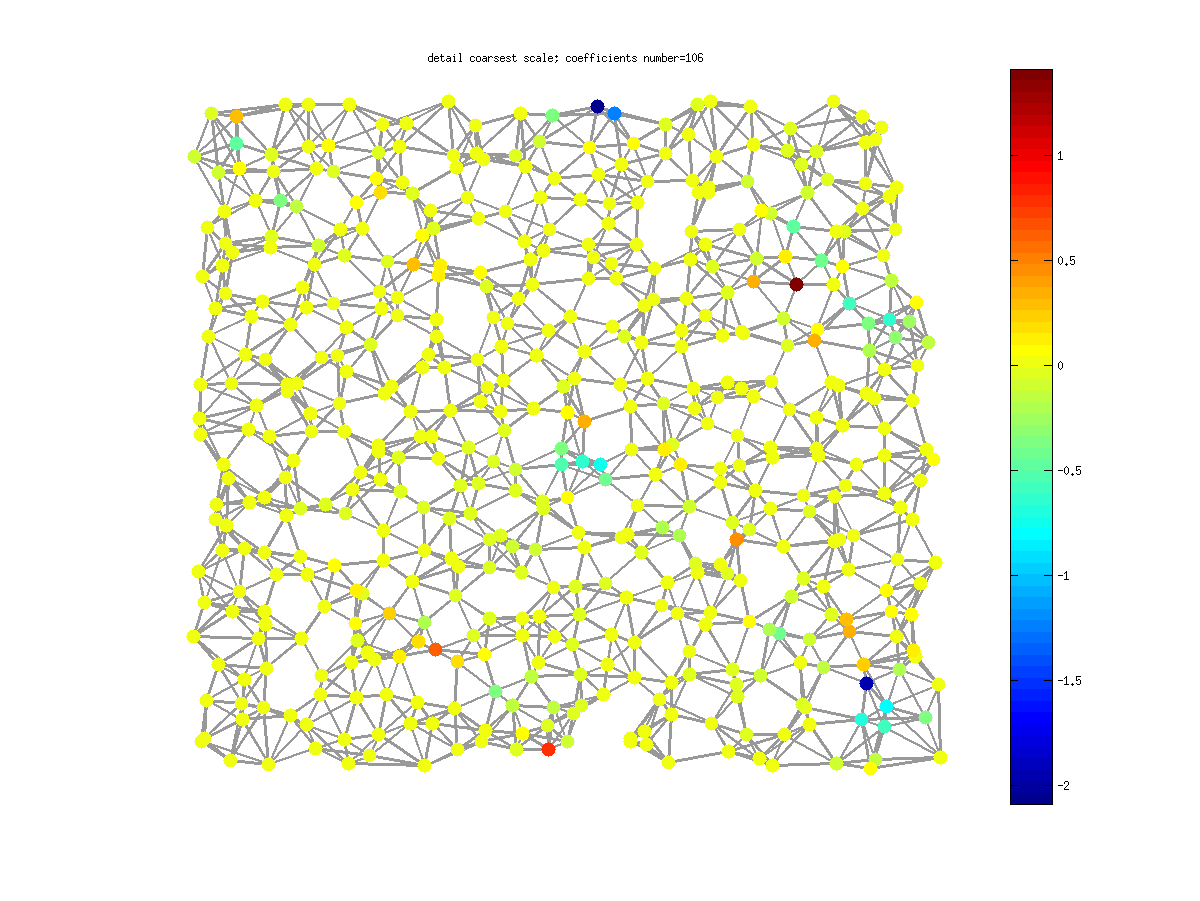}}
\end{multicols}
\caption{Two steps of the intertwining wavelets multiresolution upsampled to the original graph: (a) $f_0$; (b) $\bar{R}_0 \bar{R}_1 f_2$; 
(c) $ \brR_0 g_1$;  (d) $\bar{R}_0 \brR_1 g_2$}.
\label{lunettes}
\end{figure}

\begin{figure}
\begin{tabular}{cc}
  \includegraphics[scale=0.35]{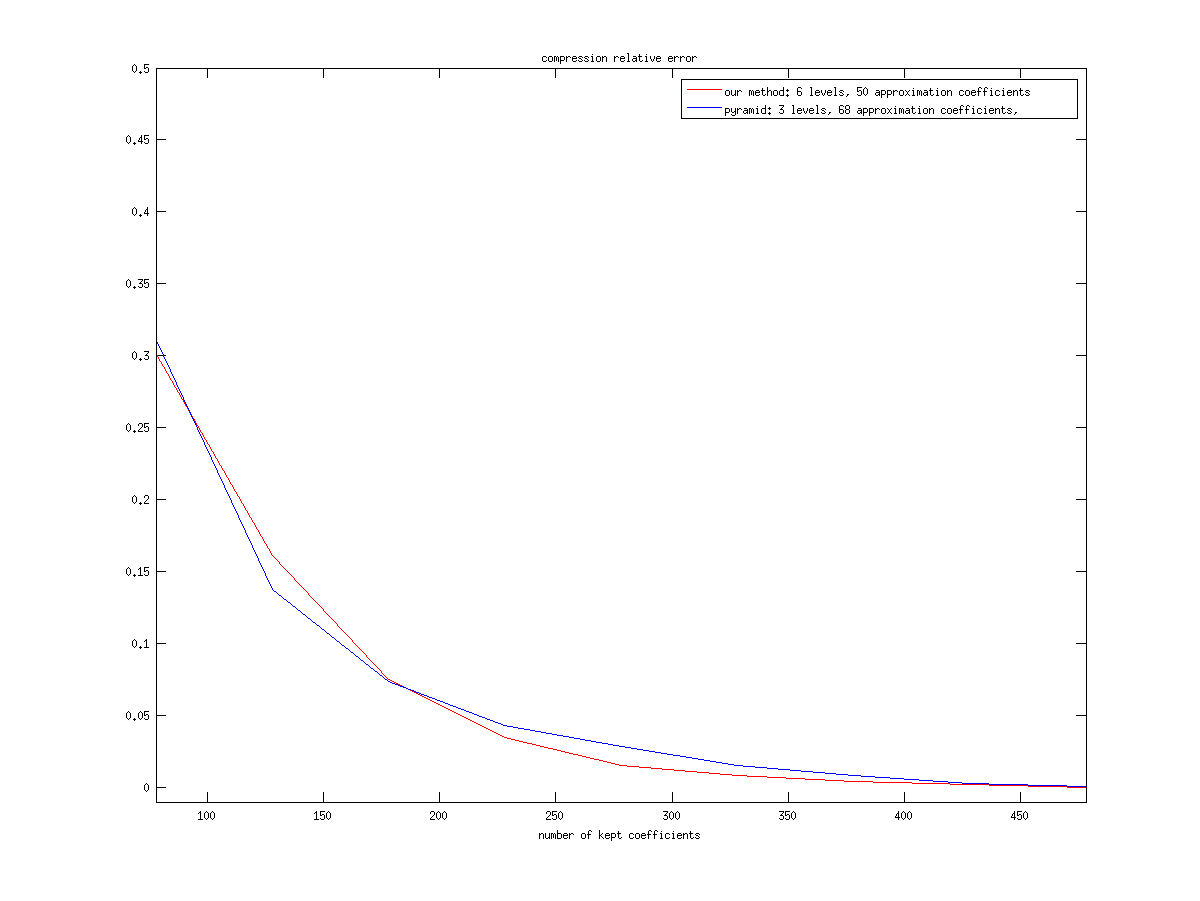}
 &
  \includegraphics[scale=0.35]{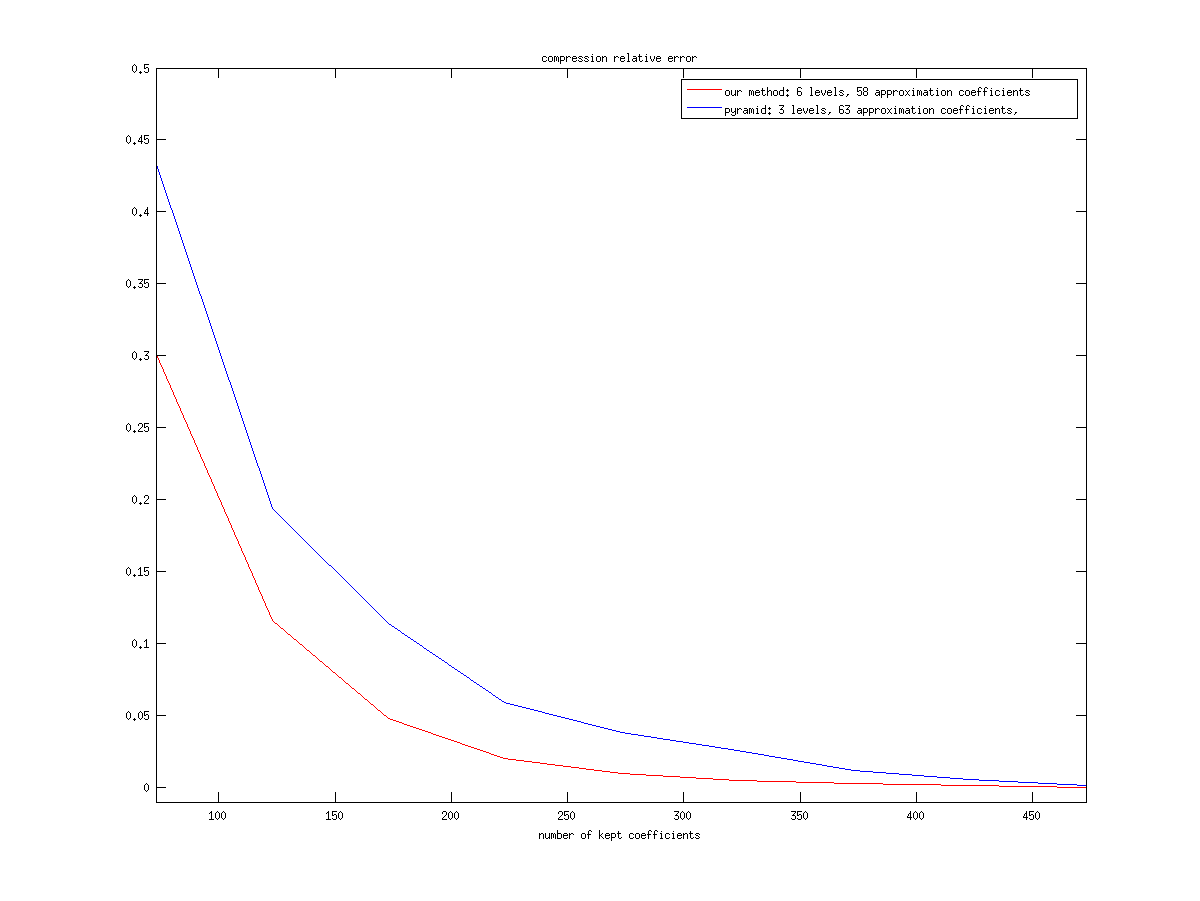}
\\
Without sparsification 
&
With sparsification 
\end{tabular}
 \caption{Relative compression error of signal in Figure~\ref{lunettes}(a), in terms of the number of kept coefficients. In red, the error using intertwining wavelets. In blue, error using the spectral graph wavelets pyramidal algorithm. In the first graph, we do not sparsify the graph, while we perform sparsification on the second. }
 \label{tortora}
 \end{figure}

We took from~\cite{SFV} and the GSP toolbox~\cite{GSP}
the sensor graph and the signal represented in 
Figure~\ref{lunettes} together with two steps
of the multiscale analysis.
In Figure~\ref{tortora}
we compare the results of the intertwining compression algorithm
with those of the spectral graph wavelet pyramidal algorithm.
Unless specified, in all the previous experiments we always included
our sparsification procedure.
In this last figure we added the result of the algorithm 
without sparsification.
It is worth to note that, in this example at least,
sparsification helps,
not only for algorithmic complexity.

\end{document}